%% file: paper_informativity.tex
\documentclass{IEEE-Con-Sys-mag}

\usepackage{amsmath}
\usepackage{amsfonts}
\usepackage{amssymb}
\usepackage{graphicx}   
\usepackage{color}
\usepackage{comment}
\usepackage{graphicx}
\usepackage{cite}
\usepackage{float}
\usepackage{pmat}
\input{ka-newcommands}

\usepackage{tikz}
\usetikzlibrary{tikzmark,patterns}
\usepackage{mathtime-2020}
\usepackage{pgf}
\usepackage{pgfplots}
\usepackage{subfigure}

\usepackage{bm}
\newcommand{\bmf}{\bm{f}}
\newcommand{\bmu}{\bm{u}}
\newcommand{\bmx}{\bm{x}}
\newcommand{\bmv}{\bm{v}}
\newcommand{\bmw}{\bm{w}}
\newcommand{\bmy}{\bm{y}}
\newcommand{\bmd}{\bm{d}}
\newcommand{\bmr}{\bm{r}}
\newcommand{\bmz}{\bm{z}}

\newcommand{\xm}{X_-}
\newcommand{\xp}{X_+}
\newcommand{\um}{U_-}

\newcommand{\xmt}{X_-^\top}
\newcommand{\xpt}{X_+^\top}
\newcommand{\umt}{U_-^\top}

\newcommand{\bpi}{\bm{\Pi}}
\newcommand{\schur}{\!\mid\!}
\newcommand{\mpi}{\bbm \Pi_{11}&\Pi_{12}\\\Pi_{21}&\Pi_{22}\ebm}
\newcommand{\gi}{^\dagger}
\renewcommand{\S}[1]{\mathbb{S}^{#1}}

\DeclareMathOperator{\In}{In}

\definecolor{marine}{rgb}{0.149,0.149,0.804}
\definecolor{fbrick}{rgb}{0.804,0.149,0.149}
\definecolor{apricot}{rgb}{0.98,0.81,0.69}
\definecolor{Blue}{cmyk}{1.,1.,0,0}
\definecolor{amber}{rgb}{1.0, 0.75, 0.0}
\definecolor{applegreen}{rgb}{0.55, 0.71, 0.0}
\definecolor{bananayellow}{rgb}{1.0, 0.88, 0.21}
\definecolor{bittersweet}{rgb}{1.0, 0.44, 0.37}
\definecolor{candyapplered}{rgb}{1.0, 0.03, 0.0}
\definecolor{carminered}{rgb}{1.0, 0.0, 0.22}
\definecolor{carrotorange}{rgb}{0.93, 0.57, 0.13}
\definecolor{goldenyellow}{rgb}{1.0, 0.87, 0.0}
\definecolor{greenpigment}{rgb}{0.0, 0.65, 0.31}
\definecolor{greenhtml}{rgb}{0.0, 0.5, 0.0}

\newcommand{\magenta}[1]{\textcolor{magenta}{#1}}

\newcommand{\green}[1]{\textcolor{greenpigment}{#1}}
\newcommand{\gray}[1]{\textcolor{gray}{#1}}

\newtheorem{assumption}[theorem]{Assumption}
\newtheorem{problem}{Problem}

\newcommand{\None}{
	\left[\begin{array}{c|c}
		I & \begin{array}{c}
			X_+\\Y_- 
		\end{array}
		\\\hline
		0 & \begin{array}{c}
			-X_-\\-U_- 
		\end{array}
	\end{array}\right]
	\!\!\!
	\bbm
	\Phi_{11} & \Phi_{12}\\
	\Phi_{21} & \Phi_{22}
	\ebm\!\!\!
	\left[\begin{array}{c|c}
		I & \begin{array}{c}
			X_+\\Y_- 
		\end{array}
		\\\hline
		0 & \begin{array}{c}
			-X_-\\-U_- 
		\end{array}
	\end{array}\right]^\top\!\!
}

\newcommand{\sysone}{
	\bbm
	I\\\hline\\[-3mm]
	\begin{matrix}
		A & B\\
		C & D
	\end{matrix}
	\ebm
}

\jvol{XX}
\jnum{XX}
\paper{XX}
\jmonth{February}
\jname{submitted to IEEE CONTROL SYSTEMS}
\pubyear{2023}

\begin{document}

\title{The informativity approach\stitle{to data-driven analysis and control}}

\author{HENK J. VAN WAARDE, JAAP EISING, M. KANAT CAMLIBEL, and HARRY L. TRENTELMAN}
\affil{}

\maketitle

\dois{}{}

\chapterinitial{R}oughly speaking, systems and control theory deals with the problem of
making a concrete physical system behave according to certain desired
specifications.  In order to achieve this desired behavior, the system can be 
interconnected with a physical device, called a controller. The problem of finding
a mathematical description of such a controller is called the \emph{control design problem}. 

In order to obtain a mathematical description of a controller for a to-be-controlled 
physical system, a possible first  step is to obtain a mathematical model of 
the physical system. Such a mathematical model can take many forms. For example, 
the model could be in terms of ordinary or partial differential equations, difference equations, 
or transfer matrices. 

There are several ways to obtain a mathematical model for the physical system. 
The usual way is to apply the basic physical laws that are satisfied
by the variables appearing in the system. This method is called \emph{first principles modeling}. 
For example, for electro-mechanical systems, the set of basic physical laws that govern the 
behavior of the variables in the system (conservation laws, Newton's laws, Kirchoff's laws, etc.)  form 
a mathematical model. 

An alternative way to obtain a model is to do experiments on the physical system: certain external 
variables in the physical system are set to take particular values, while at  the same time other variables are
measured.  In this way, one obtains \emph{data} on the system that can be used to find
mathematical descriptions of laws that are obeyed by the system variables, thus obtaining a model. 
This method is called \emph{system identification}. 

The second step in a control system design problem is to decide which desired behavior
we would like the physical system to have.  Very often, this desired behavior can be formalized
by requiring the mathematical model to have certain qualitative or quantitative mathematical properties.
Together, these properties form the \emph{design objective}.

Based on the mathematical model of the physical system and the design objective, the third, 
ultimate, step is to design a mathematical model of a suitable controller. 
This approach, leading from a model and a design objective (or list of design specifications) 
to a model of a controller is an important paradigm in systems and control, and is often 
called \emph{model-based control}. 

An approach that has recently gained popularity is to design controllers without the step of 
finding a mathematical model of the to-be-controlled physical system. This alternative approach deals with the 
problem of synthesizing control laws directly on the basis of measured data, and is called the \emph{data-driven 
approach} to control design. Of course, one can argue that also the combination of system identification followed by model based control as described above is an instance of data driven control design. Indeed, methods using this combination are often called \emph{indirect} methods of data-driven control, consisting of the two-step process of data-driven modeling (i.e., system identification \cite{Ljung1999,vanOverschee1996,Verhaegen2007}) followed by model-based control. 

\begin{pullquote}
regardless of whether a direct or indirect approach is used, any certified data-driven method requires data that contain ``sufficient" information. 
\end{pullquote}

Early contributions to direct data-driven control include PID control \cite{Ziegler1942}, direct adaptive control \cite{Astrom1989}, iterative feedback tuning \cite{Hjalmarsson1994,Hjalmarsson1998}, virtual reference feedback tuning 
\cite{Campi2002} and unfalsified control \cite{Safonov1997}. Recently, direct optimal control design \cite{Bradtke1993,Skelton1994,Furuta1995,Shi1998,Favoreel1999b,Aangenent2005,Markovsky2008,Pang2018,daSilva2019,
Baggio2019,Mukherjee2018,Alemzadeh2019,DePersis2020,Yuan2022} and predictive control \cite{Favoreel1999,Salvador2018,Coulson2019,
Huang2019,Hewing2020,Berberich2021,Alanwar2021,
Yin2021,Coulson2022} have received considerable attention. Some of these approaches, such as \cite{Coulson2019} and \cite{DePersis2020}, strongly rely on the so-called \emph{fundamental lemma} by Willems and co-authors \cite{Willems2005}. This result provides a convenient parameterization of all trajectories of a linear time-invariant system in terms of data. Originally developed in a behavioral context, the fundamental lemma is also instrumental for direct output matching control \cite{Markovsky2008} and control by interconnection \cite{Maupong2017}. The result has been extended in various ways to different model classes and data setups \cite{Berberich2019,vanWaarde2020c,Ferizbegovic2021,
Yin2021b,Yu2021,vanWaarde2021,Coulson2022b,Schmitz2022,Lopez2022}. 
Although initially developed for linear systems, the fundamental lemma has been applied in the context of data-driven control of nonlinear dynamics, such as polynomial and Lur'e systems \cite{Bisoffi2020,Guo2020,Strasser2021,Luppi2022,
Guo2022,Markovsky2022}.
In addition to control problems, also \emph{analysis} problems have been studied within a direct data-based framework. Some examples include the analysis of stability \cite{Park2009}, controllability and observability \cite{Wang2011,Liu2014,Niu2017,Zhou2018}, and dissipativity \cite{Maupong2017b,Romer2017,Romer2019,Koch2020}.

Investigating the different pros and cons of indirect and direct methods is an area of active research \cite{Krishnan2021,Dorfler2022}. However, regardless of whether a direct or indirect approach is used, any certified data-driven method requires data that contain ``sufficient" information. Indeed, as an extreme example one can consider the zero input-output trajectory generated by an unknown linear time-invariant system: clearly this trajectory does not reveal much about the system and would be a poor basis for an identification or control scheme. The purpose of this paper is to introduce a framework in which we can systematically study the richness of data that is required for various system analysis and control problems. 

The concept of ``data informativity" plays a central role in this paper. 
This concept finds its roots in system identification \cite{Ljung1999,Gevers2009,Gevers2013,Colin2020}, where informativity is usually understood as a condition on the data under which it is possible to distinguish between different models in a (parametric) model class. Here, however, we will define a general notion of data informativity \emph{for system analysis and control design}. We are thus not necessarily interested in distinguishing between different models on the basis of the data, but rather want to understand whether it is possible to assess a system-theoretic property, or to synthesize a controller for the physical system from the data. If this is possible, we will say that the data are informative for the system property, or for the control design problem. Although we will introduce the concept of informativity in general terms, it is important to note that the conditions for informative data \emph{depend} on the particular analysis or control problem at hand. For example, as we will see, the conditions under which stabilizing controllers can be obtained from data are less stringent than those for obtaining optimal controllers. This motivates a rigorous analysis of data informativity for different system analysis and control problems. 

In some situations the data obtained from the physical system contain sufficient information to identify the system model uniquely. For our purposes, this situation is the least interesting because informativity for system analysis and control design simply boil down to properties of the (unique) model of the physical system. In general however, it is not possible to uniquely identify the physical system because the data set may contain a small number of samples or the data may be corrupted by noise. In this case, we will make use of \emph{consistent systems sets} that comprise all dynamical system models that are unfalsified by the data. Such sets also play an important role in \emph{set membership identification} \cite{Milanese1991}, where they are called \emph{feasible systems sets}.

In the case that unique system identification is impossible, it will turn out that there are still many relevant cases in which the data are informative for system analysis or control design. In fact, our main results can be interpreted as \emph{robust analysis and control} methods for sets of consistent systems. For example, in our study of noise-free data the approach leads to a robust theory for \emph{affine sets of systems}, which, to the best of our knowledge has not received much attention before. In the noisy data setting, our methods draw inspiration from classical robust control results like Yakubovich's S-lemma \cite{Yakubovich1977} and extensions thereof.

\begin{pullquote}
\emph{informative data} is the cornerstone of this paper, since any data-driven identification, analysis or control task is impossible without it.
\end{pullquote} 

Besides conditions for data informativity, this paper also puts forward a number of \emph{control design methods} that enable the synthesis of a controller from informative data. In many cases, these design methods are based on data-based linear matrix inequalities. The methods thus contribute to the aforementioned lines of work on \emph{direct} data-driven control, since no (explicit) model identification is needed. Although we believe that the direct approach is appealing from a conceptual point of view (why focus on models while the goal is control design?), we also acknowledge that it would be possible to formulate indirect alternatives to our methods by describing the set of systems consistent with the data as a model plus an uncertainty description around this model. In view of this observation, we are inclined to believe that the discussion on direct versus indirect control is, perhaps, not the most fundamental one. Instead, the concept of \emph{informative data} is the cornerstone of this paper, since any data-driven identification, analysis or control task is impossible without it.

\section{Data informativity framework} \label{sec:prob}
In this section we will introduce the concept of data informativity for verifying a given system property or solving a certain control design problem. 

To start with, we fix a certain model class $\mathcal{M}$. This model class is a given set of systems that is assumed to contain the `true' system (i.e., a mathematical model of the underlying unknown physical system), denoted by $\calS$. We assume that the true system $\calS$ is not known but that we do have access to a set of data, $\calD$, which is generated by this system. As explained in the introduction, we are interested in assessing system-theoretic properties of $\calS$ and designing control laws for it from the data $\calD$.
Given the set of data $\calD$, we define $\Sigma_\calD\subseteq\mathcal{M}$ to be the set of all systems in $\calM$ that are consistent with the data $\calD$, i.e., that could also have generated the same data. In other words,  it is impossible to distinguish the true system $\calS$ from any other system in $\Sigma_\calD$ on the basis of the given data $\calD$ alone. As noted before, the setup introduced above is in line with \emph{set membership identification (SMI) methods} \cite{Milanese1991}, where sets of systems consistent with the data also play an important role (in SMI these are typically called \emph{feasible system sets}).

\begin{figure}[h!]
		\centering
			\begin{tikzpicture}[scale=1]
			\node[draw,rectangle,minimum width=8cm,minimum height = 4cm] (1) at (0,0) {};
			\node[] (2) at (3.6,1.7) {$\mathcal{M}$};
			\node[draw,fill,style=circle,inner sep=0pt,minimum size=2pt,color=magenta] (3) at (-2.5,-1) {};
			\node[color=magenta] (4) at (-2.7,-1) {$\mathcal{S}$};
			\draw [fill=greenpigment,fill opacity=0.3] (-1.3,-0.8) ellipse (1.8cm and 1cm);
			\draw [fill=gray,fill opacity=0.3] (0.5,-0.4) ellipse (2.5cm and 1.5cm);
			\node[] at (-1,.4) {\green{$\Sigma_\mathcal{D}$}};
			\node[] at (1.6,1.2) {\gray{$\Sigma_\mathcal{P}$}};
		\end{tikzpicture}
		\caption{The data are not informative for property $\calP$. Depending on the situation, either $\calS \in \Sigma_\calP$ or $\calS \not\in \Sigma_\calP$. On the basis of the given data $\calD$, it is impossible to distinguish these two cases.}
		\label{fig:notinformativedata}
		\end{figure}
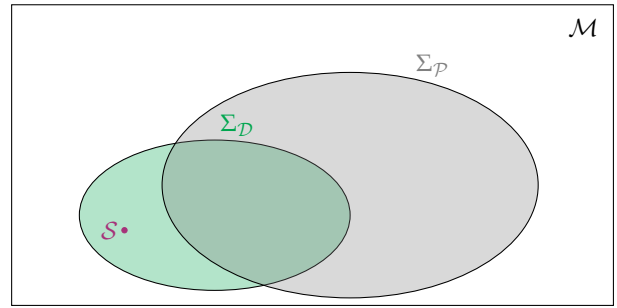

We will now first focus on data-driven analysis of system theoretic properties. Let $\calP$ be some system theoretic property. We will denote the set of all systems within $\mathcal{M}$ having this property by $\Sigma_\calP$.
Suppose we are interested in the question whether our true system $\calS$ has the property $\calP$. Since the only information we have to base our answer on are the data $\calD$ obtained from the true system, we can only conclude from the data that the true system has property $\calP$ if \textit{all} systems consistent with the data $\calD$ have the property $\calP$. If this is the case, we call the data informative for the system property. This leads to the following definition, see also Figures~\ref{fig:informativedata} and \ref{fig:notinformativedata}. 

\begin{definition}[Informativity for analysis]\label{ch1:def:informativity}
	We say that the data $\calD$ are \textit{informative} for property $\calP$ if $\Sigma_\calD \subseteq\Sigma_\calP$, i.e., all systems that are consistent with the data have property $\calP$. 
\end{definition}

\begin{figure}[h!]
		\centering
			\begin{tikzpicture}[scale=1]
			\node[draw,rectangle,minimum width=8cm,minimum height = 4cm] (1) at (0,0) {};
			\node[] (2) at (3.6,1.7) {$\mathcal{M}$};
			\node[draw,fill,style=circle,inner sep=0pt,minimum size=2pt,color=magenta] (3) at (-2.5,-1) {};
			\node[color=magenta] (4) at (-2.7,-1) {$\mathcal{S}$};
			\node[label=right:{$\mathcal{M}$: model class}] at (-4,-2.7) {};
			\node[label=right:{\magenta{$\mathcal{S}$: unknown system}}] at (-4,-3.2) {};
			\node[label=right:{$\mathcal{D}$: given data set}] at (-4,-3.7) {};
			\node[label=right:{\green{$\Sigma_\mathcal{D}$: data consistent systems}}] at (-0.3,-2.7) {};
			\node[label=right:{\gray{$\Sigma_\mathcal{P}$: systems with property $\calP$}}] at (-0.3,-3.7) {};
			\node[label=right:{$\mathcal{P}$: system property}] at (-0.3,-3.2) {};
			\draw [fill=greenpigment,fill opacity=0.3] (-1.3,-0.8) ellipse (1.8cm and 1cm);
			\draw [fill=gray,fill opacity=0.3] (-1,-0.4) ellipse (2.5cm and 1.5cm);
			\node[] at (-1,.4) {\green{$\Sigma_\mathcal{D}$}};
			\node[] at (0.1,1.2) {\gray{$\Sigma_\mathcal{P}$}};
		\end{tikzpicture}
		\caption{The data are informative for property $\calP$ as $\Sigma_\calD \subseteq\Sigma_\calP$.}
		\label{fig:informativedata}
		\end{figure}
		
In general, if the true system $\calS$ can be uniquely determined from the data $\calD$, that is $\Sigma_\calD=\pset{\calS}$ and $\calS$ has the property $\calP$, then it is evident that the data $\calD$ are informative for $\calP$. However, the converse may not be true: $\Sigma_\calD$ might contain many systems, all of which have property $\calP$. In the data informativity framework, we are interested in necessary \textit{and} sufficient conditions for informativity of the data. Such conditions reveal the minimal amount of information required to assess the property $\calP$. A natural problem statement is therefore the following:

\begin{problem}[Informativity problem for analysis]\label{ch3:prob:general}
	Provide necessary and sufficient conditions on the data $\calD$ under which these data are informative for property $\calP$. 
\end{problem}

The above gives us a general framework to deal with data-driven analysis problems. We will also deal with data-driven control problems. The objective in such problems is the data-based design of controllers such that the closed loop system, obtained from the interconnection of the true system $\calS$ and the controller, satisfies the given control objective. 
As for the analysis problem, we have only the information from the data to base our design on. Therefore, we can only guarantee that our control objective is achieved if the designed controller achieves the design objective when interconnected with \textit{any} system from the set $\Sigma_\calD$. 

For the framework to allow for data-driven control problems, we will consider a given control objective 
$\cal{O}$ (for example, a system theoretic property or a guaranteed performance of the closed loop system).  Denote by $\Sigma_{\calO}$ the set of all systems that satisfy the control objective $\calO$. For a given controller $\calK$, denote by $\Sigma_{\calD}(\calK)$ the set of all systems obtained as the interconnection of a system in $\Sigma_{\calD}$ with the controller $\calK$. We then have the following variant of informativity:
\begin{definition}[Informativity for control]\label{ch1:def:par informativity}
We say that the data $\calD$ are \textit{informative} for the control objective $\calO$ if there exists a controller $\calK$ such that $\Sigma_{\calD}(\calK) \subseteq \Sigma_{\calO}$. 
\end{definition}

Obviously, the first step in any data-driven control problem is to determine whether it is possible to obtain, from the given data, a suitable controller. This leads to the following informativity problem:

\begin{problem}[Informativity problem for control]\label{ch1:prob:parametrized}
Provide necessary and sufficient conditions on $\calD$ under which the data are informative for the control objective $\calO$. 
\end{problem}

The second step of data-driven control involves the design of a suitable controller. In terms of our framework, this can be stated as:

\begin{problem}[Control design problem]\label{ch1:prob:design}
	Under the assumption that the data $\calD$ are informative for the control objective $\calO$, find a controller $\calK$ such that $\Sigma_\calD (\calK) \subseteq \Sigma_{\calO}$. 
\end{problem}

The data-driven control design problem as formulated here has a rather natural interpretation as a problem of robust control. Indeed, the aim is to find one single controller that achieves the design objective for all systems  that are consistent with the data. In other words, the `system uncertainty` is determined directly by the given data, and no attempt is made to identify in any sense an uncertainty description that is suitable for existing methods in robust control design. The given data are called informative for a given design objective if the associated robust control problem allows a solution for the system uncertainty imposed by the data.

The data informativity framework has already been applied to various analysis and design problems. Generally speaking, there is a dichotomy between two main directions. On the one hand there are analysis and design problems based on exact, i.e., noiseless data. On the other hand we consider the more realistic situation that the data are noisy, in the sense that they are obtained from a true, unknown, system that is corrupted by additive noise. Table~\ref{tab:summary} provides an overview of the results. The first column of the table states the considered system property or control design problem. The second column refers to the type of data that are used. Here, `E' refers to exact data, and `N' to noisy data. State, input-state, input-state-output and input-output are denoted by `S', `IS', `ISO' and `IO', respectively. The results in the table apply to discrete-time systems, and most of the results are for linear time-invariant dynamics.

\begin{table}[h!]  
\begin{center}
{\small\begin{tabular}{l|c|r}
Problem & Data & References\\
\hline
\textcolor{red}{controllability} & E-IS &\cite{vanWaarde2020}\\
observability & E-S & \cite{Eising2020}\\
\textcolor{red}{stabilizability} & E-IS &\cite{vanWaarde2020}\\
stability & E-S &\cite{vanWaarde2020}\\
\textcolor{red}{state feedback stabilization} & E-IS &\cite{vanWaarde2020,vanWaarde2021-finsler,Werner2022,vanWaarde2022}\\
deadbeat controller & E-IS &\cite{vanWaarde2020}\\
\textcolor{red}{LQR} & E-IS &\cite{vanWaarde2020}\\
suboptimal LQR &E-IS &\cite{vanWaarde2020b}\\
suboptimal $\calH_2$&E-IS &\cite{vanWaarde2020b}\\
synchronization & E-IS & \cite{Jiao2021}\\
dynamic feedback stabilization & E-ISO &\cite{vanWaarde2020}\\
dynamic feedback stabilization & E-IO &\cite{vanWaarde2020}\\
\textcolor{red}{dissipativity} & E-ISO & \cite{vanWaarde2022-dissipativity}\\
\textcolor{red}{tracking and regulation} & E-IS & \cite{Trentelman22}\\
model reduction (moment matching) & E-IO & \cite{Burohman2020}\\
reachability (conic constraints) & E-IO & \cite{Eising2022}\\
stability & N-S & \cite{vanWaarde2022-CLF}\\
stabilizability & N-IS & \cite{vanWaarde2022-CLF}\\
\textcolor{red}{state feedback stabilization} & N-IS&\cite{vanWaarde2022b,vanWaarde2021-finsler,vanWaarde2022}\\
 state feedback $\calH_2$ control & N-IS &\cite{vanWaarde2022b,vanWaarde2022}\\
dynamic feedback $\calH_2$ control & N-IO & \cite{Steentjes2022}\\
\textcolor{red}{state feedback $\calH_\infty$ control} & N-IS &\cite{vanWaarde2022b,vanWaarde2022}\\
dynamic feedback $\calH_\infty$ control & N-IO &\cite{Steentjes2022}\\
stability& N-IO & \cite{vanWaarde2022-QDF}\\
\textcolor{red}{dynamic feedback stabilization} & N-IO & \cite{Steentjes2022,vanWaarde2022-QDF}\\
\textcolor{red}{dissipativity} & N-ISO & \cite{vanWaarde2022-dissipativity}\\
model reduction (balancing) & N-ISO & \cite{Burohman2021}\\
structural properties & N-ISO & \cite{Eising2021}\\
absolute stabilization Lur'e systems & N-ISO & \cite{vanWaarde2022}
\end{tabular}}
\end{center}
\caption{Summary of results within the informativity approach to data-driven analysis and control.}
\label{tab:summary}
\end{table}

The purpose of this paper is to highlight the strength of the informativity framework by reviewing a selection of analysis and design problems, indicated in Table~\ref{tab:summary} in red. Both exact and noisy data will be discussed in the present paper.

This paper is divided into three main sections. These sections are divided into subsections, each devoted to a particular analysis or control design problem. In the first main section our model class will be chosen as the set of all discrete-time linear input-state systems with given state and input dimensions. The data are measurements of the state and input obtained from a true, unknown, system on a given finite time-interval. In this first section it is assumed that the data are noiseless, in the sense that the true system does not contain any noise input. In this noiseless framework we discuss the problem of informativity for the system properties controllability and stabilizability. Next, as a first control design problem we discuss the problem of stabilization by static state feedback, and take a look at the corresponding informativity problem. In the third subsection we study informativity in the context of the linear quadratic regulator problem, and, finally, in the fourth subsection we look at the classical problem of tracking and regulation. 

In the second main section we incorporate noise into the models and data. The model class consists of all input-state systems with additive noise and the input-state data are assumed to be obtained from the noisy true system. The additive noise is unknown, but its samples on the data sampling interval are assumed to satisfy a given quadratic matrix inequality. In this framework we discuss a number of analysis and control design problems. In the first subsection we again look at the problem of stabilization by state feedback, this time in a noisy setting. The next subsection deal with informativity in the context of the well known ${\calH}_{\infty}$ control problem. The final subsection of this part deals with the problem of determining from noisy data whether an unknown system is dissipative with respect to a given supply rate.

In the third, final, main section we abandon the state space framework and shift to input-output systems represented by higher order difference equation, also called autoregressive (AR) systems. As model class we take all AR systems with additive noise, of a given order and with given input and output dimensions. The data are now input-output data that are obtained from the true noisy input-output system. We study data-driven stabilization by dynamic output feedback and discuss informativity in this framework.

\subsection*{Notation}
The set of nonnegative integers will be denoted by $\mathbb{Z}_+$. We will denote by $\mathbb{R}^n$ the $n$-dimensional Euclidean space. For given positive integers $m$ and $n$ the linear space of all real $m \times n$ matrices will be denoted by $\mathbb{R}^{m \times n}$. The subset of $\mathbb{R}^{n \times n}$ consisting of all symmetric  matrices will be denoted by $\S{n}$. For vectors $x$ and $y$ we will denote $\bbm x^\top & y^\top  \ebm^\top$ by $\col(x,y)$. For given integer $n$ we denote by $I_n$ the $n \times n$ identity matrix and $0_n$ the $n \times n$ zero matrix. In order to enhance readability, we sometimes denote  the $n \times m$ zero matrix by $0_{n \times m}$.  Given a real matrix $M$, we will denote its  Moore-Penrose pseudo-inverse by $M\gi$. For a given matrix $M$ of full row rank we denote any right-inverse by $M^{\sharp}$.

\section{Analysis and control using exact input-state data}

In order to provide a solid foundation for more complex problems, in this section we will first consider the essential model class of linear, time-invariant input-state systems. Our goal is to analyze and control these systems on the basis of input-state data, consisting of a finite number of measurements of the input and state trajectories. Moreover, we will assume that these measurements are exact, that is, not corrupted by any noise. 

Given this situation, we can make the abstract framework that was introduced in the introduction more tangible. To be precise, the unknown system $\calS$ is assumed to be the following:
\begin{equation} \label{ch2:eq:true system once more}
	\bmx(t+1) = A_s\bmx(t) + B_s\bmu(t),
\end{equation}
where $\bmx$ denotes the $n$-dimensional state and $\bmu$ the $m$-dimensional input. In the following, we assume that the dimensions $n$ and $m$ are known, but the matrices $A_s$ and $B_s$ are unknown. As such, we see that $\calS$ is contained in the model class $\calM$ given by the set of all discrete-time linear input-state systems of the form
\begin{equation}  \label{ch2:e: is-system}
	\bmx(t+1) = A\bmx(t) + B\bmu(t),
\end{equation}
with given state space and input dimensions $n$ and $m$. 

Suppose that we collect input-state data from the true system \eqref{ch2:eq:true system once more} on a set of time instances $\pset{0,1,\ldots,T}$, in the sense that we excite the true system with an input sequence $u(0)$, $u(1)$, $\ldots$, $u(T-1)$ and obtain measurements of a corresponding state sequence $x(0)$, $x(1)$, $\ldots$, $x(T)$. We can collect these measurements in matrices by defining:
\begin{subequations}\label{ch2:eq: UXdata}
	\begin{align}
		U_-& := \bbm u(0) & u(1) & \cdots & u(T-1)\ebm, \\
		X& := \bbm x(0) & x(1) & \cdots & x(T)\ebm.
	\end{align}
\end{subequations}
If, additionally, we define the matrices
\begin{subequations}\label{ch2:eq: def of X- X+}
	\begin{align}
		X_-& := \bbm x(0) & x(1) & \cdots & x(T-1) \ebm, \\
		X_+& := \bbm x(1) & x(2) & \cdots & x(T) \ebm,
	\end{align}	
\end{subequations}
we have $X_+=A_sX_-+B_sU_-$, since the data were assumed to be generated by the true system. Moreover, this is \textit{all} the information we have regarding the true system on the basis of the data $\calD:=(U_-,X)$. As explained in the introduction, we are interested in the set $\Sigma_\calD$ containing all systems in $\calM$ that are consistent with these data. Obviously, this set is given by 
\begin{equation}\label{ch2:eq:SigmaD}
	\Sigma_\calD = \left\{ (A,B) \in \calM \mid X_+= \bbm A&B \ebm
	\begin{bmatrix}
		X_-\\U_-
	\end{bmatrix} \right\}. 
\end{equation}

In fact, there is no need for the data to be collected sequentially. It is straightforward to adapt the results above to the situation of measurements on multiple sets of time instances.

First, we consider the problem of system identification. In the terminology of this paper, we say that the data $(U_-,X)$ are \textit{informative for system identification} if the set $\Sigma_\calD$ contains exactly one element. Since by definition the true system $(A_s,B_s) \in \Sigma_\calD$ the data are informative for identification if and only if $\Sigma_\calD = \{(A_s,B_s)\}$. Moreover it is the solution set of the affine equation appearing in \eqref{ch2:eq:SigmaD}. Therefore, the data are informative for system identification if and only if the full rank condition
\begin{equation}\label{ch2:eq:inf for sys ident} 
\rank \begin{bmatrix} X_- \\ U_- \end{bmatrix} = n+m. 
\end{equation}
holds. There exists a unique system in $\Sigma_\calD$ if and only if the condition \eqref{ch2:eq:inf for sys ident} holds, and this true system can then be obtained from the data as
\[
\bbm A_s & B_s \ebm = X_+ \begin{bmatrix}
		X_-\\U_-
	\end{bmatrix}^\sharp.
\]	
Once we have identified the true system is this way, we can apply model-based methods in order to verify its properties or to achieve the desired control objective. In the Sidebar ``Willems' fundamental lemma'' we discuss the problem of designing inputs such that the resulting measurements are guaranteed to be informative for system identification.

\begin{sidebar}{Willems' fundamental lemma}\label{sidebar:fundamental}
	\setcounter{sequation}{0}
	\renewcommand{\thesequation}{S\arabic{sequation}}
	\setcounter{stable}{0}
	\renewcommand{\thestable}{S\arabic{stable}}
	\setcounter{sfigure}{0}
	\renewcommand{\thesfigure}{S\arabic{sfigure}}
	
	\sdbarinitial{A}s explained in the main text, the input-state data satisfy the full row rank condition \eqref{ch2:eq:inf for sys ident} if and only if the data are informative for identification, meaning that the true system can be uniquely determined from the data by solving a linear equation. This brings up the question whether it is possible to {\em choose} a time horizon $T$ and a finite input sequence $u(0), u(1), \ldots, u(T -1)$ such that condition \eqref{ch2:eq:inf for sys ident} is guaranteed to hold for any resulting state sequence. Indeed, if this could be done, then the true system could be uniquely determined from the data by choosing a suitable input sequence. In this sidebar we will discuss this question in the context of \textit{persistently exciting} inputs and Willems' fundamental lemma.
	
In the sequel, denote any given finite sequence $f(0)$, $f(1)$, $\ldots, f(T-1)$ by $f_{[0,T-1]}$. For a given finite input sequence $u_{[0,T-1]}$, define the associated Hankel matrix of depth $k$ by	
\[
H_k(u_{[0,T-1]}) := \begin{bmatrix}
		u(0) 	& u(1) 		& \cdots 	& u(T - k) \\
		u(1) 	& u(2) 		& \cdots 	& u(T - k +1) \\
		\vdots  & \vdots	&    		& \vdots  \\
		u(k-1)  & u(k) 		& \cdots  	& u(T-1) 
	\end{bmatrix}.
 \]
The input sequence $u_{[0,T-1]}$ is called \textit{persistently exciting} of order $k$ if $H_k(u_{[0,T-1]})$ has full row rank. 
	
	A special case of Willems' fundamental lemma \cite[Thm. 1]{vanWaarde2020c} (originally proven in a behavioral context in \cite[Thm. 1]{Willems2005}) is relevant in the context of informativity for identification.  Consider a linear input-state-output  system, defined by the quadrupel of matrices $(A,B,C,D)$, with input, state and output denoted by $u$, $x$ and $y$ respectively. Let $x_{[0,T-1]}$ and $y_{[0,T-1]}$ denote the finite length state and output trajectories corresponding to the input sequence $u_{[0,T-1]}$. Suppose that the system is controllable and observable, and that the input sequence $u_{[0,T-1]}$ is persistently exciting of order $n+L$. Denote $X_L =\begin{bmatrix} x(0) & \cdots & x(T-L) \end{bmatrix}$. Then a consequence of Willems' fundamental lemma is that the matrix
\begin{sequation}
\label{stateinputHankel}
\begin{bmatrix} X_L \\ H_L(u_{[0,T-1]}) \end{bmatrix}, 
\end{sequation}
has full row rank. A special case of this arises for $L=1$, which shows that \eqref{ch2:eq:inf for sys ident} holds if we take $u_{[0,T-1]}$ to be persistently exciting of degree $n+1$. This resolves the problem of designing inputs for informativity for identification in the context of input-state measurements. More generally, in the case of input-output measurements, full row rank of \eqref{stateinputHankel} enables the identification of the system matrices $(A,B,C,D)$ up to similarity transformation if $L$ is larger than the so-called \emph{lag} of the system. This can be done, for example, by using subspace identification methods \cite{vanOverschee1996,Verhaegen2007}.

On the other hand, in many applications an identified (state space) model might not be the most convenient representation to work with. One of these applications is Data-enabled Predictive control (DeePC), introduced in \cite{Coulson2019}. As a data-driven variant of Model Predictive Control (MPC), the central problem is to minimize a cost function over all length-$L$ trajectories. A consequence of the fundamental lemma is that if the input sequence $u_{[0,T-1]}$ is persistently exciting of order $n+L$, then any  $\bar{u}_{[0,L-1]},\bar{y}_{[0,L-1]}$ is an input/output trajectory of the system if and only if 
	\[ \begin{bmatrix} \bar{u}_{[0,L-1]} \\ \bar{y}_{[0,L-1]} \end{bmatrix} \in \im \begin{bmatrix} H_L(u_{[0,T-1]}) \\ H_L(y_{[0,T-1]}) \end{bmatrix}. \] 
	As such, the fundamental lemma allows this problem to be formulated directly in terms of measurements, without explicitly finding a system model.  
	
	As it turns out, there are many other applications where avoiding the modeling step, and dealing with data directly is convenient. For example, the paper \cite{Markovsky2008} treats data-based simulation and output matching control, while \cite{DePersis2020} provides methods for stabilization on the basis measurements for which \eqref{ch2:eq:inf for sys ident} holds. 
\end{sidebar}

\subsection{Controllability and stabilizability}

As we will show in this subsection, condition \eqref{ch2:eq:inf for sys ident} is not necessary to perform data-driven analysis in general. As an illustration, we will establish necessary and sufficient conditions in terms of the data for verifying controllability and stabilizability, which do not require the rank condition to hold.

Recall from Definition~\ref{ch1:def:informativity} the definition of informativity of data for a given system property. In accordance with Definition~\ref{ch1:def:informativity} we say that the data $(U_-,X)$ are {\em informative for controllability} if all systems in $\Sigma_\calD$ are controllable. Likewise, we call the data  {\em informative for stabilizability} if all systems in $\Sigma_\calD$ are stabilizable.

In order to establish tests for these notions of informativity, the well known Hautus test 
for controllability can be used: a system $(A,B)$ is controllable if and only if
\begin{equation}\label{ch2:eq:Hautus} \rank \begin{bmatrix} A-\lambda I & B \end{bmatrix} =n \end{equation} 
for all $\lambda\in \mathbb{C}$. For stabilizability, the Hautus test requires that \eqref{ch2:eq:Hautus} holds for all $\lambda$ outside the open unit disc.


The following theorem gives necessary and sufficient conditions on the input-state data to be informative for these two properties. The result provides tests on the given data matrices.

\begin{theorem}[Data-driven Hautus tests \cite{vanWaarde2020}]
	\label{ch2:t:contstab}
	The data $(U_-,X)$ are informative for controllability if and only if 
	\begin{equation}
		\label{ch2:eq:rank for cont}
		\rank ( X_+ -\lambda X_-) = n \quad \forall \lambda \in \mathbb{C}.
	\end{equation}
	Similarly, the data $(U_-,X)$ are informative for stabilizability if and only if 
	\begin{equation}\label{ch2:eq:rank for stab}
		\rank ( X_+ -\lambda X_-) = n \quad \forall \lambda \in \mathbb{C}\textrm{ with } |\lambda|\geq 1.
	\end{equation}
\end{theorem}

As announced at the beginning of this section, there are situations in which we can conclude controllability or stabilizability from the data without being able to identify the true system uniquely. This can be seen by the fact that in order for \eqref{ch2:eq:inf for sys ident} to hold we require at least $m+n$ separate measurements, whereas the conditions of Theorem~\ref{ch2:t:contstab} can hold for $n$ measurements. This is illustrated in the following example.

\begin{example}[Full rank is not necessary]
	Suppose that $n=2$ and $m=1$. Assume we collect data on the single time interval $\{0,1 \ldots, T\}$ with $T = 2$ to obtain
	\[ X= \begin{bmatrix} 0 & 1 & 0 \\ 0 & 0 & 1 \end{bmatrix} \textrm{ and } U_-= \begin{bmatrix} 1 & 0 \end{bmatrix}. \] 
	This implies that 
	\[X_+ = \begin{bmatrix} 1 & 0 \\ 0 & 1 \end{bmatrix} \textrm{ and } X_- = \begin{bmatrix} 0 & 1 \\ 0& 0\end{bmatrix}. \]
	Using Theorem~\ref{ch2:t:contstab} we see that these data are informative for controllability, as 
	\[ \rank \begin{bmatrix} 1 & -\lambda \\ 0 & 1 \end{bmatrix}= 2 \quad \forall \lambda \in \mathbb{C}.\] 
	Recall that this means that all systems consistent with the data are controllable. Therefore we can conclude that also the true system is controllable. Moreover, note that the rank condition \eqref{ch2:eq:inf for sys ident} does not hold, and therefore $\Sigma_\calD$ is not a singleton.  To be precise:
	\begin{equation}\label{ch2:eq:structural} \Sigma_\calD = \left\lbrace \left( \begin{bmatrix} 0 & a_1 \\ 1 & a_2\end{bmatrix} ,\begin{bmatrix} 1 \\ 0 \end{bmatrix}\right) \mid a_1,a_2 \in\mathbb{R} \right\rbrace.
	\end{equation}
	This means that there are multiple systems consistent with the data.
\end{example}

Computationally, the conditions \eqref{ch2:eq:rank for cont} and \eqref{ch2:eq:rank for stab} might seem daunting, since these require to test the rank of a matrix for each $\lambda\in\mathbb{C}$. However, it is well known that in order for the classical Hautus test to be satisfied, it suffices to test the rank of \eqref{ch2:eq:Hautus} for only $\lambda\in \sigma(A)$, where $\sigma(A)$ denotes the set of eigenvalues of the matrix $A$. 

In a similar fashion, the conditions of Theorem~\ref{ch2:t:contstab} can be verified in a finite number of steps. Indeed, \eqref{ch2:eq:rank for cont} is equivalent to 
\[ \rank(X_+)=n \textrm{ and } \rank ( X_+ -\lambda X_-) = n \]
for all $\lambda\neq 0$ with $\lambda\inv\in\sigma(X_-X_+^\sharp)$, where $X_+^\sharp$ is any right inverse of $X_+$. Regarding stabilizability, we obtain that \eqref{ch2:eq:rank for stab} is equivalent to 
\[\rank(X_+-X_-)=n \textrm { and } \rank ( X_+ -\lambda X_-) = n \]
for all $\lambda\neq 1$ with $(\lambda-1)\inv\in\sigma(X_-(X_+-X_-)^\sharp)$.


\subsection{Stabilization}
After considering data-driven controllability and stabilizability analysis in the previous subsection, we now turn attention to data-driven control design. In particular, we consider the quintessential control problem of stabilization by state feedback. 

Recall the definition of informativity for control as given in Definition~\ref{ch1:def:par informativity}. We take the model class $\calM$ and data $\calD$ as before, and take as the control objective $\calO$: `interconnection with a  state feedback controller yields a stable\footnote{meaning Schur, that is, all its eigenvalues $\lambda$ satisfy $\abs{\lambda}<1$.}, closed loop system'. This means that the set $\Sigma_{\calO}$ of all systems that satisfy the control objective is equal to the set of all stable $n \times n$ matrices
$$
M^{n \times n}_{\rm stab}:=  \{ A \in \R^{n \times n} \mid A \text{ is stable} \}.
$$
For a given state feedback controller $K\in\mathbb{R}^{m\times n}$, the corresponding set of closed loop systems consistent with the data is equal to 
\[ 
\Sigma_\calD(K) = \{ A+ BK \mid (A,B) \in \Sigma_\calD \}. 
\] 
In line with Definition \ref{ch1:def:par informativity} we say that the data $(U_-,X)$ are {\em informative for stabilization by state feedback} if there exists a $K \in \mathbb{R}^{m \times n}$ such that $\Sigma_\calD(K) \subseteq M^{n \times n}_{\rm stab}$.

In other words, the input-state data $(U_-,X)$ are informative for stabilization by state feedback if there exists a single real $m \times n$ matrix $K$ such that $A +BK$ is stable for all $(A,B) \in \calM$ that are consistent with the data.

At this point, one may wonder about the relation between informativity for stabilizability and informativity for stabilization. It is clear that the data $(U_-,X)$ are informative for stabilizability if $(U_-,X)$ are informative for stabilization by state feedback. However, the reverse statement does not hold in general. This is due to the fact that all systems $(A,B)$ in $\Sigma_\calD$ may be stabilizable, but there may not exist a \emph{common} feedback gain $K$ such that $A+BK$ is stable for all of these systems. 

The following example further illustrates the difference between informativity for stabilizability and informativity for stabilization.

\begin{example}[Stabilizability and stabilization]
	Consider the scalar system 
	$
	\bmx(t+1)= \bmu(t),
	$
	where $\bmx, \bmu \in \mathbb{R}$. Suppose that we collect data on the single time interval $\{0,1\}$, specifically, $x(0) = 0$, $u(0) = 1$ and $x(1) = 1$. This means that $U_- = \begin{bmatrix}1 \end{bmatrix}$ and $X = \begin{bmatrix}0 & 1	\end{bmatrix}$. It can be shown that $\Sigma_\calD = \{(a,1) \mid a \in \mathbb{R}\}$. Clearly, all systems in $\Sigma_\calD$ are stabilizable. Nonetheless, the data are not informative for \emph{stabilization}. This is because the systems $(-1,1)$ and $(1,1)$ in $\Sigma_\calD$ cannot be stabilized by the \emph{same} controller of the form $u(t) = K x(t)$. 
	We conclude that informativity of the data for stabilizability does not imply informativity for stabilization by state feedback. 
\end{example}

Having defined the notion of informativity for stabilization, we now take the steps described in the introduction. First, we resolve Problem~\ref{ch1:prob:parametrized}, that is, we find necessary and sufficient conditions for informativity for stabilization by state feedback. After this, we design a corresponding controller, as described in Problem~\ref{ch1:prob:design}. 

In order to do this, we first state a useful lemma. Recall from \eqref{ch2:eq:SigmaD} that $(A,B)\in\Sigma_\calD$ if and only if it is a solution of the the corresponding affine equation. Now let $\Sigma_\calD^0$ denote the solution set of the corresponding homogeneous equation. That is, \begin{equation}\label{ch3:eq:sigma is 0}
	\Sigma_\calD^0 := \left\{(A_0,B_0) \mid 0 = \bbm A_0 & B_0 \ebm
	\begin{bmatrix}
		X_- \\ U_-
	\end{bmatrix}\right\}.
\end{equation}
This allows us to state the following lemma. 
\begin{lemma}[A necessary condition \cite{vanWaarde2020}]
	\label{ch3:lemmaF0=0}
	Suppose that the data $(U_-,X)$ are informative for stabilization by state feedback, and let $K$ be a  feedback gain such that $A+BK$ is stable for all $(A,B) \in \Sigma_\calD$. Then $A_0 + B_0 K = 0$ for all $(A_0,B_0) \in \Sigma_\calD^0$. Equivalently, 
	$$
	\im \begin{bmatrix}
		I \\ K
	\end{bmatrix} \subseteq 
	\im \begin{bmatrix}
		X_- \\ U_-
	\end{bmatrix}.
	$$
\end{lemma}

The solution set of an affine equation is equal to the sum of any solution and the solution set of the corresponding homogenous equation. Since we know that $(A_s,B_s)\in\Sigma_\calD$ by definition, this means that we can write 
\[ \Sigma_\calD = (A_s,B_s) + \Sigma_\calD^0.\]
Using this, as a consequence of Lemma~\ref{ch3:lemmaF0=0} we have that if $K$ is a  feedback gain such that $A+BK$ is stable for all $(A,B) \in \Sigma_\calD$, then
\[ \Sigma_\calD(K) = \{A_s+B_sK\},  \] 
that is, the set of closed-loop systems consistent with the data is a singleton. It is important to note, however, that this does not mean that $\Sigma_\calD$ is necessarily a singleton.

The above observation turns out to be instrumental in proving the following theorem, which gives necessary and sufficient conditions for informativity for stabilization by state feedback. 

\begin{theorem}[Conditions for stabilization \cite{vanWaarde2020}]
	\label{ch3:t:algstab}
	The data $(U_-,X)$ are informative for stabilization by state feedback if and only if the matrix $X_-$ has full row rank and there exists a right inverse $X_-^\sharp$ of $X_-$ such that $X_+ X_-^\sharp$ is stable. 
	
	Moreover, $K$ is such that $A +BK$ is stable for all $(A,B) \in \Sigma_\calD$ if and only if $K = U_- X_-^\sharp$, where $X_-^\sharp$ satisfies the above properties. In that case, $A +BK = X_+ X_-^\sharp$ for all $(A,B) \in \Sigma_\calD$.
\end{theorem}

Theorem \ref{ch3:t:algstab} gives a characterization of all input-state data that are informative for stabilization by state feedback and provides a stabilizing controller. Nonetheless, the procedure to compute this controller might not be entirely satisfactory since it is not clear how to find a right inverse of $X_-$ that makes $X_+ X_-^\sharp$ stable. In general, $X_-$ has many right inverses, and $X_+ X_-^\sharp$ can be stable or unstable depending on the particular right inverse $X_-^\sharp$. To deal with this problem and to solve the design problem, we give a characterization of informativity for stabilization in terms of linear matrix inequalities (LMIs). The feasibility of such LMIs can be verified using standard tools. 

LMI conditions for data-driven stabilization were first studied in \cite{DePersis2020}. In that work, the following conditions were presented under the additional assumption that the measurements satisfy \eqref{ch2:eq:inf for sys ident}, that is, the measurements are informative for system identification. As it turns out, this assumption can be removed, leading to the following result. 

\begin{theorem}[LMI conditions for stabilization \cite{vanWaarde2020}]\label{ch3:t:lmistab}
	The data $(U_-,X)$ are informative for stabilization by state feedback if and only if there exists a matrix $\Theta \in \mathbb{R}^{T \times n}$ satisfying
	\begin{equation}
		\label{ch3:LMI/E}
		X_- \Theta = (X_- \Theta)^\top\quad\text{ and }\quad
		\begin{bmatrix}
			X_- \Theta & X_+ \Theta \\ \Theta^\top X_+^\top & X_- \Theta
		\end{bmatrix} > 0.
	\end{equation}
	Moreover, $K$ is such that $A +BK$ is stable for all $(A,B) \in \Sigma_\calD$ if and only if $K = U_- \Theta (X_-\Theta)^{-1}$ for some matrix $\Theta$ satisfying \eqref{ch3:LMI/E}.
\end{theorem}

The following example provides a simple illustration of the above results. 
\begin{example}[Full rank not necessary for informativity] 
	Consider an unstable system $(A_s,B_s)$, where $A_s$ and $B_s$ are given by
	$$
	A_s = \begin{bmatrix}
		1.5 & 0 \\ 1 & 0.5
	\end{bmatrix}, \quad B_s = \begin{bmatrix}
		1 \\ 0
	\end{bmatrix}.
	$$
	We collect data from this system on a single time interval from $t = 0$ until $t = 2$, which results in the data matrices
	$$
	X = \begin{bmatrix}
		1 & 0.5 & -0.25 \\
		0 & 1 & 1
	\end{bmatrix}, \quad U_- = \begin{bmatrix}
		-1 & -1
	\end{bmatrix}.
	$$
	Clearly, the matrix $X_-$ is square and invertible, and it can be verified that
	$$
	X_+ X_-^{-1} = \begin{bmatrix}
		0.5 & -0.5 \\
		1 & 0.5
	\end{bmatrix}
	$$
	is stable, since its eigenvalues are $\half(1 \pm \sqrt{2}i)$. We conclude by Theorem \ref{ch3:t:algstab} that the data $(U_-,X)$ are informative for stabilization by state feedback. The same conclusion can be drawn from Theorem \ref{ch3:t:lmistab} since $$\Theta = \begin{bmatrix}
		1 & -1 \\ 0 & 2
	\end{bmatrix}$$ solves \eqref{ch3:LMI/E}. Next, we can conclude from either Theorem \ref{ch3:t:algstab} or Theorem \ref{ch3:t:lmistab} that the stabilizing feedback gain in this example is unique, and given by $K = U_-X_-^{-1} = \begin{bmatrix}
		-1 & -0.5
	\end{bmatrix}$. Finally, it is worth noting that the data are not informative for identification. In fact, $(A,B) \in \Sigma_\calD$ if and only if
	$$
	A = \begin{bmatrix}
		1.5+a_1 & 0.5a_1 \\ 1+ a_2 & 0.5+ 0.5a_2
	\end{bmatrix}, \quad B = \begin{bmatrix}
		1+a_1 \\ a_2
	\end{bmatrix}
	$$
	for some $a_1$ and $a_2 \in \mathbb{R}$.
\end{example}

\begin{sidebar}{The linear quadratic regulator problem}
\label{sidebar:LQR}

Consider the discrete time linear system 
\begin{sequation} \label{ch4:e:disc}
\bmx(t+1) = A \bmx(t) + B \bmu(t),
\end{sequation}
where $A$  and $B$ are matrices of dimensions $n \times n$ and $n \times m$, and where 
$\bmx$ is the $n$-dimensional state and $\bmu$ the $m$-dimensional input. In the linear quadratic regulator problem we quantify the performance of the system using a quadratic cost functional $J(x_0,\bmu)$ involving the state trajectory $\bmx$ and the input $\bmu$. The optimal linear quadratic regulator problem is then the problem of finding, for each initial state $x_0$ of the system, an optimal input, i.e. an input that minimizes the cost functional. In this sidebar the basics of discrete-time linear quadratic optimal control are reviewed. In the sequel,  the abbreviation `LQR'  will be used for  `linear quadratic regulator'. 

For an initial state $x_0$, let $\bmx_{x_0,\bmu}$ be the state sequence of \eqref{ch4:e:disc} resulting from the input $\bmu$ and initial condition $\bmx(0) = x_0$. We omit the subscript and simply write $\bmx$ whenever the dependence on $x_0$ and $\bmu$ is clear from the context. 

Associated to system \eqref{ch4:e:disc}, we define the quadratic cost functional
\begin{sequation}\label{ch4:e:cost}
J(x_0,\bmu)=\sum_{t=0}^\infty  \bmx^\top(t) Q \bmx(t) + \bmu^\top(t) R \bmu(t),
\end{sequation}
where $Q \in \S{n}$ is positive semidefinite and $R \in \S{m}$ is positive definite. Then, the optimal LQR problem is the following: 
\begin{problem}[The LQR problem]
	Determine for every initial condition $x_0$ an input $\bmu^*$, such that $\lim_{t\to\infty} \bmx_{x_0,\bmu^*}(t) = 0$, and the cost functional $J(x_0,\bmu)$ is minimized under this constraint. 
\end{problem}
\noindent Such an input $\bmu^*$ is called optimal for the given $x_0$. Of course, an optimal input does not necessarily exist for all $x_0$. We say that the optimal LQR problem is {\em solvable\/} for $(A,B,Q,R)$ if for every $x_0$ there exists an input ${\bmu}^*$ such that
\begin{enumerate}
	\item The cost $J(x_0,\bmu^*)$ is finite.
	\item The limit $\lim_{t\to\infty}\bmx_{x_0,\bmu^*}(t)=0$. 
	\item The input $\bmu^*$ minimizes the cost functional, i.e., 
	\[J(x_0,{\bmu}^*)\leq J(x_0,\bar{\bmu})\]
	for all $\bar{\bmu}$ such that $\lim_{t\to\infty}\bmx_{x_0,\bar{\bmu}}(t)=0$.
\end{enumerate}
In the sequel, we will require the notion of observable eigenvalues. 
An eigenvalue $\lambda$ of $A$ is called $(Q,A)$-observable if 
\[ 
\rank \begin{pmatrix} A-\lambda I \\ Q \end{pmatrix}=n.
\] 
The following theorem provides necessary and sufficient conditions for the solvability of the optimal LQR problem for $(A,B,Q,R)$. 
This theorem is the discrete-time analogue to the continuous-time case stated in \cite[Thm. 10.18]{Trentelman2001}. 
\begin{theorem}[Conditions for LQR]\label{ch4:t:Harry}
	Let $Q \geq 0$ and $R >0$. Then the following statements hold:	
	\begin{enumerate}
		\item If $(A,B)$ is stabilizable, there exists a unique largest real symmetric solution $P^+$ to the discrete-time algebraic Riccati equation (DARE) 
		\begin{sequation}
		\label{ch4:dare}
		P = A^\top PA-A^\top PB(R+B^\top P B)\inv B^\top  P A+Q,
		\end{sequation}
		in the sense that $P^+ \geq P$ for every real symmetric $P$ satisfying \eqref{ch4:dare}. The matrix $P^+$ is positive semidefinite.
		\item If, in addition to stabilizability of $(A,B)$, every eigenvalue of $A$ on the unit circle is $(Q,A)$-observable then for every $x_0$ a unique optimal input $\bmu^*$ exists. Furthermore, this input sequence is generated by the feedback law $\bmu = K \bmx$, where
		\begin{sequation}
		\label{ch4:optgain}
		K := -(R+B^\top P^+ B)\inv B^\top  P^+ A.
		\end{sequation}
		Moreover, the matrix $A+BK$ is stable. 
		\item In fact, the optimal LQR problem is solvable for $(A,B,Q,R)$ if and only if $(A,B)$ is stabilizable and every eigenvalue of $A$ on the unit circle is $(Q,A)$-observable. 
	\end{enumerate}
\end{theorem}

If the optimal LQR problem is solvable for $(A,B,Q,R)$, we say that the matrix $K$ given by \eqref{ch4:optgain} is the optimal feedback gain for $(A,B,Q,R)$. 

\end{sidebar}

\subsection{The linear quadratic regulator problem}

An important classical control design problem is the optimal linear quadratic regulator (LQR) problem. In this subsection we will study the data-driven version of this problem within the informativity framework.  

For given state and input dimensions $n$ and $m$, again consider the model class $\mathcal{M}$ of all discrete-time linear input-state systems \eqref{ch2:e: is-system}.
Assume we have input-state data on multiple time intervals, leading to data $\calD:= (U_-,X)$ as given in \eqref{ch2:eq: UXdata}. As before, the set $\Sigma_{\calD}$ of all systems in $\calM$ that are consistent with the data is then given by \eqref{ch2:eq:SigmaD}.
We assume that the data are generated by the true (but unknown) system $(A_s,B_s)$, which is therefore in $\Sigma_\calD$ itself.

In the context of the optimal LQR problem the control objective $\calO$ is: `the system must be  controlled using the optimal feedback gain'.  In order to formalize this, we introduce the following notation.
For any given $K$, let $\Sigma^{Q,R}_{K}$ denote the set of all systems of the form \eqref{ch2:e: is-system} for which $K$ is the optimal feedback gain corresponding to $Q$ and $R$, that is,
\[ 
\Sigma_K^{Q,R}:=\set{(A,B) \in \calM }{K \text{ is optimal for }(A,B,Q,R)}.
\]
This gives rise to yet another notion of informativity in line with Definition~\ref{ch1:def:par informativity}. Indeed, informativity requires the existence of a single feedback gain that is optimal for all systems consistent with the data. For the definition of solvability of the optimal LQR problem we refer to the sidebar `The linear quadratic regulator problem`.
\begin{definition}[Informativity for LQR]
	Given matrices $Q$ and $R$, we say that the data $\calD = (U_-,X)$ are \emph{informative for optimal linear quadratic regulation} if the optimal LQR problem is solvable for all $(A,B) \in \Sigma_{\calD} $ and there exists $K$ such that $\Sigma_{\calD} \subseteq \Sigma^{Q,R}_{K}$.
\end{definition}
An instrumental result in obtaining necessary and sufficient conditions for informativity for optimal linear quadratic regulation is the following lemma.
\begin{lemma}[Common solution of the Riccati equation]\label{ch4:l: same P works for all}
Let $Q=Q^\top$ be positive semidefinite and $R=R^\top$ be positive definite. Suppose the data $(U_-,X)$ are informative for optimal linear quadratic regulation. Let $K$ be such that $\Sigma_{\calD} \subseteq \Sigma^{Q,R}_{K}$. Then, there exist a square matrix $M$ and a positive semidefinite matrix $P^+$ such that for all $(A,B)\in\Sigma_{\calD}$
\begin{align}
&\!\!\!\!\!M=A+BK,\label{ch4:e: define M}\\
&\!\!\!\!\!P^+\!= A^\top\! P^+\!A\! -\! A^\top\! P^+\!B(R + B^\top\! P^+\! B)\inv B^\top\!  P^+\! A + Q,\!\!\label{ch4:e:dare aux}\\
&\!\!\!\!\!P^+-M^\top P^+M=K^\top RK+Q,\label{ch4:e:lyap aux}\\
&\!\!\!\!\!K=-(R+B^\top P^+ B)\inv B^\top  P^+ A.\label{ch4:e:uni K aux}
\end{align}
\end{lemma}
Statement \eqref{ch4:e:dare aux} of the lemma says that if the data are informative, there exists a common solution $P^+ \geq 0$ to the whole collection of AREs associated with systems $(A,B)$ that are consistent with the data. Statement \eqref{ch4:e:uni K aux} 
says that if $K$ is the common optimal gain for all systems that are consistent with the data, then it must be of the expected form \eqref{ch4:e:uni K aux} for all  $(A,B)$ consistent with the data. According to \eqref{ch4:e: define M}, the optimal closed loop system matrices $A + BK$ are identical for all consistent pairs $(A,B)$.

The following theorem gives necessary and sufficient conditions for informativity for optimal linear quadratic regulation. 
\begin{theorem}[Conditions for informativity \cite{vanWaarde2020}]\label{ch4:t:LQinform}
	Let $Q \geq 0$ and $R > 0$. Then the data $(U_-,X)$ are informative for optimal linear quadratic regulation if and only if at least one of the following two conditions hold:
	\begin{enumerate}
		\item\label{ch4:cond:a} The data $(U_-,X)$ are informative for identification, that is, $\Sigma_{\calD}=\pset{(A_s,B_s)}$, and the optimal LQR problem is solvable for $(A_s,B_s,Q,R)$. In this case, the optimal feedback gain $K$ is of the form \eqref{ch4:e:uni K aux} where $P^+$ is the largest real symmetric solution to \eqref{ch4:e:dare aux} with $A = A_s$ and $B = B_s$.
		\item\label{ch4:cond:b} For all $(A,B) \in \Sigma_{\calD}$ we have $A=A_s$. Moreover, $A_s$ is stable, $QA_s = 0$, and the optimal feedback gain is given by $K = 0$. 
	\end{enumerate}
\end{theorem}

This theorem should be interpreted as follows. Condition \ref{ch4:cond:b}) of Theorem \ref{ch4:t:LQinform} can be considered as a  pathological case in which the only $A$ consistent with the data is the true one, namely $A_s$. This matrix $A_s$ is stable and $QA_s = 0$. Since $\bmx(t) \in \im A_s$ for all $t > 0$, we have $Q \bmx(t) = 0$ for all $t > 0$ if the input function is chosen as $\bmu = 0$. Additionally, since $A_s$ is stable, this shows that the optimal input is equal to $\bmu^* = 0$. If we set aside the pathological case \ref{ch4:cond:b}), the main message of Theorem \ref{ch4:t:LQinform} is the following: if the data are informative for optimal linear quadratic regulation they are also informative for system identification, in the sense that the set of systems consistent with the data contains only one element, i.e., $\Sigma_{\calD} = \{(A_s,B_s) \}$. This observation is consistent with the paper \cite{Polderman1986} that showed the necessity of identifiability of the true system in adaptive LQ control.

At first sight, this might seem like a negative result in the sense that data-driven LQR is only possible with data that are also informative enough to uniquely identify the system. However, at the same time, Theorem \ref{ch4:t:LQinform} can be viewed as a positive result in the sense that it provides fundamental justification for the data conditions imposed in e.g. \cite{DePersis2020}. Indeed, in \cite{DePersis2020} the data-driven infinite horizon LQR problem\footnote{Note that the authors of \cite{DePersis2020} formulate this problem as the minimization of the $H_2$-norm of a certain transfer matrix.} is solved using input-state data under the assumption that the input is persistently exciting of sufficiently high order. Under the latter assumption, the input-state data are informative for system identification, i.e., the matrices $A_s$ and $B_s$ can be uniquely determined from data. Theorem \ref{ch4:t:LQinform} justifies such a strong assumption on the richness of data in data-driven linear quadratic regulation.
The data-driven \emph{finite} horizon LQR problem was solved under a persistency of excitation assumption in \cite{Markovsky2007}. Our results suggest that also in this case informativity for system identification is necessary for data-driven LQR, although further analysis is required to prove this claim.

Although Theorem \ref{ch4:t:LQinform} gives necessary and sufficient conditions under which the data are informative for optimal linear quadratic regulation, it might not be directly clear how these conditions can be verified given the input-state data. Therefore, in what follows we rephrase the conditions of Theorem \ref{ch4:t:LQinform} in terms of the data matrices $X$ and $U_-$.

\begin{theorem}[Alternative conditions for informativity \cite{vanWaarde2020}]
	\label{ch4:t:LQinform2}
	Let $Q \geq 0$ and $R >0$. Then the data $(U_-,X)$ are informative for optimal linear quadratic regulation if and only if at least one of the following two conditions hold:
	\begin{enumerate}
		\item\label{ch4:cond:a2} The data $(U_-,X)$ are informative for identification, equivalently, there exists $\begin{bmatrix} V_1 & V_2\end{bmatrix}$ such that 
\begin{equation}\label{ch2:eq:V1 V2} \begin{bmatrix} X_- \\ U_- \end{bmatrix} \begin{bmatrix} 	V_1 & V_2	\end{bmatrix} = \begin{bmatrix} I_n & 0 \\ 0& I_m\end{bmatrix}. \end{equation}		
Moreover, the optimal LQR problem is solvable for $(A_s,B_s,Q,R)$, where $A_s=X_+V_1$ and $B_s=X_+V_2$.
		\item\label{ch4:cond:b2} There exists $\Theta \in \mathbb{R}^{T \times n}$ such that $X_- \Theta = (X_- \Theta)^\top$, 	$U_- \Theta = 0 $, 
		\begin{equation}
		\label{ch4:e:LMI/E/K/Q}
		\begin{bmatrix}
		X_- \Theta & X_+ \Theta \\ \Theta^\top X_+^\top & X_- \Theta
		\end{bmatrix} > 0.
		\end{equation}
		and $ QX_+\Theta = 0$. 
	\end{enumerate}
\end{theorem}

It is also possible to directly compute the optimal LQR feedback gain $K$ from the given data. 
%
Indeed, the following theorem asserts that $P^+$ as in Lemma~\ref{ch4:l: same P works for all} can be found as the unique solution to an optimization problem involving only the data. Furthermore, the optimal feedback gain $K$ can subsequently be found by solving a set of linear equations. In the sequel, for a given square matrix $M$, $\trace(M)$ will denote the trace of $M$.

\begin{theorem}[A semi-definite programming approach \cite{vanWaarde2020}]
	\label{ch4:t:LQgaindata}
Let $Q \geq 0$ and $R > 0$. Suppose that the data $(U_-,X)$ are informative for optimal linear quadratic regulation. Consider the linear operator $P\mapsto\calL(P)$ defined by
$$
\mathcal{L}(P) := \xmt P\xm-\xpt P\xp -\xmt Q\xm-\umt R\um.
$$
Let $P^+$ be as in Lemma~\ref{ch4:l: same P works for all}. The following statements hold:
	\begin{enumerate}
		\item \label{ch4:semidefiniteprogram} The matrix $P^+$ is equal to the unique solution to the optimization problem
		\begin{align*}
		\text{ maximize } \: &\trace(P) \\
		\text{ subject to } \: &P  \geq 0 
		\,\,\text{ and }\,\,  \mathcal{L}(P) \leq 0.
		\end{align*} 
		\item There exists a right inverse $X_-^\sharp$ of $X_-$ such that
			\begin{align}
			\label{ch4:eq:1}
			\mathcal{L}(P^+) X_-^\sharp &= 0.
			\end{align}
		Moreover, if $X_-^\sharp$ satisfies \eqref{ch4:eq:1}, then the optimal feedback gain is given by $K = U_- X_-^\sharp$.
	\end{enumerate}
\end{theorem}
From a design viewpoint, the optimal feedback gain $K$ can be found in the following way. First solve the semidefinite program in Theorem \ref{ch4:t:LQgaindata}. Subsequently, compute a solution $X_-^\sharp$ to the linear equations $X_- X_-^\sharp = I$ and \eqref{ch4:eq:1}. Then, the optimal feedback gain is given by $K = U_- X_-^\sharp$.

\subsection{The problem of tracking and regulation}


Yet another important classical control design problem is the problem of tracking and regulation, also called the algebraic regulator problem, as studied, for example, in \cite{Davison1975,Francis1977,Francis1975,Isidori1990} and the textbooks \cite{Trentelman2001,Saberi2000}). This is the problem of finding a feedback controller (called a regulator) such that the output of the resulting controlled system tracks a given reference signal, regardless of the disturbance input entering the system and the initial state. The relevant reference signals and disturbances (such as step functions, ramps or sinusoids) are assumed to be solutions of a suitable autonomous linear system. Given a class of reference and disturbance signals, one first constructs a suitable autonomous system (called the exosystem) that has these reference and disturbance signals as solutions. Next, this exosystem is interconnected to the system to be controlled (called the endosystem) and the difference between the original system output and the reference signal is taken as output. Finally, a regulator should be designed to make the output of the interconnection converge to zero for all disturbances and initial states.

In a data-driven context, the true endosystem is assumed to be unknown, and no mathematical model is available. Instead, we collect data on the input, endosystem state, and exosystem state in the form of samples on a finite time-interval. Whereas the true endo-system is unknown, the exosystem is assumed to be known, since this system models the reference signals and possible disturbance inputs. Also, the matrices in the output equations are assumed to be known, since these specify the design specification (namely the output that should converge to zero) on the controlled system.  A given set of data will then be called informative for regulator design if the data contain sufficient information to design a single regulator for the entire family of systems that are consistent with this set of data. In this section we will study this data-driven regulator problem, and provide necessary and sufficient conditions for informativity for regulator design.

Consider a true, unknown, endosystem represented by
\begin{equation} \label{e:endo-general}
\bmx_2(t+1) = A_{2s}\bmx(t) + B_{2s}\bmu(t) + A_{3}\bmx_1(t).
\end{equation}
Here, $\bmx_2$ is the $n_2$-dimensional state, $\bmu$ the $m$-dimensional input, and $\bmx_1$ the $n_1$-dimensional state of the exosystem 
\begin{equation} \label{e:exo}
\bmx_1(t+1) = A_1\bmx_1(t).
\end{equation}
that generates all possible reference signals and disturbance inputs. The dimensions $n_1$, $n_2$ and $m$ are known, but the matrices $A_{2s}$ and $B_{2s}$ are unknown. Since $A_3$ specifies how the disturbances and reference signals enter the system, we assume that it is known. Also the exo-system matrix $A_1$ is known. The output to be regulated is specified by
\begin{equation} \label{e:outequation}
\bmz(t) = D_1 \bmx_1(t) + D_2 \bmx_2(t) + E \bmu(t),
\end{equation}
where the matrices $D_1,D_2$ and $E$ are known. By interconnecting the endosystem with the state feedback controller
\begin{equation} \label{e:feedback}
\bmu(t) = K_1 \bmx_1(t) + K_2 \bmx_2(t), 
\end{equation}
we obtain the controlled system 
$$
\begin{bmatrix} \bmx_1(t+1) \\ \bmx_2(t+1)
\end{bmatrix}
   = 
\begin{bmatrix} A_1   & 0 \\ A_{3} + B_2 K_1& A_{2s} + B_{2s} K_2
\end{bmatrix}
\begin{bmatrix} \bmx_1(t) \\ \bmx_2(t) 
\end{bmatrix},
$$
$$
 \bmz(t)  =  ~(D_1 + E K_1) \bmx_1(t) + (D_2  + E K_2)\bmx_2(t).
 $$
If $\bmz(t) \rightarrow 0$ as $t \rightarrow \infty$ for all initial states $\bmx_1(0)$ and $\bmx_2(0)$, we say that the controlled system is {\em output regulated}.  If $A_{2s} + B_{2s} K_2$ is a stable matrix we call the controlled system {\em endo-stable}. If the control law \eqref{e:feedback} makes the controlled system both output regulated and endo-stable, we call it a {\em regulator}.

Since we do not know the true endosystem \eqref{e:endo-general}, the design of a regulator can only be based on available data. These data are finite sequences of samples of $\bmx_1(t), \bmx_2(t)$ and $\bmu(t)$ on a given time interval $\{0, 1, \ldots, T \}$ given by 
\begin{align*}
U_-& := \bbm u(0) & u(1) & \cdots & u(T-1)\ebm, \\
X_{1-}  & := \bbm x_1(0) & x_1(1) & \cdots & x_1(T-1)\ebm, \\
X_2  & := \bbm x_2(0) & x_2(1) & \cdots & x_2(T)\ebm.
\end{align*}
An endosystem with (unknown) system matrices $(A_2,B_2)$ is called {\em consistent} with these data if $A_2$ and $B_2$ satisfy the equation
\vspace{-2mm}
\begin{equation}  \label{e:data}
X_{2+} = A_{2} X_{2-} + A_3 X_{1-} + B_{2} U_- , 
\end{equation}
where we denote
\begin{align*}
X_{2-}& := \bbm x_2(0) & x_2(1) & \cdots & x_2(T -1) \ebm, \\
X_{2+}& := \bbm x_2(1) & x_2(2) & \cdots & x_2(T) \ebm.
\end{align*}	
The set of all $(A_2,B_2)$ that are consistent with the data is denoted by $\Sigma_\calD$, i.e.,
\begin{equation} \label{e:SigmaD}
\Sigma_\calD := \left\{ (A_2,B_2) \mid  \mbox{\eqref{e:data} holds} \right\}.
\end{equation}
We assume that the true endosystem $(A_{2s},B_{2s})$ is in $\Sigma_\calD$, i.e. the true system is consistent with the data.  In general, the equation \eqref{e:data} does not specify the true system uniquely, and many endosystems $(A_2,B_2)$ may be consistent with the same data.

Now we turn to controller  design based on the data $(U_-,X_{1-},X_2)$. Since on the basis of the given data we can not distinguish between the true endosystem and any other endosystem consistent with these data, a controller will be a regulator for the true system only if it is a regulator for any endosystem with $(A_2,B_2)$ in $\Sigma_\calD$. If such regulator exists, we call the data {\em informative for regulator design}. More precisely:
\begin{definition}[Informativity for regulator design]\label{def:informativity-regulator}
	We say that the data $(U_-,X_{1-},X_2)$ are informative for regulator design if there exists $K_1$ and $K_2$ such that the control law $\bmu(t) = K_1 \bmx_1(t) + K_2 \bmx_2(t)$ is a regulator for any endosystem with $(A_2,B_2)$ in $\Sigma_\calD$. 
\end{definition}
In this subsection we will present necessary and sufficient conditions on the data $(U_-,X_{1-},X_2)$ to be informative for regulator design. Also, in case that these conditions are satisfied, we will explain how to compute a regulator using only these data. 

The following theorem gives necessary and sufficient conditions on the data to be  informative for regulator design, and explains how suitable regulators are computed using only these data.
\begin{theorem}[Conditions for informativity \cite{Trentelman22}] \label{th:mainresult}
Assume that the matrix $A_1$ is anti-stable~\footnote{Anti-stable means that all its eigenvalues $\lambda$ satisfy $|\lambda| \geq 1$}. Then the data 
$(U_-,X_{1-},X_2)$ are informative for regulator design if and only if at least one of the following two conditions hold~\footnote{We denote by $\im M$ the image of the matrix $M$.}:
\begin{enumerate}
\item
$X_{2-}$ has full row rank, and there exists a right-inverse $X_{2-}^\sharp$ of $X_{2-}$ such that $(X_{2+} - A_3 X_{1-}) X_{2-}^\sharp$ is stable and $D_2 +E U_-X_{2-}^\sharp=0$.  Moreover, $\im D_1 \subseteq \im E$. In this case, a regulator is found as follows: choose $K_1$ such that $D_1 + EK_1 =0$ and define $K_2 := U_- X_{2-}^\sharp$.
\item
$X_{2-}$ has full row rank and there exists a right-inverse $X_{2-}^\sharp$ of $X_{2-}$ such that $(X_{2+} - A_3 X_{1-}) X_{2-}^\sharp$ is stable. Moreover, there exists a solution $W$ to the linear equations
\begin{subequations} \label{e:W}
\begin{align}
X_{2-}W A_1 - (X_{2+} - A_3 X_{1-}) W = A_3 ,\\
D_1 +(D_2 X_{2-} + E U_-) W  = 0,
\end{align} 
\end{subequations}
In this case, a regulator is found as follows: choose $K_1 := U_- (I - X_{2-}^\sharp X_{2-}) W$ and $K_2 := U_- X_{2-}^\sharp$.
\end{enumerate}
\end{theorem}
This theorem can be applied as follows. What we know about the system are the system matrices $A_1, A_3, D_1,D_2$ and $E$ and the data $(U_-,X_{1-},X_2)$.
The aim is to use this knowledge to compute a {\em single} regulator $(K_1,K_2)$ that works {\em for all} endosystems $(A_2,B_2)$ in the set $\Sigma_D$ defined by \eqref{e:SigmaD}.

In order to check the existence of such regulator, we verify  the two conditions 1) and 2) in Theorem \ref{th:mainresult}. If neither of the two conditions holds, then the data are not informative.  On the other hand, if condition 1) holds then a regulator $(K_1,K_2)$ is computed as follows: 
\begin{itemize}
\item
find a right-inverse $X_{2-}^{\sharp}$ of $X_{2-}$ such that the matrix $(X_{2+} - A_3 X_{1-}) X_{2-}^\sharp$ is stable and $D_2 +E U_-X_{2-}^\sharp =0$,
\item
compute $K_1$ as a solution of $D_1 + EK_1 =0$,
\item
define $K_2 := U_- X_{2-}^\sharp$. 
\end{itemize}
If condition 2) holds then a regulator is computed as follows: 
\begin{itemize}
\item
find a right-inverse $X_{2-}^{\sharp}$ of $X_{2-}$ such that the matrix $(X_{2+} - A_3 X_{1-}) X_{2-}^\sharp$ is stable,
\item
find a solution $W$ of the data-driven regulator equations \eqref{e:W},
\item
define $K_1 := U_- (I - X_{2-}^\sharp X_{2-}) W$,
\item
define $K_2 := U_- X_{2-}^\sharp$.
\end{itemize}

Although Theorem \ref{th:mainresult} gives a characterization of all data that are informative for regulator design and gives a method to design a suitable regulator, the procedure to compute
this regulator is not entirely satisfactory. Indeed, in the case that condition 2) holds it is not
clear how to find a right inverse of $X_{2-}$  such that $(X_{2+} - A_3 X_{1-}) X_{2-}^\sharp$ is stable. In the case of condition 1), the additional constraint $D_2 +E U_-X_{2-}^\sharp =0$ needs to be satisfied. In general, $X_{2-}$ has many right inverses, and $(X_{2+} - A_3 X_{1-}) X_{2-}^\sharp$ can be stable, with or without $D_2 +E U_-X_{2-}^\sharp =0$, depending on the choice of the particular right inverse $X_{2-}^\sharp$.
To deal with this problem and to solve the problem of regulator design, the problem of finding a suitable right inverse can be reformulated in terms of feasibility of an LMI, drawing some inspiration from Theorem~\ref{ch3:t:lmistab}.

\begin{theorem}[An LMI approach \cite{Trentelman22}] \label{th:regulator-with-LMI's}
Let $(U_-,X_{1-},X_2)$ be given data. Then the following hold:
\begin{enumerate}
\item
$X_{2-}$ has full row rank and has a right inverse $X_{2-}^{\sharp}$ such that  $(X_{2+} - A_3 X_{1-}) X_{2-}^\sharp$ is stable if and only if there exists a matrix $\Theta \in \mathbb{R}^{T \times n}$ such that 
\begin{equation} \label{e:regulator-sym}
X_{2-} \Theta = (X_{2-} \Theta)^\top 
\end{equation}
and
\begin{equation} \label{e:regulator-LMI}
 \bbm X_{2-} \Theta & (X_{2+} - A_3 X_{1-})\Theta \\
                                                                          \Theta^\top (X_{2+} - A_3 X_{1-})^\top & X_{2-} \Theta \ebm > 0.
\end{equation}    
\item $X_{2-}$ has full row rank and has a right inverse $X_{2-}^{\sharp}$ such that  $(X_{2+} - A_3 X_{1-}) X_{2-}^\sharp$ is stable with, in addition, $D_2 +E U_-X_{2-}^\sharp =0$ if and only if there exists a solution $\Theta \in \mathbb{R}^{T \times n}$ of \eqref{e:regulator-sym} and \eqref{e:regulator-LMI} that satisfies the linear equation $$(D_2X_{2-} + EU_-)\Theta = 0.$$     
\end{enumerate}
In both cases, a suitable right-inverse is given by $X_{2-}^{\sharp}: =   \Theta (X_{2-}\Theta)^{-1}$.                                                 
\end{theorem}

It is also possible to consider the situation that, in addition to $A_2$ and $B_2$, the matrix $A_3$ (representing how the exosignal $\bmx_1$ enters the endosystem) is unknown. In that case, the set all endosystems consistent with the data $(U_-,X_2,X_-)$ is given by:
\[
\Sigma_{D} = \{(A_2,B_2,A_3) \mid X_{2+} = A_2 X_{2-} + B_2U_- + A_3 X_{1-} \}.
\]
The data are then called informative for regulator design if there exists a single regulator $u = K_1 x_1 + K_2 x_2$ for all endosystems in $\Sigma_D$.  
The analogue of Theorem \ref{th:mainresult} for this situation is as follows. 
\begin{theorem}[Conditions for informativity \cite{Trentelman22}] \label{unknown A3}
Assume that the matrix $A_1$ is anti-stable. Then the data 
$(U_-,X_{1-},X_2)$ are informative for regulator design if and only if at least one of the following two conditions hold:
\begin{enumerate}
\item
$X_{2-}$ has full row rank, and there exists a right-inverse $X_{2-}^\sharp$ of $X_{2-}$ such that $X_{1-}X_{2-}^\sharp=0$, $(X_{2+} - A_3 X_{1-}) X_{2-}^\sharp$ is stable and $D_2 +E U_-X_{2-}^\sharp=0$.  Moreover, $\im D_1 \subseteq \im E$. In this case, a regulator is found as follows: choose $K_1$ such that $D_1 + EK_1 =0$ and define $K_2 := U_- X_{2-}^\sharp$.
\item
$X_{2-}$ has full row rank and there exists a right-inverse $X_{2-}^\sharp$ of $X_{2-}$ such that $X_{1-}X_{2-}^\sharp =0$ and $(X_{2+} - A_3 X_{1-}) X_{2-}^\sharp$ is stable. Moreover, there exists a solution $W$ to the linear equations
\begin{subequations} \label{e:Wnew}
\begin{align}
X_{2-}W A_1 - X_{2+} W = 0 ,\\
X_{1-} W = I, \\
D_1 +(D_2 X_{2-} + E U_-) W  = 0.
\end{align} 
\end{subequations}
In this case, a regulator is found as follows: choose $K_1 := U_- (I - X_{2-}^\sharp X_{2-}) W$ and $K_2 := U_- X_{2-}^\sharp$.
\end{enumerate}
\end{theorem}
Note that, as expected, $A_3$ no longer appears in the equations (it is unknown). In both cases, the formulas for $K_1$ and $K_2$ are the same as in Theorem \ref{th:mainresult}. 

Finally, we illustrate the application of Theorem \ref{th:mainresult} in the following example.
\begin{example}[Illustration of the theory] \label{ex:nonid}
Consider the two-dimensional endosystem
$$
\bmx_2(t+1) = A_{2s}\bmx_2(t) + B_{2s}\bmu(t) + \begin{bmatrix} 0 \\ 1 \end{bmatrix} \bmd(t),
$$
where $A_{2s}$ and $B_{2s}$ are unknown $2\times 2$  and $2 \times 1$ matrices, respectively. Let 
$\bmx_2 = \begin{bmatrix} \bmx_{21} & \bmx_{22} \end{bmatrix}^T$. The disturbance input $\bmd$ is assumed to be a constant signal with finite amplitude, so is generated by $\bmd(t+1) = \bmd(t)$. We want to design a regulator so that $2\bmx_{21} + \frac{1}{2} \bmx_{22}$ tracks a given reference signal.
 In this example, the reference signals $\bmr$ are assumed to be generated by a given autonomous linear system with state space dimension, say, $n_1$. Its representation will be irrelevant here. The total exosystem will then have state space dimension $n_1 + 1$, and our output equation is given by
$
\bmz(t) = D_1 \bmx_1(t) + D_2 \bmx_2(t) + E \bmu(t),
$
with $D_1$ a $1 \times (n_1 + 1)$ matrix such that $D_1\bmx_1 = - \bmr$ and $D_2 = \bbm 2 & \frac{1}{2} \ebm$. We take $E = 2$. Also note that 
$
A_3 = {\small \bbm  0_{1 \times n_1}  &  0 \\
                    0_{1 \times n_1}  & 1   \ebm.}
$
Here, $0_{1 \times n_1}$ denotes $1 \times n_1$ zero matrix.  Suppose that $T = 2$ and assume we have the following data:
$$
U_- = \bbm -1 & -1 \ebm, ~D_-  =  \bbm  1 &  1    \ebm,~
X_2   = \bbm 1 & \frac{1}{2} & -\frac{1}{4} \\
                        0  &      2         &   \frac{5}{2}   \ebm. 
$$
These data can be seen to be generated by the true endosystem
$
A_{2s} = {\small \bbm 2 & ~ \frac{1}{8} \vspace{.6mm}\\ 
                 4          & ~ \frac{5}{4} \ebm,} ~ 
B_{2s} = {\small \bbm \frac{3}{2} \vspace{.6mm}\\ 3 \ebm}.
$
We now check condition 1) of Theorem \ref{th:mainresult}. First note that, indeed, $\im D_1 \subseteq \im E$.  Also, $X_{2-}$ is non-singular and 
$
(X_{2+} - A_3 X_{1-}) X_{2-}^{-1} = {\small \bbm  \frac{1}{2}   & -\frac{1}{4} \\
                                                                    1                &  \frac{1}{2} \ebm}.
$
This matrix has eigenvalues $\frac{1}{2} \pm \frac{1}{2}i$, so is stable. Finally,
$D_2 +E U_-X_{2-}^{-1} =0$. According to Theorem \ref{th:mainresult}, a regulator for all endosystems consistentwith the given data is given by
\begin{equation} \label{e:reg1}
K_2 = U_1 X_{2-}^{-1} = \bbm -1 & -\frac{1}{4} \ebm,~ K_1 = - \frac{1}{2} D_1.
\end{equation}
It can be verified that the set of endosystems consistentwith our data is equal to the affine set
$$
\Sigma_{\calD} = \{ \left( \bbm a & \frac{1}{4} a - \frac{3}{8} \vspace{.6mm}\\
                                        b  & \frac{1}{4} b + \frac{1}{4} \ebm, 
                                        \bbm a - \frac{1}{2}  \vspace{.6mm} \\ b -1 \ebm \right )\mid a,b \in \mathbb{R} \}.
$$
The controller given by \eqref{e:reg1} is a regulator for all these endosystems.
\end{example}

Data-driven regulator design has also been studied in \cite{deCarolis2018} and \cite{Carnevale2017}. The perspective of these contributions is however quite different from the one discussed in this subsection.  We also mention alternative methods that deal with tracking objectives, such as iterative feedback tuning (IFT) and virtual reference feedback tuning (VRFT) as developed in \cite{Hjalmarsson1998} and \cite{Campi2002}, respectively. Also these methods do not address the classical regulator problem, and are very different from the work discussed here.

\section{Analysis and control using noisy input-state data}

So far, we have focused on analysis and design of input-state systems using \emph{exact data}. In this section, we shift our attention to input-state systems with noise. We will first introduce the model class that we will be using, and discuss the assumptions that will be made on the noise samples.

Suppose that the unknown, true system is given by 
\begin{equation}
    \label{ch2:eq:system1}
    \bmx(t+1) = A_s \bmx(t) + B_s \bmu(t) +\bmw(t),
 \end{equation}
where $\bmx \in \mathbb{R}^n$ is the state, $\bmu \in \mathbb{R}^m$ is the control input and $\bmw \in \mathbb{R}^n$ is an unknown noise term. The matrices $A_s \in \mathbb{R}^{n \times n}$ and $B_s \in \mathbb{R}^{n \times m}$ denote the unknown state and input matrices.
We embed this unknown system into the model class $\mathcal{M}$ of all  input-state systems with unknown process noise, with fixed dimensions $n$ and $m$, of the form
\begin{equation} \label{ch2:eq:input state with noise}
\bmx(t+1) = A\bmx(t) + B\bmu(t) +\bmw(t).
\end{equation}
Suppose that we obtain input-state data from the true system \eqref{ch2:eq:system1}. These data are given in the matrices
\begin{align*}
U_- &= \begin{bmatrix}
u(0) & u(1) & \cdots & u(T-1)
\end{bmatrix}, \\
X &= \begin{bmatrix}
x(0) & x(1) & \cdots & x(T)
\end{bmatrix}.
\end{align*}
We denote the submatrix of $X$ consisting of its first (respectively, last) $T$ columns by $X_-$ (respectively, $X_+$). The noise $\bmw$ is unknown, so $w(0),w(1),\dots,w(T-1)$ are not measured, and therefore are not part of the data. 
We do however assume that we have the following information on the noise during the data sampling period.
\begin{assumption}[Noise model] \label{ch2:assumption on noise samples}
The unknown noise samples $w(0),w(1),\dots,w(T-1)$, collected in the matrix 
$$
W_- := \bbm w(0) ~ w(1) ~ \cdots ~ w(T-1) \ebm,
$$
satisfy the \emph{quadratic matrix inequality}
\begin{equation} 
    \label{ch2:asnoise}
    \begin{bmatrix}
    I \\ W_-^\top 
    \end{bmatrix}^\top 
    \Phi
    \begin{bmatrix}
    I \\ W_-^\top 
    \end{bmatrix} \geq 0,
\end{equation}
where $\Phi \in \S{n + T}$ is a given partitioned matrix 
\begin{equation} \label{ch3:eq:Phi}
\Phi = \bbm \Phi_{11}  & \Phi_{12} \\ \Phi_{21} & \Phi_{22} \ebm,
\end{equation}
with $\Phi_{11} \in \S{n}$, $\Phi_{12} \in \mathbb{R}^{n \times T}$, $\Phi_{21} = \Phi_{12}^\top$ and $\Phi_{22} \in \S{T}$, and where we assume that 
$\Phi \in \bpi_{n,T}$ (as defined in \eqref{Piqr} of the Sidebar `Quadratic matrix inequalities'). 
\end{assumption}
In other words, the data $\calD$ consist of the measurements $(U_-,X)$ together with the information that the noise on the sampling interval $\{0, \ldots, T\}$ satisfies the inequality \eqref{ch2:asnoise} for a partitioned matrix $\Phi \in \bpi_{n,T}$.

\begin{sidebar}{Quadratic matrix inequalities}
\label{sidebar:QMIs}

\sdbarinitial{A} great deal can be said about sets defined in terms of quadratic matrix inequalities (QMIs). Such sets play an important role in robust control, where they are used to describe parameter uncertainty \cite{Petersen1987,Petersen1986,Megretski1997,Iwasaki1998,Scherer2001,Scherer2005,Scherer2006}. Also in the context of this paper they are important, since they describe sets of systems consistent with the data. We will briefly discuss sets of the form
\begin{sequation} \label{ch0:e:Zr}
\calZ_{r}(\Pi):=\left\{ Z\in\R^{r\times q} \mid \bbm I_q\\Z\ebm^\top\Pi\bbm I_q\\Z\ebm\geq 0\right\},
\end{sequation}
for $\Pi\in\S{q+r}$, and we refer to \cite{vanWaarde2022} for more details. In what follows, we assume that $\Pi \in\S{q+r}$ is partitioned as
$$
\Pi = \mpi,
$$
where $\Pi_{11} \in \S{q}$, $\Pi_{12} = \Pi_{21}^\top \in \mathbb{R}^{q\times r}$ and $\Pi_{22} \in \S{r}$. The very first question one may ask is: under what conditions on $\Pi$ is the set $\calZ_{r}(\Pi)$ nonempty? An immediate necessary condition is that $\Pi$ must have at least $q$ nonnegative eigenvalues. However this is not sufficient in general. It turns out that for particular matrices $\Pi$, a Schur complement argument on the matrix $\Pi$ leads to a simple characterization of nonemptiness of the set $\calZ_{r}(\Pi)$. 
Specifically, suppose that $\Pi_{22}\leq 0$ and $\ker\Pi_{22}\subseteq\ker\Pi_{12}$. Since the latter condition is equivalent to $\Pi_{12}\Pi_{22}\Pi_{22}\gi=\Pi_{12}$, we have that
\begin{equation*}
\mpi\!=\!\bbm I_q & \Pi_{12}\Pi_{22}\gi\\0&I_r\ebm\!\!\bbm \Pi\schur\Pi_{22}& 0 \\0 & \Pi_{22}\ebm\!\!\bbm I_q & 0\\\Pi_{22}\gi\Pi_{21} & I_r\ebm\!,
\end{equation*}
where $\Pi\schur\Pi_{22}:=\Pi_{11}-\Pi_{12}\Pi_{22}\gi\Pi_{21}$ is the (generalized) Schur complement of $\Pi$ with respect to $\Pi_{22}$. This can be used to prove the following conditions for nonemptiness of $\calZ_{r}(\Pi)$.
\begin{theorem}[Nonemptiness of $\calZ_{r}(\Pi)$ \cite{vanWaarde2022}] \label{ch0:t:Z-r nonempty}
Let $\Pi \in \mathbb{S}^{q+r}$ and assume that $\ker \Pi_{22} \subseteq \ker\Pi_{12}$. Then $\calZ_r(\Pi)$ is nonempty if $\Pi \schur \Pi_{22} \geq 0$. Moreover, under the assumption that $\Pi_{22} \leq 0$  we have that $\calZ_r(\Pi)$ is nonempty if and only if $\Pi \schur \Pi_{22} \geq 0$.
\end{theorem}
Motivated by this, we define the set
\begin{sequation}
\label{Piqr}
\bpi_{q,r}\!:=\!\set{\Pi\! \in\!\S{q+r}}{\!\Pi_{22}\leq 0, \Pi\schur\Pi_{22} \geq 0 \text{ and }\ker\Pi_{22}\!\subseteq\!\ker\Pi_{12}}.
\end{sequation}
Next, for $\Pi\in\bpi_{q,r}$ we will investigate properties of the sets $\calZ_{r}(\Pi)$ and the following closely related set 
\begin{sequation} \label{ch0:e:Zr+}
\calZ_r^+(\Pi) := \left\{ Z \in \mathbb{R}^{r\times q} \mid \begin{bmatrix}
I_q \\ Z
\end{bmatrix}^\top \Pi \begin{bmatrix}
I_q \\ Z
\end{bmatrix} > 0 \right\},
\end{sequation}
involving a strict inequality.

\section{Matrix versions of Yakubovich's S-lemma}
\label{s:FinslerandSlemma}
Yakubovich' S-lemma \cite{Yakubovich1977} is a classical result with a wide range of applications, most notably the problem of absolute stability of Lur'e systems. Roughly speaking, this result says that one quadratic inequality implies another one if and only if a certain \emph{linear matrix inequality} (LMI) is feasible. A seemingly difficult implication involving quadratic functions is thereby replaced by a convex problem which can be solved using computational tools such as Sedumi and Mosek. In this section we deal with matrix versions of the S-lemma, i.e., with the question under what conditions all solutions to one quadratic \emph{matrix} inequality also satisfy another QMI. In other words, we state necessary and sufficient conditions for the inclusion $\calZ_r(N) \subseteq \calZ_r(M)$, where $M,N\in \mathbb{S}^{q+r}$. We will also consider a similar inclusion with $Z_r^+(M)$ instead of $Z_r(M)$.  This leads to non-strict and strict versions of Yakubovich's S-lemma.

\subsection{Non-strict inequalities}
The following theorem states the matrix S-lemma for non-strict inequalities.
\begin{theorem}[Matrix S-lemma \cite{vanWaarde2022}]
\label{t:nonstrictS-lemma}
Let $M,N \in \mathbb{S}^{q+r}$. If there exists a real $\alpha \geq 0$ such that $M - \alpha N \geq 0$, then $\calZ_r(N) \subseteq \calZ_r(M)$. Next, assume that $N \in \bpi_{q,r}$ and $N$ has at least one positive eigenvalue. Then $\calZ_r(N) \subseteq \calZ_r(M)$ if and only if there exists a real $\alpha \geq 0$ such that $M-\alpha N \geq 0$. 
\end{theorem}

Similar to the `standard' S-lemma we note that the matrix S-lemma requires $N$ to have at least one positive eigenvalue, an assumption known as the \emph{Slater condition}. It turns out, however, that under additional assumptions on $M$ and $N$, a result similar to Theorem~\ref{t:nonstrictS-lemma} holds when $N$ has \emph{no} positive eigenvalues. This leads to a matrix version of Finsler's lemma, which will not be further discussed here.

\subsection{A strict and non-strict inequality}

Next, we consider {\em strict} versions of the matrix S-lemma. This means that we consider the set $\calZ_r^+(M)$ instead of $\calZ_r(M)$, i.e., a strict inequality on the QMI induced by $M$. Note that in this case, the Slater condition on $N$ is not required.

\begin{theorem}[Strict matrix S-lemma \cite{vanWaarde2022}]
\label{t:strictS-lemmaN22}
Let $M,N \in \mathbb{S}^{q+r}$. If there exists a real $\alpha \geq 0$ such that $M - \alpha N > 0$, then $\calZ_r(N) \subseteq \calZ_r^+(M)$. Next, assume that $N \in \bpi_{q,r}$ and $N_{22} < 0$. Then $\calZ_r(N) \subseteq \calZ_r^+(M)$ if and only if there exists a real $\alpha \geq 0$ such that $M-\alpha N > 0$. 
\end{theorem}

It is also possible to proceed if $N_{22}$ is not necessarily negative definite, but an extra condition on $M$ holds. In that case two real numbers $\alpha \geq 0$ and $\beta > 0$ are required to arrive at a necessary and sufficient condition.

\begin{theorem}[Strict matrix S-lemma with $\alpha$ and $\beta$ \cite{vanWaarde2022}]
\label{c:combinedstrictS-lemmaFinslerslemma}
Let $M,N \in  \mathbb{S}^{q+r}$. Then $\calZ_r(N) \subseteq \calZ_r^+(M)$
if there exist scalars  $\alpha \geq 0$ and $\beta >0$ such that
\begin{sequation}
\label{ineqalphabeta}
M-\alpha N \geq \begin{bmatrix}
\beta I & 0 \\ 0 & 0
\end{bmatrix}.
\end{sequation} 
Next, assume that $N \in \bpi_{q,r}$ and $M_{22} \leq 0$. Then $\calZ_r(N) \subseteq \calZ_r^+(M)$ if and only if there exist $\alpha \geq 0$ and $\beta >0$ such that \eqref{ineqalphabeta} holds.
\end{theorem}

\end{sidebar}

Of course, an issue is whether the set of noise matrices $W_-$ defined by \eqref{ch2:asnoise} is nonempty. This is equivalent to the nonemptiness of the set $\calZ_T(\Phi)$, as defined in \eqref{ch0:e:Zr}. This issue is discussed in more detail in the Sidebar `Quadratic matrix inequalities'. Indeed, under the assumption $\Phi \in \bpi_{n,T}$ the set $\calZ_T(\Phi)$ is nonempty and  convex. Consequently then, the set of of noise matrices $W_-$ satisfying \eqref{ch2:asnoise} is nonempty and convex.

In order to make the above quadratic inequality constraint on the matrix of noise samples more concrete, we will now look at a number of special cases.
\begin{enumerate}
\item In the special case $\Phi_{12} = 0$ and $\Phi_{22} = -I$, the quadratic inequality \eqref{ch2:asnoise} reduces to 
\begin{equation}
    \label{ch2:redasnoise}
    W_- W_-^\top = \sum_{t=0}^{T-1} w(t) w(t)^\top \leq \Phi_{11}.
\end{equation}
The inequality \eqref{ch2:redasnoise} can be interpretated as saying that the energy of $\bmw$ has a given upper bound on the time interval $\{0,\ldots, T-1\}$.
\item
Norm bounds on the individual noise samples $w(t)$ also give rise to bounds of the form \eqref{ch2:asnoise}, although this 
does introduce some conservatism in general.
Indeed, note that for all $t$ the pointwise norm bound $\norm{w(t)}_2^2 \leq \epsilon$ is equivalent to  the matrix inequality $w(t) w(t)^\top \leq \epsilon I$.  As such, the bound \eqref{ch2:redasnoise} is satisfied for $\Phi_{11} = T \epsilon I$.
\item
In some cases, we may know a priori that the noise $\bmw$ does not directly affect the entire state-space, but is contained in a subspace, say $\im E$, with $E$ a known $n \times d$ matrix.  This prior knowledge can be captured by the noise model in Assumption \ref{ch2:assumption on noise samples}. Indeed, suppose that $w(t) = E \hat{w}(t)$ for all $t = 0,1,2\dots,T-1$, where $\hat{w}(t) \in \mathbb{R}^d$ and $E \in \mathbb{R}^{n \times d}$ is a given matrix of full column rank. The matrix $
\hat{W}_- = \begin{bmatrix}
\hat{w}(0) & \hat{w}(1) & \cdots & \hat{w}(T-1)
\end{bmatrix}$ captures the noise. As before, $\hat{W}_-$ is unknown but is assumed to satisfy $\hat{W}_-^\top \in \calZ_T(\hat{\Phi})$, where $\hat{\Phi} \in \bpi_{d,T}$ is such that $\hat{\Phi}_{22} < 0$. It can then be shown that $W_- = E\hat{W}_-$ for some $\hat{W}_-^\top \in \calZ_T(\hat{\Phi})$ if and only if $W_-^\top \in \calZ_T(\Phi)$, where
\begin{equation}
\label{PhiwithE}
\Phi := \begin{bmatrix}
E \hat{\Phi}_{11} E^\top & E \hat{\Phi}_{12} \\
\hat{\Phi}_{21} E^\top & \hat{\Phi}_{22} 
\end{bmatrix} \in \bpi_{n,T}.
\end{equation}
The conclusion is that Assumption \ref{ch2:assumption on noise samples} also covers the case in which the noise is constrained to a known subspace, which is captured by the noise bound \eqref{ch2:asnoise} with $\Phi$ in \eqref{PhiwithE}.
\item As shown in \cite[Section 5.4]{vanWaarde2022}, these noise models can also be applied in settings of Gaussian noise. To be precise, such sets can be employed as confidence intervals corresponding to a given probability.
\end{enumerate}

\subsection{Quadratic stabilization}

In this subsection we will take a look at the problem of quadratic stabilization using noisy input-state data. Quadratic stabilization means that all systems in the set $\Sigma_{\calD}$ of systems consistent with the data can be stabilized by the same state feedback gain, with \emph{a common Lyapunov function} for all closed loop systems. In particular then, this feedback gain will stabilize the unknown system.
Conditions for the existence of such feedback gain will be in terms of feasibility of certain linear matrix inequalities involving the data $(X,U_-)$ and the (known) matrix $\Phi$ representing the quadratic inequality constraint on the matrix of noise samples. In addition, the controller gain will be computed in terms of solution to these linear matrix inequalities.

As explained before, we have access to the input-state data $\mathcal{D} = (U_-,X)$. The possible matrices $W_-$ of noise samples satisfy the quadratic inequality \eqref{ch2:asnoise} for a given matrix $\Phi \in  \bpi_{n,T}$. This means that the set $\Sigma_{\calD}$ of all systems consistent with the data is equal to the set of all systems $(A,B)$ satisfying 
\begin{equation}
    \label{ch3:eq:dataeq}
    X_+ = A X_- + B U_- + W_-
\end{equation}
for some $W_-$ satisfying \eqref{ch2:asnoise}, i.e.,
\begin{equation} \label{ch3:def:SigmaD}
 \Sigma_\calD = \{ (A,B) \mid \eqref{ch3:eq:dataeq} \text{ holds for some } W_- \text{ satisfying } \eqref{ch2:asnoise} \}.
\end{equation}

\begin{definition}[Informativity for quadratic stabilization]
\label{ch3:def:informativity quad stabilization}
The data $(U_-,X)$ are called \emph{informative for quadratic stabilization} if there exists a feedback gain $K \in \mathbb{R}^{m \times n}$ and a matrix $P \in \S{n}$ such that $P > 0$ and
\begin{equation}
    \label{ch3:eq:closed loop lyapunovineq}
P - (A+BK) P (A+BK)^\top > 0
\end{equation}
for all $(A,B) \in \Sigma_{\calD}$.
\end{definition}
We are interested in \emph{quadratic stabilization} in the sense that we ask for a \emph{common} Lyapunov matrix $P$ for all $(A,B) \in \Sigma_{\calD}$.  Note that $P>0$ satisfies \eqref{ch3:eq:closed loop lyapunovineq} if and only $Q:= P^{-1}$ satisfies $Q - (A +BK)^\top Q (A +BK) >0$, which expresses that $V(x) = x^T Q x$ is a Lyapunov function for the system $\bmx(t + 1) = (A + BK)\bmx(t)$.

Definition~\ref{ch3:def:informativity quad stabilization} leads to two natural problems. First, we are interested in the question under which conditions the data are informative. The second problem is a design issue: we are interested in procedures to come up with a feedback gain that stabilizes all systems in $\Sigma_{\calD}$. By making use of the linear equation \eqref{ch3:eq:dataeq} and the assumption \eqref{ch2:asnoise} on the noise, it is straightforward to see that $(A,B) \in \Sigma_{\calD}$ if and only if 
\begin{equation}
    \label{ch3:eq:ineqAB}
    \begin{bmatrix}
    I \\ A^\top \\ B^\top 
    \end{bmatrix}^\top 
    \begin{bmatrix}
    I & X_+ \\ 0 & -X_- \\ 0 & -U_-
    \end{bmatrix}
    \begin{bmatrix}
    \Phi_{11} & \Phi_{12} \\
    \Phi_{21} & \Phi_{22}
    \end{bmatrix}
    \begin{bmatrix}
    I & X_+ \\ 0 & -X_- \\ 0 & -U_-
    \end{bmatrix}^\top
    \begin{bmatrix}
    I \\ A^\top \\ B^\top 
    \end{bmatrix} \geq 0.
\end{equation}
Next, suppose that we fix a Lyapunov matrix $P> 0$ and a feedback gain $K$. The inequality \eqref{ch3:eq:closed loop lyapunovineq} is equivalent to 
\begin{equation}
\label{ch3:eq:ineqABPK}
    \begin{bmatrix}
    I \\ A^\top \\ B^\top 
    \end{bmatrix}^\top
    \begin{bmatrix}
    P & 0 & 0 \\
    0 & -P & -PK^\top \\
    0 & -KP & -KPK^\top
    \end{bmatrix}
    \begin{bmatrix}
    I \\ A^\top \\ B^\top 
    \end{bmatrix} > 0,
\end{equation}
which is also a quadratic matrix inequality in $A$ and $B$. Therefore, finding conditions for quadratic stabilization amounts to finding conditions under which the quadratic matrix inequality \eqref{ch3:eq:ineqABPK} holds for all $(A,B)$ satisfying the quadratic matrix inequality \eqref{ch3:eq:ineqAB}. Let $N$ be defined by 
\begin{equation}
\label{ch2:Nstab}
N \!=\! \begin{pmat}[{|}]
N_{11} & N_{12} \cr\- N_{12}^\top & N_{22} \cr
\end{pmat} 
\!:=\! \begin{pmat}[{.}]
    I & X_+ \cr\- 0 & -X_- \cr 0 & -U_- \cr
    \end{pmat}
    \!\!
    \begin{bmatrix}
    \Phi_{11} & \Phi_{12} \\
    \Phi_{21} & \Phi_{22}
    \end{bmatrix}
    \!\!
    \begin{pmat}[{.}]
    I & X_+ \cr\- 0 & -X_- \cr 0 & -U_- \cr
    \end{pmat}^\top\!,
\end{equation}
 and define
\begin{align}
    M :=& 
\begin{pmat}[{|}]
M_{11} & M_{12} \cr\-
M_{12}^\top & M_{22} \cr
\end{pmat} := \begin{pmat}[{|.}]
    P & 0 & 0 \cr\-
    0 & -P & -PK^\top \cr
    0 & -KP & -KPK^\top \cr
    \end{pmat} \label{ch3:Mstab}.
\end{align}
Then we need to find conditions on the data such that there exist $P >0$ and $K$ such that the inclusion 
\begin{equation} \label{ch3:eq:QMI inclusion}
\calZ_{n +m}(N) \subseteq \calZ_{n +m}^+(M)
\end{equation}
holds. In order to find such conditions, we can use a matrix version of the S-lemma, as reported in Theorem~\ref{c:combinedstrictS-lemmaFinslerslemma} of the Sidebar `Quadratic matrix inequalities'. It is straightforward to verify the assumptions of this result, i.e., $M_{22} \leq 0$ and $N \in \bpi_{n,T}$.
Theorem~\ref{c:combinedstrictS-lemmaFinslerslemma} then asserts that \eqref{ch3:eq:QMI inclusion} holds if and only if there exist scalars $\alpha \geq 0$ and $\beta > 0$ such that
\begin{equation}
\label{ch3:ineqMNab}
M - \alpha N \geq \begin{bmatrix}
\beta I & 0  & 0 \\ 0 & 0 & 0 \\ 0  &  0  &  0
\end{bmatrix}.
\end{equation}
From a design point of view, the matrices $P$ and $K$ that appear in $M$ are not given. However, the idea is now to \emph{compute} matrices $P$, $K$ and scalars $\alpha$ and $\beta$ such that \eqref{ch3:ineqMNab} holds. In fact, by the above discussion, the data $(U_-,X)$ are informative for quadratic stabilization \emph{if and only if} there exist $P \in \S{n}$, $P > 0$, $K \in \mathbb{R}^{m \times n}$ and two scalars $\alpha \geq 0$ and $\beta > 0$ such that \eqref{ch3:ineqMNab} holds. We note that \eqref{ch3:ineqMNab} (in particular, $M$) is not linear in $P$ and $K$. Nonetheless, by a rather standard change of variables and a Schur complement argument, we can transform \eqref{ch3:ineqMNab} into a linear matrix inequality. 
Moreover, it turns out that the scalar $\alpha$ is necessarily positive. By a scaling argument then, it can be chosen to be equal to 1.  We summarize our result in the following theorem.

\begin{theorem}[Informativity for quadratic stabilization \cite{vanWaarde2022}]
\label{ch3:th:theoremstab}
Suppose that the data $(U_-,X)$ are collected from system \eqref{ch2:eq:system1} with noise as in Assumption~\ref{ch2:assumption on noise samples}. The data $(U_-,X)$ are informative for quadratic stabilization if and only if there exist $P \in \S{n}$, $P > 0$,   $L \in \mathbb{R}^{m \times n}$ and a scalar 
$\beta > 0$ satisfying 
\begin{equation}
\begin{pmat}[{..|}]
    P-\beta I & 0 & 0 & 0 \cr
    0 & -P & -L^\top & 0 \cr
    0 & -L & 0 & L \cr\-
    0 & 0 & L^\top & P \cr
    \end{pmat} - \begin{pmat}[{|}]
    N & 0 \cr\- 0 & 0 \cr
\end{pmat} \geq 0. \label{ch3:eq:LMIstab}
\end{equation}
Moreover, if $P$ and $L$ satisfy \eqref{ch3:eq:LMIstab} then $K := L P^{-1}$ is a stabilizing feedback gain for all $(A,B) \in \Sigma_{\calD}$.
\end{theorem}

Theorem~\ref{ch3:th:theoremstab} provides a necessary and sufficient condition under which quadratically stabilizing controllers can be obtained from noisy data. The theorem leads to an effective design procedure for obtaining stabilizing controllers directly from data. Indeed, the approach entails solving the linear matrix inequality \eqref{ch3:eq:LMIstab} for $P, L$ and $\beta$ and computing a controller as $K = LP^{-1}$. Below, we discuss some of the features of our control design procedure. 
\begin{enumerate}
    \item First of all, we stress that the procedure is \emph{non-conservative} since Theorem~\ref{ch3:th:theoremstab} provides a necessary and sufficient condition for obtaining quadratically stabilizing controllers from data. 
    \item The variables $P, L$ and $\beta$ are \emph{independent} of the time horizon $T$ of the experiment. In fact, note that $P \in \mathbb{R}^{n \times n}$, $L \in \mathbb{R}^{m \times n}$ and $\beta \in \mathbb{R}$. Also, the LMI \eqref{ch3:eq:LMIstab} is of dimension $(3n+m) \times (3n+m)$ and thus independent of $T$. This $T$-independent design method can play a crucial role in control design from larger data sets. We note that collections of big data sets are often unavoidable, for example because the signal-to-noise ratio is small, or because the data-generating system is large-scale. 
\end{enumerate}

We note that under the extra assumptions $\Phi_{22} < 0$ and  
\begin{equation}
    \label{ch3:fullrank}
\rank \begin{bmatrix}
X_- \\ U_-
\end{bmatrix} = n+m
\end{equation}
it is possible to prove a variant Theorem~\ref{ch3:th:theoremstab} in which the non-strict inequality is replaced by a strict inequality, and the term $-\beta I$ is removed. This can be done by invoking Theorem~\ref{t:strictS-lemmaN22} of Sidebar `Quadratic matrix inequalities', which is possible since the conditions $\Phi_{22} < 0$ and \eqref{ch3:fullrank} yield $N_{22} < 0$.  Thus we obtain the following theorem.
\begin{theorem}[Informativity via a strict inequality \cite{vanWaarde2022}] \label{ch3:th:stabilizationN22 <0}
Suppose that the data $(U_-,X)$ are collected from system \eqref{ch2:eq:system1} with noise as in Assumption~\ref{ch2:assumption on noise samples}. In addition, assume that $\Phi_{22} < 0$ and the rank condition \eqref{ch3:fullrank} holds. Then the data $(U_-,X)$ are informative for quadratic stabilization if and only if there exist $P \in \S{n}$, $P > 0$  and $L \in \mathbb{R}^{m \times n}$ satisfying 
\begin{equation}
\begin{pmat}[{..|}]
    P & 0 & 0 & 0 \cr
    0 & -P & -L^\top & 0 \cr
    0 & -L & 0 & L \cr\-
    0 & 0 & L^\top & P \cr
    \end{pmat} - \begin{pmat}[{|}]
    N & 0 \cr\- 0 & 0 \cr
\end{pmat} > 0. \label{ch3:eq:LMIstabN22}
\end{equation}
Moreover, if $P$ and $L$ satisfy \eqref{ch3:eq:LMIstabN22} then $K := L P^{-1}$ is a stabilizing feedback gain for all $(A,B) \in \Sigma_{\calD}$.
\end{theorem}


Assume now that $\Phi_{22} < 0$ in \eqref{ch3:eq:Phi}. Under this assumption, it turns out that if the data $(U_-,X)$ are informative for quadratic stabilization and if $K$ stabilizes all systems in $\Sigma_{\calD}$ with a common Lyapunov matrix $P >0$, then, in fact, $X_-$ must have full row rank, and $K$ must be of the form $K = U_- X_-^{\sharp}$ for some right inverse $X_-^{\sharp}$ of
$X_-$. Thus, the following theorem extends Lemma \ref{ch3:lemmaF0=0} to the noisy case.
\begin{theorem}[Necessary conditions for informativity \cite{vanWaarde2022b}] \label{ch3:th:generalization}
Suppose that the data $(U_-,X)$ are collected from system \eqref{ch2:eq:system1} with noise as in Assumption~\ref{ch2:assumption on noise samples}. In addition, assume $\Phi_{22} < 0$. Let the data $(U_-,X_-)$ be informative for quadratic stabilization and suppose that $P >0$ and $K$ are such that \eqref{ch3:eq:QMI inclusion} holds. Then 
\begin{equation}
    \label{ch5:IKXU}
\im \begin{bmatrix}
I \\ K
\end{bmatrix} \subseteq 
\im
\begin{bmatrix}
X_- \\ U_-
\end{bmatrix}.
\end{equation}
Consequently, $X_-$ has full row rank $n$ and there exists a right-inverse $X_-^{\sharp}$ of
$X_-$ such that $K = U_- X_-^{\sharp}$.
\end{theorem}

\subsection{Related conditions for quadratic stabilization}

Theorem~\ref{ch3:th:stabilizationN22 <0} gives a necessary and sufficient LMI condition under which all systems consistent with the data are quadratically stabilizable by a single feedback gain $K$. In this section we compare this result to other conditions within the literature on data-driven control. 

We begin with \cite[Thm. 6]{DePersis2020}. This result works under the assumption that \eqref{ch3:fullrank} holds, and $X_+$ has full row rank. Moreover, it is assumed that 
\begin{equation}
\label{asnoiseDePersis}
W_- W_-^\top \leq \gamma X_+ X_+^\top
\end{equation}
for some $\gamma > 0$. Under these assumptions, \cite[Thm. 6]{DePersis2020} states that the data $(U_-,X)$ are informative for quadratic stabilization if there exists a matrix $Q \in \mathbb{R}^{T \times n}$ and a scalar $\alpha > 0$ such that $X_- Q$ is symmetric and 
\begin{align}
\label{DePersis1}
\begin{bmatrix}
X_- Q - \alpha X_+ X_+^\top & X_+ Q \\ Q^\top X_+^\top & X_- Q
\end{bmatrix} >0, \:\: \begin{bmatrix}
I & Q \\ Q^\top & X_- Q
\end{bmatrix} > 0, \\
\label{DePersis2}
\frac{\alpha^2}{4+2\alpha} > \gamma.
\end{align} 
If $(Q,\alpha)$ solve \eqref{DePersis1},\eqref{DePersis2} then $K:= U_- Q (X_- Q)^{-1}$ quadratically stabilizes all systems in $\Sigma_{(U_-,X)}$. We note that the inequality \eqref{asnoiseDePersis} can be interpreted as a special case of \eqref{ch2:asnoise} with $\Phi_{11} = \gamma X_+ X_+^\top$, $\Phi_{12} = 0$ and $\Phi_{22} = -I$. 

Yet another condition for quadratic stabilization is given in \cite{Berberich2019c}. This paper works with a noise model that can be interpreted as the dual of \eqref{ch2:asnoise}. More precisely, it is assumed that 
\begin{equation}
\label{asnoiseBerberich}
\begin{bmatrix}
 W_- \\ I
\end{bmatrix}^\top \begin{bmatrix}
Q_w & S_w \\ S_w^\top & R_w
\end{bmatrix}
\begin{bmatrix}
 W_- \\ I
\end{bmatrix} \geq 0
\end{equation}
for known matrices $Q_w \in \mathbb{S}^{n}$, $S_w \in \mathbb{R}^{n \times T}$ and $R_w \in \mathbb{S}^{T}$ with $R_w > 0$. To make a meaningful comparison, we will assume the same bound as in \eqref{asnoiseDePersis}. This can also be stated equivalently in terms of the noise model \eqref{asnoiseBerberich} by choosing the specific matrices $Q_w = -(\gamma X_+ X_+^\top)^{-1}$, $S_w = 0$ and $R_w = I$. Then, the main result of \cite[Cor. 6, Rem. 7]{Berberich2019c} is that the data $(U_-,X)$ are informative for quadratic stabilization if there exist matrices $\mathcal{Y} \in \mathbb{S}^{n}$ and $M \in \mathbb{R}^{T \times n}$ satisfying 
\begin{align}
\label{Berberich1}
\begin{bmatrix}
-\mathcal{Y} & 0 & M^\top X_+^\top & M^\top \\
0 & Q_w & I & 0 \\
X_+ M & I & -\mathcal{Y} & 0 \\ 
M & 0 & 0 & -R_w^{-1}
\end{bmatrix} &< 0, \\
\label{Berberich2}
X_- M &= \mathcal{Y}.
\end{align}

We note that the conditions from \cite{DePersis2020} and \cite{Berberich2019c} are stated as sufficient conditions for quadratic stabilization. It is an interesting question whether these conditions are also \emph{necessary} for quadratic stabilization. Indeed, in this case, they would then be equivalent to those of Theorem~\ref{ch3:th:stabilizationN22 <0} (for the noise model in \eqref{asnoiseDePersis}). It turns out, however, that this is not the case.

To show this, one can consider the true system described by the matrices $A_s = 1$ and $B_s = 1$. Suppose that $T = 3$ and the noise matrix is given by
$$
W_- = \begin{bmatrix}
\frac{1}{2} & \frac{1}{2} & \frac{1}{2}
\end{bmatrix}.
$$
We collect the data samples 
\begin{align*}
X &= \begin{bmatrix}
0 & 0 & 1 & 0
\end{bmatrix}, \\
U_- &= \begin{bmatrix}
-\frac{1}{2} & \frac{1}{2} & -\frac{3}{2}
\end{bmatrix}.
\end{align*}
Throughout the example, we assume that we have access to the noise bound $W_- W_-^\top \leq 1$. Note that this bound is indeed satisfied, and that it can be captured using Assumption~\ref{ch2:assumption on noise samples} by the choices $\Phi_{11} = 1$, $\Phi_{12} = 0$ and $\Phi_{22} = -I$. We also note that this is equivalent to noise model \eqref{asnoiseBerberich} with $Q_w = -1$, $S_w = 0$ and $R_w = I$, and to noise model \eqref{asnoiseDePersis} with $\gamma = 1$. As such, we can compare the design methods reported in Theorem~\ref{ch3:th:stabilizationN22 <0} of this paper with the approaches in \cite[Cor. 6, Rem. 7]{Berberich2019c} and \cite[Thm. 6]{DePersis2020}.

For this example, it can be shown analytically \cite{vanWaarde2022b} that only the LMI condition of Theorem~\ref{ch3:th:stabilizationN22 <0} is feasible while those in \eqref{DePersis1}-\eqref{DePersis2} and \eqref{Berberich1}-\eqref{Berberich2} are not. At a high level, the reason for this is that the approach of \cite{DePersis2020} relies on a number of possibly conservative bounds, while the method from \cite{Berberich2019c} utilizes an overparameterization of the set of consistent systems.

\subsection{The $\mathcal{H}_\infty$ control problem}

The informativity framework also allows a treatment of the data-driven $\mathcal{H}_\infty$ control problem. This will be the topic of the current subsection. We first review some basic material that will be needed in order to formulate the problem.
 
Denote by $\ell_2^q(\mathbb{Z}_+)$ the linear space of all sequences $\bmv$ with $\bmv(t) \in \mathbb{R}^q$ and $t \in \mathbb{Z}_+$ such that $\sum_{t = 0}^{\infty} \|\bmv(t)\|^2 < \infty$. For any such sequence $\bmv$, define its $\ell_2$-norm as
$$
\|\bmv\|_2 := \left( \sum_{t = 0}^{\infty} \|\bmv(t)\|^2 \right)^{\frac{1}{2}}.
$$
Next, consider the discrete-time input-state-output system
\begin{equation} \label{ch5:eq:isosystem}
\begin{aligned}
    \bmx(t+1) &= A\bmx(t) + E\bmw(t),  \\
    \bmz(t) &= C\bmx(t) + D\bmw(t), 
\end{aligned} 
\end{equation}
with $\bmw(t) \in \mathbb{R}^q$ and $\bmz(t) \in \mathbb{R}^p$. Let its transfer matrix be denoted by $G(z)  := C (zI - A)^{-1}E + D$. If we take as initial state $\bmx(0) = 0$, then each input sequence $\bmw$ on $\mathbb{Z}_+$ yields a unique output sequence $\bmz$ on $\mathbb{Z}_+$. If $A$ is stable, then this output sequence $\bmz$ is in $\ell_2^p(\mathbb{Z}_+)$ whenever $\bmw$ is in $\ell_2^q(\mathbb{Z}_+)$. The $\mathcal{H}_{\infty}$ performance of \eqref{ch5:eq:isosystem} is now defined as 
$$
J_{\mathcal{H}_\infty}:= \sup_{\|\bmw\|_2 \leq 1} \|\bmz\|_2.
$$
Due to the fact that $A$ is stable, the $\mathcal{H}_{\infty}$ performance is indeed a finite number, and is in fact equal to the 
 $\mathcal{H}_{\infty}$-norm of the transfer matrix $G(z)$, which is given by 
 $$
 \|G\|_{\mathcal{H}_\infty} : = \max_{|z| =1} \|G(z)\|.
 $$
As is well known, the famous bounded real lemma gives necessary and sufficient conditions for the $\mathcal{H}_{\infty}$ performance to be strictly less than a given tolerance:
\begin{proposition}[Discrete-time bounded real lemma]  \label{ch5:prop:BRlemma}
Consider the system \eqref{ch5:eq:isosystem}. Let $\gamma > 0$. Then $A$ is stable and $J_{\mathcal{H}_\infty} < \gamma$ if and only if there exists $P > 0$ such that
\begin{equation}
\bbm P - A^\top P A - C^\top C &  -A^\top P E - C^\top D \\
-E^\top P A - D^\top C & \gamma^2 I - E^\top P E - D^\top D \ebm > 0.
\end{equation}
\end{proposition}
The data-driven ${\cal H}_{\infty}$ control problem in the context of noisy input-state data deals with the true (but unknown) system
\begin{equation} 
    \bmx(t+1) = A_s\bmx(t) + B_s\bmu(t) + \bmw(t), \label{ch4:eq:unknowndisnoisy}
\end{equation}
where $\bmx \in \mathbb{R}^n$ is the state, $\bmu \in \mathbb{R}^m$ is the control input and $\bmw \in \mathbb{R}^n$ is an unknown noise input. The matrices $A_s$ and $B_s$ denote the unknown state and input matrices. As model class $\mathcal{M}$ we take the set of all input-state systems with unknown noise inputs, with given, known, dimensions $n$ and $m$, of the form
\begin{equation} \label{ch5:eq:input state with noise}
\bmx(t+1) = A\bmx(t) + B\bmu(t) +\bmw(t).
\end{equation}
We assume that data $(U_-,X)$ have been collected on the time interval $\{0,1, \ldots,T\}$. Since the noise input $\bmw$ is unknown, the noise samples $w(0),w(1),\dots,w(T-1)$ are not measured, and are therefore not part of the data. However, we adopt the noise model specified in Assumption \ref{ch2:assumption on noise samples}, and assume that the (unknown) matrix $W_- = \bbm w(0) ~ w(1) ~ \cdots ~ w(T-1) \ebm$  satisfies the quadratic matrix inequality \eqref{ch2:asnoise} for a given, known,  partitioned matrix $\Phi \in \bpi_{n,T}$. 

As before, the set $\Sigma_{\calD}$ of all systems in ${\cal M}$ that are consistent with the data $(U_-,X)$ is then equal to the set of all $(A,B)$ that satisfy the QMI 
\begin{equation}
    \label{ch5:ineqAB}
    \begin{bmatrix}
    I \\ A^\top \\ B^\top 
    \end{bmatrix}^\top 
    \begin{bmatrix}
    I & X_+ \\ 0 & -X_- \\ 0 & -U_-
    \end{bmatrix}
    \begin{bmatrix}
    \Phi_{11} & \Phi_{12} \\
    \Phi_{12}^\top & \Phi_{22}
    \end{bmatrix}
    \begin{bmatrix}
    I & X_+ \\ 0 & -X_- \\ 0 & -U_-
    \end{bmatrix}^\top
    \begin{bmatrix}
    I \\ A^\top \\ B^\top 
    \end{bmatrix} \geq 0.
\end{equation}
A standing assumption remains that the unknown system \eqref{ch4:eq:unknowndisnoisy} is consistent with the data, i.e.,  is in $\Sigma_{\calD}$, i.e., $(A_s,B_s)$ satisfies the QMI \eqref{ch5:ineqAB}.

We associate to \eqref{ch5:eq:input state with noise} the output equation
\begin{equation}
\bmz(t) = C\bmx(t) + D\bmu(t), \label{ch5:CD}
\end{equation}
where $\bmz(t) \in \mathbb{R}^p$, and $C$ and $D$ are known matrices. For any $(A,B) \in \Sigma_{\cal D}$, the feedback law $\bmu = K \bmx$ yields the closed-loop system 
\begin{equation}
\label{ch5:closedloop}
\begin{aligned}
    \bmx(t+1) &= (A+BK)\bmx(t) + \bmw(t), \\
    \bmz(t) &= (C+DK)\bmx(t).
\end{aligned}
\end{equation}
Denote the transfer matrix of the closed loop system \eqref{ch5:closedloop} by $G_K(z)$. For any $K$ such that $A+BK$ is stable, the $\mathcal{H}_{\infty}$ performance associated with \eqref{ch5:closedloop} is then given by 
$$
J_{\mathcal{H}_{\infty}}(K) := \|G_K\|_{\mathcal{H}_{\infty}}.
$$
Let $\gamma > 0$. By applying Proposition \ref{ch5:prop:BRlemma} to the closed loop system \eqref{ch5:closedloop}, the matrix $A+BK$ is stable and $J_{\mathcal{H}_{\infty}}(K)  < \gamma$ if and only if there exists a matrix $P > 0$ such that 
\begin{equation} \label{ch5:eq:BRineq}
\bbm P - A_K^\top P A_K - C_K^\top C_K &  -A_K^\top P \\
                     -P A_K                 & \gamma^2 I - P \ebm > 0,
\end{equation}
where $A_K := A+BK$ and $C_K := C+DK$. In order to make this applicable to data-driven $\mathcal{H}_{\infty}$ control design, \eqref{ch5:eq:BRineq} will be restated in a different form.
Clearly, $P>0$ satisfies \eqref{ch5:eq:BRineq} if and only if 
\begin{align}
     P - A_K^\top \left( P +  P (\gamma^2 I - P)^{-1}  P \right) A_K  - C_K^\top C_K &> 0, \label{ch5:eq:Hinf1alt}  \\
    \gamma^2 I - P&> 0.\label{ch5:eq:Hinf2alt}
\end{align}
Since $P^{-\frac{1}{2}} ( \gamma^2 I - P)P^{-\frac{1}{2}} =  \gamma^2 (P^{-1} -\frac{1}{\gamma^2} I)$ and $P +  P (\gamma^2 I - P)^{-1}  P = (P^{-1} -\frac{1}{\gamma^2} I)^{-1}$, the inequalities \eqref{ch5:eq:Hinf1alt} and \eqref{ch5:eq:Hinf2alt} can be reformulated as 
\begin{align} \label{ch5:eq:Hinf1}
    P - A_K^\top (P^{-1} -\frac{1}{\gamma^2} I)^{-1} A_K - C_K^\top C_K &> 0,  \\
    P^{-1} -\frac{1}{\gamma^2} I &> 0. \label{ch5:eq:Hinf2}
    \end{align}
This leads to the following definition of informativity for $\mathcal{H}_{\infty}$ control. 
\begin{definition}[Informativity for $\mathcal{H}_{\infty}$ control] \label{ch5:def:infoHinf}
Let $\gamma >0$. The data $(U_-,X)$ are \emph{informative for $\mathcal{H}_\infty$ control} with performance $\gamma$ if there exist matrices  $P> 0$ and $K$  such that \eqref{ch5:eq:Hinf1} and \eqref{ch5:eq:Hinf2} hold for all $(A,B) \in \Sigma_{\calD}$. 
\end{definition}
Of course, if $K$ satisfies the conditions of Definition~\ref{ch5:def:infoHinf}, then it is a suitable control gain for all $(A,B) \in \Sigma_{\calD}$, in the sense that $A + BK$ is stable and $J_{\mathcal{H}_{\infty}}(K)  < \gamma$ for all $(A,B) \in \Sigma_{\calD}$.

The problem is now to derive necessary and sufficient conditions for informativity, and to find a suitable control gain. By pre- and postmultiplication of \eqref{ch5:eq:Hinf1} by $P^{-1}$ we obtain that \eqref{ch5:eq:Hinf1} and \eqref{ch5:eq:Hinf2} are equivalent to
\begin{equation}  \label{ch5:eq:twoLMI's}
\begin{aligned}  
    Y - A_{Y,L}^\top (Y -\frac{1}{\gamma^2} I)^{-1} A_{Y,L} - C_{Y,L}^\top C_{Y,L} &> 0, \\
    Y -\frac{1}{\gamma^2} I &> 0,
\end{aligned}
\end{equation}
where we define $Y := P^{-1}$, $L := KY$, $A_{Y,L} := AY + BL$ and $C_{Y,L}:=CY+DL$. 
Next, note that \eqref{ch5:eq:twoLMI's} holds if and only if
\begin{equation} \label{ch5:eq:altLMI}
\bbm Y - C_{Y,L}^\top C_{Y,L}   &   A_{Y,L}^\top  \\
            A_{Y,L}                          & Y  -\frac{1}{\gamma^2} I  \ebm > 0,
\end{equation}
which in turn is equivalent to  
\begin{align}
    Y - C_{Y,L}^\top C_{Y,L} &> 0, \label{ch5:CLPCL}\\
    Y - \frac{1}{\gamma^2} - A_{Y,L} (Y-C_{Y,L}^\top C_{Y,L})^{-1} A_{Y,L}^\top &> 0. \label{ch5:ALPAL}
\end{align}
Note that \eqref{ch5:CLPCL} is independent of $A$ and $B$. In addition, we can write \eqref{ch5:ALPAL} as
\begin{equation}\label{ch5:ineqABYL}
\small
\begin{pmat}[{.}]
I \cr\- A^\top \cr B^\top \cr
\end{pmat}^\top 
\hspace{-4pt}
\begin{pmat}[{|}]
Y - \frac{1}{\gamma^2}I & 0 \cr\-
0 & - \begin{bmatrix}
Y \\ L
\end{bmatrix}
(Y - C_{Y,L}^\top C_{Y,L})^{-1}
\begin{bmatrix}
Y \\ L
\end{bmatrix}^\top \cr
\end{pmat}
\hspace{-4pt}
\begin{pmat}[{}]
I \cr\- A^\top \cr B^\top \cr
\end{pmat} \hspace{-1pt} > \hspace{-1pt} 0.
\end{equation}
A crucial observation is now that the inequality \eqref{ch5:ineqABYL} is of a form where $A$ and $B$ appear on the left and their transposes appear on the right, analogous to \eqref{ch3:eq:ineqAB} and \eqref{ch3:eq:ineqABPK}. As before, let
\begin {equation} \label{ch5:eq:N}  
\begin{aligned}
   N = &\begin{pmat}[{|}]
N_{11} & N_{12} \cr\- N_{12}^\top & N_{22} \cr
\end{pmat} 
:= \\
& \begin{pmat}[{.}]
    I & X_+ \cr\- 0 & -X_- \cr 0 & -U_- \cr
    \end{pmat}
    \begin{bmatrix}
    \Phi_{11} & \Phi_{12} \\
    \Phi_{21} & \Phi_{22}
    \end{bmatrix}
    \begin{pmat}[{.}]
    I & X_+ \cr\- 0 & -X_- \cr 0 & -U_- \cr
    \end{pmat}^\top. 
\end{aligned} 
\end{equation}
and let $M$ be defined by
\begin{equation} \label{ch5:eq:M}
\begin{aligned}
M =&  \begin{pmat}[{|}]
M_{11} & M_{12} \cr\- M_{12}^\top & M_{22} \cr
\end{pmat} 
:= \\
& \begin{pmat}[{|}]
Y - \frac{1}{\gamma^2} I& 0 \cr\-
0 & - \begin{bmatrix}
Y \\ L
\end{bmatrix}
(Y - C_{Y,L}^\top C_{Y,L})^{-1}
\begin{bmatrix}
Y \\ L
\end{bmatrix}^\top \cr
\end{pmat}.
\end{aligned}
\end{equation}
Then, for given $\gamma > 0$,  informativity for $\mathcal{H}_{\infty}$ control with performance $\gamma$ holds if and only if 
there exist matrices $Y> 0$ and $L$ that satisfy the inequality $Y - C_{Y,L}^\top C_{Y,L}> 0$ with in addition
\begin{equation} \label{ch5:implicationMNdatafirst}
\begin{aligned}
& \begin{bmatrix} I \\ Z \end{bmatrix}^\top M \begin{bmatrix} I \\ Z \end{bmatrix} > 0 \:\:  \text{ for all }  Z \in \mathbb{R}^{(n+m) \times n} \\
&\text{ such that } \begin{bmatrix} I \\ Z \end{bmatrix}^\top N \begin{bmatrix} I \\ Z \end{bmatrix} \geq 0,
\end{aligned}
\end{equation}
with $Z$ given by 
$$
Z := \begin{bmatrix}
A^\top \\ B^\top 
\end{bmatrix}.
$$
Moreover, in that case a suitable control gain is given by $K = L Y^{-1}$. 

Using the sets \eqref{ch0:e:Zr} and \eqref{ch0:e:Zr+} introduced in 'Quadratic Matrix Inequalities', condition \eqref{ch5:implicationMNdatafirst} is equivalent to
\begin{equation} \label{ch5:eq:QMI inclusion}
\calZ_{n +m}(N) \subseteq \calZ_{n +m}^+(M).
\end{equation}
This observation brings us in a position to apply Corollary~\ref{c:combinedstrictS-lemmaFinslerslemma}. In fact, combining Corollary~\ref{c:combinedstrictS-lemmaFinslerslemma} with some suitable Schur complement arguments leads to the following necessary and sufficient conditions for informativity for $\mathcal{H}_{\infty}$ control with a given performance. In addition, a control gain is computed that achieves $\mathcal{H}_{\infty}$ performance strictly less that $\gamma$ for all systems consistent with the data.
\begin{theorem}[Conditions for informativity]
\label{ch5:theoremHinf}
Suppose that the data $(U_-,X)$ are collected from system \eqref{ch4:eq:unknowndisnoisy} with noise as in Assumption~\ref{ch2:assumption on noise samples}. In addition, let $\gamma > 0$. Then the data $(U_-,X)$ are informative for $\mathcal{H}_{\infty}$ control with performance $\gamma$  if and only if there exist matrices $Y \in \S{n}$, $L \in \mathbb{R}^{m \times n}$ and scalars $\alpha \geq 0$ and $\beta > 0$ satisfying

\begin{equation}
\label{ch5:LMIHinf}
\begin{aligned}
&
\begin{bmatrix}
    Y- \frac{1}{\gamma^2} I - \beta I & 0 & 0 & 0 & 0 \\
    0 & 0 & 0 & Y & 0  \\
    0 & 0 & 0 & L & 0 \\
    0 & Y & L^\top & Y & C_{Y,L}^\top \\
    0 & 0 & 0 & C_{Y,L} & I
    \end{bmatrix}    \hspace{-1mm}  \\
    & - \alpha \hspace{-1mm} \begin{bmatrix}
    I & X_+ \\ 0 & -X_- \\ 0 & -U_- \\ 0 & 0 \\ 0 & 0
    \end{bmatrix} \hspace{-1mm}
    \begin{bmatrix}
    \Phi_{11} & \Phi_{12} \\
    \Phi_{21}& \Phi_{22}
    \end{bmatrix} \hspace{-1mm}
    \begin{bmatrix}
    I & X_+ \\ 0 & -X_- \\ 0 & -U_- \\ 0 & 0 \\ 0 & 0
    \end{bmatrix}^\top   \hspace{-2mm} \geq 0, 
    \\
   & \begin{bmatrix}
    Y & C_{Y,L}^\top \\
    C_{Y,L} & I
    \end{bmatrix}  > 0. 
\end{aligned}
\end{equation}
Moreover, if $Y$ and $L$ satisfy \eqref{ch5:LMIHinf} then $K := L Y^{-1}$ is such that $A+BK$ is stable and $J_{\mathcal{H}_{\infty}}(K) < \gamma$ for all $(A,B) \in \Sigma_{\calD}$.
\end{theorem}
If $Y$ and $L$ satisfy \eqref{ch5:LMIHinf} then $K := L Y^{-1}$ and $P := Y^{-1}$ satisfy \eqref{ch5:eq:BRineq} for all $(A,B)$ in the set $\Sigma_{\calD}$ of systems consistent with the data. Clearly, \eqref{ch5:eq:BRineq} implies that 
\begin{align*} 
& \bmx(t+1)^\top P \bmx(t+1)  - \bmx(t)^\top P \bmx(t)\\ 
& \leq \bbm \bmz(t) \\ \bmw(t) \ebm^\top   \bbm -I  &  0 \\ 0  &  \gamma^2 I\ebm \bbm \bmz(t) \\ \bmw(t) \ebm 
\end{align*}
for all $t \in \mathbb{Z}_+$, where $\bmx$, $\bmw$ and $\bmz$ satisfy the closed loop system equations \eqref{ch5:closedloop}.
This can be interpreted as saying that the system \eqref{ch5:closedloop} is dissipative with respect to the supply rate
\[
s(z,w) : = \bbm z \\ w \ebm^\top \bbm -I  &  0 \\ 0  &  \gamma^2I \ebm \bbm z \\ w \ebm.
\]
with storage function $x^\top P x$. 
In other words: the control law $\bmu = K \bmx$ with $K := L Y^{-1}$ makes all systems in $\Sigma_{\calD}$ dissipative with {\em common} storage function given by $P := Y^{-1}$.

\subsection{Dissipativity analysis}

In this subsection, we study dissipativity of linear finite-dimen\-sional input-state-output systems from a data-driven perspective. This problem has received considerable attention, and we mention the papers \cite{Maupong2017a,Romer2019,Koch2020b} as the approaches that are closest to the one taken here. In \cite{Maupong2017a}, the notion of (finite-horizon) $L$-dissipativity was introduced. This was further studied in \cite{Romer2019}. Both contributions rely on the notion of persistently exciting input data (see \cite{Willems2005} and the sidebar ``Willems' fundamental lemma"). This property of the input sequence implies that the data-generating system is uniquely identifiable from the data.  

In this paper we adopt the more classical notion of dissipativity for linear systems, rather than $L$-dissipativity. Indeed, we consider a setup similar to that of \cite{Koch2020b}, where sufficient data-based conditions were given for dissipativity. Here, we employ the informativity approach to derive necessary and sufficient conditions.

We will first review the definition of dissipativity. Consider a discrete-time linear input-state-output system
\bse\label{ch5:e:lin-sys}
\begin{align}
\bm x(t+1)&=A \bm x(t)+B \bm u(t), \\
\bm y(t)&=C \bm x(t)+D \bm u(t),
\end{align}
\ese
where $A\in\R^{n\times n}$, $B\in\R^{n\times m}$, $C\in\R^{p\times n}$, and $D\in\R^{p\times m}$ are given matrices. 
Let $S \in \S{m+p}$. The system \eqref{ch5:e:lin-sys} is said to be \emph{dissipative\/} with respect to the {\em supply rate}  
\beq\label{ch5:e:supply}
s(u,y)=\bbm u\\y\ebm^\top S \bbm u\\y\ebm
\eeq
if there exists $P \in \S{n}$ with $P \geq 0$ such that the {\em dissipation inequality\/}
\begin{equation}
\label{ch5:e:dispineq}
 \bm x(t+1)^\top P \bm x(t+1) - \bm x(t)^\top P \bm x(t) \leq s\big(\bm u(t),\bm y(t) \big) 
\end{equation} 
holds for all $t \geq 0$ and for all  trajectories $(\bm u,\bm x,\bm y): \mathbb{Z}_+ \rightarrow \mathbb{R}^{m+n+p}$ of \eqref{ch5:e:lin-sys}.
It follows from \eqref{ch5:e:dispineq} that dissipativity with respect to the supply rate \eqref{ch5:e:supply} is equivalent with the feasibility of the linear matrix inequalities $P \geq 0$ and 
\setlength\arraycolsep{2pt}
\beq\label{ch5:eq:KY_- P}
\bbm
I & 0 \\A & B
\ebm^\top
\bbm
P & 0\\0 & -P
\ebm
\bbm
I & 0 \\A & B
\ebm+
\bbm
0 & I\\C & D
\ebm^\top
S
\bbm
0 & I\\C & D
\ebm
\geq 0.
\eeq

In the framework of data-driven system analysis, the system matrices are unknown. The question we want to study then is whether we can verify dissipativity using only the input-state-output data obtained from the unknown system. In the present section we will study this question for the situation that our data are noiseless.

Consider the unknown input-state-output system
\bse\label{ch5:e:tru-sys}
\begin{align}
\bm x(t+1)&=A_{s} \bm x(t)+B_{s} \bm u(t),\\
\bm y(t)&=C_{s} \bm x(t)+D_{s} \bm u(t),  \end{align}
\ese
with $\bm u(t) \in \mathbb{R}^m$, $\bm x(t) \in \mathbb{R}^n$ and $\bm y(t) \in\R^{p}$ the input, state and output. We assume that the dimensions $m,n$ and $p$ are known, but the true system matrices $(A_{s}, B_{s},C_{s},D_{s})$ are unknown. What is known instead are a finite number of input-state-output measurements of \eqref{ch5:e:tru-sys}. 

More concrete, we suppose that we have collected input-state-output data.  Let $U_-,X,X_-,$ and $X_+$ be defined as the previous section and let $Y_-$ be defined in a similar way as $U_-$. Our data are now given by $\calD = (U_-,X,Y_-)$.
These data are assumed to be generated by the true system $(A_s,B_s, C_s,D_s)$, which means that
\begin{equation}  \label{ch5:e:true system compatible}
\begin{bmatrix}
X_+ \\ Y_-
\end{bmatrix} = \begin{bmatrix}
A_s & B_s \\ C_s & D_s
\end{bmatrix} \begin{bmatrix}
X_- \\ U_-
\end{bmatrix}.
\end{equation}
The set of all systems that are consistent with these data is then given by: 
\begin{equation}
\label{ch5:eq:sigma iso}
\Sigma_{(U_-,X,Y_-)}: = \left\{ (A,B,C,D) \mid \begin{bmatrix} X_+ \\ Y_- \end{bmatrix} = \begin{bmatrix} A & B \\ C & D \end{bmatrix} \begin{bmatrix} X_- \\ U_- \end{bmatrix} \right\}.
\end{equation}
It follows from \eqref{ch5:e:true system compatible} that the unknown system $(A_s,B_s,C_s,D_s)$ is contained in $\Sigma_{(U_-,X,Y_-)}$.
Our goal is to infer from the data $(U_- ,X,Y_- )$ whether the unknown system \eqref{ch5:e:tru-sys} is dissipative.  

On the basis of the given data we are unable to distinguish between the systems in $\Sigma_{(U_-,X,Y_-)}$, in the sense that any of these systems could have generated the data. Nonetheless, if  all of these systems are dissipative, then we can also conclude that the true data-generating system \eqref{ch5:e:tru-sys} is dissipative.  
With this in mind, we now define the  property of \emph{informativity for dissipativity} for the case of noiseless data.
\begin{definition}[Informativity of noiseless data]\label{def:dd diss}
The data $(U_-,X,Y_-)$ are \emph{informative for dissipativity\/} with respect to the supply rate \eqref{ch5:e:supply} if there exists a matrix $P \in \S{n}$, $P \geq0$, such that the LMI \eqref{ch5:eq:KY_- P} holds for every system $(A,B,C,D) \in \Sigma_{(U_-,X,Y_-)}$. 
\end{definition}
Note that our definition of informativity for dissipativity requires the systems in $\Sigma_{(U_-,X,Y_-)}$ to be dissipative with a \emph{common} storage function. 

We will restrict ourselves to the case that the number of negative eigenvalues of the matrix $S$ representing the supply rate is equal to the output dimension $p$ and the number of positive eigenvalues of $S$ is equal to the input dimension $m$. In particular then, $S$ is nonsingular. In other words, we will impose the following assumption on the inertia of $S$:
\begin{equation} \label{ch5:e:inertia}
\In(S)=(p,0,m).
\end{equation}
It is a well-known fact that a necessary condition for dissipativity of any system of the form \eqref{ch5:e:lin-sys} is that the input dimension does not exceed the positive signature of $S$. Our assumption requires that the input dimension is equal to this positive signature and in addition that the matrix $S$ is nonsingular. This assumption is satisfied, for example, for the positive-real and bounded-real case.  Indeed, in the positive-real case we have that $m = p$ and
$$
S = \begin{bmatrix}
0 & I_m \\ I_m & 0
\end{bmatrix},
$$
so that $\In(S) = (m,0,m)$. In the bounded-real case we have 
$$
S = \begin{bmatrix}
\gamma^2 I_m & 0 \\ 0 & -I_p
\end{bmatrix}
$$
for some $\gamma > 0$, which implies that $\In(S)=(p,0,m)$. 

Before establishing conditions for informativity for dissipativity, we note that $\Sigma_{(U_-,X,Y_-)}$ contains exactly one element if and only if
\beq\label{ch5:e:full row rank}
\rank \bbm X_-\\U_- \ebm = n + m.
\eeq
in this case, we say the data  $(U_-,X,Y_-)$ are \emph{informative for system identification}.

As the main result of this part we will now show that the noiseless input-state-output data $(U_-,X,Y_-)$ are informative for dissipativity if and only if they are informative for system identification and the unique system consistent with these data is dissipative. In addition, dissipativity of this unknown true system can be expressed in terms of feasibility of an LMI involving the data.
\begin{theorem}[Informativity of noiseless data \cite{vanWaarde2022-dissipativity}] \label{ch5:th:info diss}
Assume that $\In(S)=(p,0,m)$. Then the data $(U_- ,X,Y_- )$ are informative for dissipativity with respect to the supply rate \eqref{ch5:e:supply} if and only if they are informative for system identification
and there exists $P=P^\top \geq0$ such that
\beq\label{ch5:e:exact cond2}
\bbm
X_-\\X_+
\ebm^\top
\bbm
P & 0\\0 & -P
\ebm
\bbm
X_-\\X_+
\ebm+
\bbm
U_- \\Y_- 
\ebm^\top
S
\bbm
U_- \\Y_- 
\ebm
\geq 0.
\eeq
\end{theorem}

Next, we proceed with studying informativity for dissipativity in the case that our input-state-output data are obtained from an unknown system subject to unknown process noise and measurement noise. We assume that the unknown system is given by 
\bse \label{ch5:e:tru-sys with noise}
\begin{align}
\bm x(t+1)&=A_{s} \bm x(t)+B_{s} \bm u(t)  + \bmw(t),\\
\bm y(t)&=C_{s} \bm x(t)+D_{s} \bm u(t) + \bmz(t),  \end{align}
\ese
where $\bm u(t) \in \mathbb{R}^m$, $\bm x(t) \in \mathbb{R}^n$ and $\bm y(t) \in\R^{p}$ are the input, state and output. The dimensions $m,n$ and $p$ are assumed to be known. The terms $\bmw(t) \in \mathbb{R}^n$ and $\bmz(t) \in \mathbb{R}^p$ represent process and measurement noise, respectively, and are assumed to be unknown. Also the system matrices $(A_{s}, B_{s},C_{s},D_{s})$ are assumed to be unknown. Again, we assume that a supply rate is represented by a given matrix $S \in \S{m + p}$, viz. \eqref{ch5:e:supply}. The problem that we will study is whether we can determine whether the unknown system \eqref{ch5:e:tru-sys with noise} is dissipative with respect to the given supply rate. 

Suppose that we obtain  input-state-output data data from the unknown system  \eqref{ch5:e:tru-sys with noise}. These data are collected in the matrices $(U_-,X,Y_-)$. The auxiliary matrices $X_-$ and $X_+$ are as defined before. The noise terms $\bmw$ and $\bmz$ are unknown, so $w(0),w(1),\dots,w(T-1)$ and  $z(0),z(1),\dots,z(T-1)$ are not measured, and are therefore not part of the data. 
We do have the following information on the noise during the data sampling period.
\begin{assumption}[Noise model] \label{ch5:assumption on noise samples}
The noise samples,
collected in the real $(n + p) \times T$ matrix 
$$
V_- := \bbm w(0) & w(1) & \cdots & w(T-1) \\  z(0) & z(1)  & \cdots & z(T-1)    \ebm
$$
satisfy the quadratic matrix inequality
\begin{equation} 
    \label{ch5:asnoise}
    \begin{bmatrix}
    I \\ V_-^\top 
    \end{bmatrix}^\top 
    \Phi
    \begin{bmatrix}
    I \\ V_-^\top 
    \end{bmatrix} \geq 0,
\end{equation}
where $\Phi \in \S{n +p + T}$ is a given partitioned matrix 
\begin{equation} \label{ch5:eq:Phi}
\Phi = \bbm \Phi_{11}  & \Phi_{12} \\ \Phi_{21} & \Phi_{22} \ebm
\end{equation}
with $\Phi_{11} \in \S{n + p}$, $\Phi_{12} \in \mathbb{R}^{(n + p) \times T}$, $\Phi_{21} = \Phi_{12}^\top$ and $\Phi_{22} \in \S{T}$. We assume that $\Phi \in \bpi_{n + p,T}$. Then $\calZ_T(\Phi)$ is nonempty and convex (see Sidebar ``Quadratic matrix inequalities").
We have that $V_-$ satisfies \eqref{ch5:asnoise} if and only if $V_-^\top \in \calZ_T(\Phi)$. 
\end{assumption}

We now turn to defining the  property of \emph{informativity for dissipativity} for noisy input-state-output data, i.e. data that are generated by the unknown system \eqref{ch5:e:tru-sys with noise} with unknown process noise and measurement noise whose samples satisfy the quadratic matrix inequality \eqref{ch5:asnoise}. 
As our model class $\calM$  we take all noisy input-state-output systems
\bse \label{ch5:e:model class}
\begin{align}
\bm x(t+1)&=A \bm x(t)+B \bm u(t)  + \bmw(t),\\
\bm y(t)&=C \bm x(t)+D \bm u(t) + \bmz(t),  \end{align}
\ese
with input dimension $m$, state space dimension $n$ and output dimension $p$. Given the input-state-output data $(U_-,X,Y_-)$ together with the information that the matrices of noise samples satisfy \eqref{ch5:asnoise}, the set of all systems consistent with the data is then given by 
\begin{equation} \label{ch5:def:SigmaD}
\Sigma_{\calD} = \left \lbrace (A,B,C,D) \! \mid \! ( \begin{bmatrix} X_+\\Y_-  \end{bmatrix} \!-\! \begin{bmatrix} A&B\\
C&D\end{bmatrix}\!\begin{bmatrix}X_-\\ U_-  \end{bmatrix})^\top \!\in\!\calZ_T(\Phi)  \right \rbrace.
\end{equation}
We assume that the data have been obtained from the unknown system \eqref{ch5:e:tru-sys with noise}, i.e., $(A_s,B_s,C_s,D_s) \in \Sigma_{\calD}$. Therefore, $\Sigma_{\calD}$ is nonempty. Define
\begin{equation} \label{ch5:eq:bigN}
N\!:= \!\begin{pmat}[{|}]
N_{11} & N_{12} \cr\- N_{12}^\top & N_{22} \cr
\end{pmat} \! = \! \left[\begin{array}{c|c}
I & \begin{array}{c}
X_+\\Y_- 
\end{array}
\\\hline
0 & \begin{array}{c}
-X_-\\-U_- 
\end{array}
\end{array}\right]
\!\!
\bbm
\Phi_{11} & \Phi_{12}\\
\Phi_{21} & \Phi_{22}
\ebm\!\!
\left[\begin{array}{c|c}
I & \begin{array}{c}
X_+\\Y_- 
\end{array}
\\\hline
0 & \begin{array}{c}
-X_-\\-U_- 
\end{array}
\end{array}\right]^\top\!\!.
\end{equation}
Note that $(A,B,C,D)\in \Sigma_{\calD}$ if and only if
\beq \label{ch5:e:char N2 model}
\bbm
I\\\hline\\[-3mm]
\begin{matrix}
A^\top & C^\top\!\\
B^\top & D^\top\!
\end{matrix}
\ebm^\top\!\! 
N
\bbm
I\\\hline\\[-3mm]
\begin{matrix}
A^\top & C^\top\!\\
B^\top & D^\top\!
\end{matrix}
\ebm
\geq 0.
\eeq
This can be restated equivalently as
$$
\bbm
A^\top & C^\top\!\\
B^\top & D^\top\!
\ebm \in \calZ_{n + m}(N).
$$
From Assumption \ref{ch5:assumption on noise samples} we have $\Phi_{22} \leq 0$ and therefore $N_{22} \leq 0$. It follows from the assumption $\ker \Phi_{22} \subseteq \ker \Phi_{12}$ that  $\ker N_{22} \subseteq \ker N_{12}$. 
Since $ \calZ_{n + m}(N)$ is nonempty it follows from Theorem \ref{ch0:t:Z-r nonempty} of Sidebar ``Quadratic matrix inequalities" that $N \schur N_{22} \geq 0$. 
Thus the matrix $N$ given by \eqref{ch5:eq:bigN} is in $\bpi_{n+p, n+m}$.

Next, we give the definition of informativity for dissipativity in the context of noisy input-state-output data. Again, we will require that all systems consistent with the data are dissipative with a {\em common storage function}.
\begin{definition}[Informativity of noisy data]\label{ch5:def:info diss noisy}
The noisy input-state-output data $(U_-,X,Y_-)$ are \emph{informative for dissipativity\/} with respect to the supply rate \eqref{ch5:e:supply} if there exists a matrix $P\geq 0$ such that the LMI \eqref{ch5:eq:KY_- P} holds for all systems $(A,B,C,D)\in\Sigma_{\calD}$. 
\end{definition}
Similar to the noiseless case as studied before, in the remainder of this section we will assume that the matrix $S$ representing the supply rate satisfies the inertia condition $\In(S)=(p,0,m)$.

The following preliminary lemma states that also in the context of noisy data, the rank condition \eqref{ch5:e:full row rank} on the input-state data is necessary for informativity.

\begin{lemma}[Necessity of full row rank condition \cite{vanWaarde2022-dissipativity}] \label{ch5:lem:necc noisy case} 
Assume that $\In(S)=(p,0,m)$. If the data $(U_- ,X,Y_- )$ are informative for dissipativity with respect to the supply rate \eqref{ch5:e:supply} then \eqref{ch5:e:full row rank} holds.
\end{lemma}
In addition, we need the following lemma which states that if the data are informative for dissipativity with all systems in $\Sigma_{\calD}$ having a given common storage function $P \geq 0$, then $P$ is necessarily \emph{positive definite}. This is true under the additional assumption that the Schur complement $N \schur N_{22} $ is positive definite. Combining this with the fact that $N \in \bpi_{n+p, n+m}$ as was already established above, this implies that the set $\Sigma_{\calD}$ has a nonempty interior.
\begin{lemma}[Necessity of positive definite storage \cite{vanWaarde2022-dissipativity}]
\label{ch5:lem:P>0}
Suppose that $\In(S)=(p,0,m)$ and that $N \schur N_{22} >0$. If $P \geq 0$ satisfies the dissipation inequality \eqref{ch5:eq:KY_- P} for all $(A,B,C,D) \in \Sigma_{\calD}$ then $P > 0$. 
\end{lemma}

Our next step is to partition 
\begin{equation} \label{ch5:eq:partitionS}
S=\bbm F & G\\G^\top & H\ebm,
\end{equation}
where  $F\in\R^{m\times m}$, $G\in\R^{m\times p}$, $H\in\R^{p\times p}$. For any $P \geq 0$ define
\begin{equation} \label{ch5:eq:partitionM}
M:=\bbm
P & 0 & 0 & 0\\
0 & F & 0 & G\\
0 & 0 & -P & 0\\
0 & G^\top & 0 & H
\ebm.
\end{equation}
Then the system $(A,B,C,D)$ can be seen to satisfy the dissipation inequality \eqref{ch5:eq:KY_- P} if and only if  
\begin{equation}
\label{ch5:eq:eqM1}
\sysone^\top\!\! M \sysone \geq 0
\end{equation} 
Moreover, with this notation in place, the problem of characterizing informativity for dissipativity is equivalent to finding conditions  for the existence of a matrix $P > 0$ such that the inequality \eqref{ch5:eq:eqM1}
holds for all $(A,B,C,D)$ satisfying the inequality \eqref{ch5:e:char N2 model}.

Our strategy to solve this problem is to invoke the nonstrict matrix S-lemma, Theorem \ref{t:nonstrictS-lemma} of Sidebar ``Quadratic matrix inequalities".
Before we can apply Theorem \ref{t:nonstrictS-lemma} however, note that the inequality \eqref{ch5:eq:eqM1} is in terms of $(A,B,C,D)$ while the inequality \eqref{ch5:e:char N2 model} is in terms of the \emph{transposed} matrices $(A^\top,C^\top,B^\top,D^\top)$. Therefore, we will need an additional dualization result that we formulate in the following lemma. 

\begin{lemma}[Dualization of dissipation inequality \cite{vanWaarde2022-dissipativity}]\label{ch5:lem:diss dual}
Let  $P > 0$ and let $(A,B,C,D)$ be any system with input dimension $m$, state space dimension $n$ and output dimension $p$. Assume that $\In(S) = (p,0,m)$. Define
\begin{equation}
\label{ch5:eq:Shat}
 \hat{S}:=\begin{bmatrix} 0&-I_p\\
 I_m&0 \end{bmatrix} S\inv \begin{bmatrix} 0&-I_m\\I_p&0 \end{bmatrix}.
\end{equation}
Then we have
\begin{equation}
\label{ch5:eq:eqL1}
\begin{bmatrix}
I&0\\
A&B
\end{bmatrix}^\top\!\!\begin{bmatrix}
P&0\\
0&-P
\end{bmatrix}\begin{bmatrix}
I&0\\
A&B
\end{bmatrix}+\begin{bmatrix}
0&I\\
C&D
\end{bmatrix}^\top\!\! S \begin{bmatrix}
0&I\\
C&D
\end{bmatrix}\geq 0
\end{equation}
if and only if 
\begin{align}\label{ch5:eq:L2}
\begin{bmatrix}
I&0 \!\\
A^\top&C^\top\!
\end{bmatrix}^\top\!\!\!\begin{bmatrix}
P^{-1}&0\\
0&-P^{-1}
\end{bmatrix}\!\!\!
\begin{bmatrix}
I&0\!\! \\
A^\top&C^\top\!\!
\end{bmatrix}\!\!+\!\!\begin{bmatrix}
0&I \!\! \\
B^\top&D^\top \!\!
\end{bmatrix}^\top\!\!\!\! \hat{S}\!\! \begin{bmatrix}
0&I \!\! \\
B^\top&D^\top \!\!
\end{bmatrix} \!\!\geq\! 0.
\end{align}
\end{lemma}

Lemma~\ref{ch5:lem:diss dual} can be interpreted as saying that the system defined by the quadruple $(A,B,C,D)$ is dissipative with respect to the supply rate $S$, with storage function $P$ if and only if the dual system $(A^\top,C^\top,B^\top,D^\top)$ is dissipative with respect to the supply rate $\hat{S}$, with storage function $P^{-1}$. A behavioral analogue of this result was obtained in \cite{Willems2002}, Proposition 12. 

Now partition
$$
-S^{-1} = \begin{bmatrix}
\hatF&\hatG\\
\hatG^\top  &\hatH
\end{bmatrix},
$$
where $\hatF=\hatF^\top \in\R^{m\times m}$, $\hatG\in\R^{m\times p}$, and $\hatH=\hatH^\top \in\R^{p\times p}$ and  
define 
\begin{equation} \label{ch5:eq:Mhat}
\hat{M} := \begin{bmatrix}
    P\inv & 0 & 0 & 0 \\
    0 & \hatH & 0 &-\hatG^\top\\
    0 & 0 & -P\inv & 0 \\
    0 & -\hatG & 0 & \hatF
    \end{bmatrix}.
\end{equation}
Then it is easily seen that $(A^\top, C^\top,B^\top, D^\top)$ satisfies the inequality \eqref{ch5:eq:L2} if and only if 
\begin{equation}
\label{ch5:eq:M1}
\bbm
I\\\hline\\[-3mm]
\begin{matrix}
A^\top & C^\top\!\\
B^\top & D^\top\!
\end{matrix}
\ebm^\top\!\! \hat{M} \bbm
I\\\hline\\[-3mm]
\begin{matrix}
A^\top & C^\top\!\\
B^\top & D^\top\!
\end{matrix}
\ebm \geq 0.
\end{equation} 
We may now observe that, under the assumptions that $\In(S)=(p,0,m)$ and $N \schur N_{22} >0$, informativity for dissipativity with respect to the supply rate given by $S$ holds if and only if there exists $P>0$ such that the quadratic inequality \eqref{ch5:eq:M1} holds for all $(A,B,C,D)$ that satisfy the the quadratic inequality \eqref{ch5:e:char N2 model}, equivalently 
\begin{equation} \label{ch5:eq:Z-inclusion}
\calZ_{n +m}(N) \subseteq \calZ_{n+m}(\hat{M}).
\end{equation}
This brings us in position to apply Theorem \ref{t:nonstrictS-lemma} and to obtain the following characterization for informativity for dissipativity for noisy input-state-output data.
\bthe[Informativity of noisy data \cite{vanWaarde2022-dissipativity}] \label{ch5:t:noise 1}
Suppose that the data $(U_-,X,Y_-)$ are collected from system \eqref{ch5:e:model class} with noise as in Assumption~\ref{ch5:assumption on noise samples}. In addition, assume that $\In(S) = (p,0,m)$ and that the data $(U_- ,X,Y_- )$ are such that $N \schur N_{22} >0$. Partition 
\begin{equation} \label{ch5:eq:partition of Sinv}
-S\inv = \begin{bmatrix}
\hatF&\hatG\\
\hatG^\top  &\hatH
\end{bmatrix},
\end{equation}
where $\hatF=\hatF^\top \in\R^{m\times m}$, $\hatG\in\R^{m\times p}$, and $\hatH=\hatH^\top \in\R^{p\times p}$. 
Then the data are informative for dissipativity with respect to the supply rate \eqref{ch5:e:supply} if and only if there exist a real $n \times n$ matrix $Q \in \S{n}$, $Q >0$ and a scalar $\alpha\geq 0$ such that 
\begin{equation}
\begin{bmatrix}
    \! Q & \!0\! & \!0\! & 0 \!\!\! \\
    \! 0 & \!\hatH\! & \!0\! &-\hatG^\top \!\!\! \\
   \! 0 & \!0\! & \!-Q\! & 0 \!\!\! \\
   \! 0 & \!-\hatG\! & \!0\! & \hatF \!\!\!
    \end{bmatrix} \!-\! \alpha\!\!
    \None
     \!\!\!\!\geq\! 0. \label{ch5:eq:LMI2}
\end{equation}
In that case $P : = Q^{-1}$ is a common storage function for all systems consistent with the data.
\ethe

Theorem~\ref{ch5:t:noise 1} provides a tractable method for verifying informativity for dissipativity of noisy data given the noise model introduced in Assumption~\ref{ch5:assumption on noise samples}. The procedure involves solving the linear matrix inequality \eqref{ch5:eq:LMI2} for $Q$ and $\alpha$. Given $Q$, a common storage function $P$ for all systems in $\Sigma_{\calD}$ is also readily computable as $P = Q^{-1}$. 

\section{Auto-regressive systems and noisy input-output data}

Whereas the first two sections of this paper have dealt with input-output systems in state space form together with input-state data, in the current section we will abandon the state space framework and consider input-output systems described by higher order difference equations, also called auto-regressive (AR) systems. Instead of input-state data we will assume to have (noisy) input-output data. In this framework we will discuss data-driven stabilization. 
Several contributions in the literature have also dealt with input-output data \cite{DePersis2020,Koch2021,Steentjes2022,Berberich2023}. A general strategy in these papers is to construct an artificial state-space representation of the system with a state comprised of shifts of the inputs and outputs. This leads to an input-state-output system to which techniques for state data (as discussed before in this paper) are applicable. A drawback of this approach is that the obtained state space systems are non-minimal and of high dimension. Thus a large amount of data can be required for control (see e.g. \cite[Section VIC]{DePersis2020}). In addition, the system matrices of the state-space representation are structured and consist of a combination of known and unknown blocks. Often, this structure is not taken fully into account, which can lead to rather conservative conditions for data-driven control design. Exploiting this prior knowledge of the system matrices is an important problem, which has recently been studied in \cite{Berberich2023}.

Motivated by these limitations of an artificial state space, the main purpose of this section is to discuss a theory on data driven design of stabilizing feedback controllers on the basis of input-output data, without relying on state construction.

\subsection{Stabilization using input-output data}

We consider input-output systems with additive noise represented by auto-regressive AR models of the form
\begin{equation}  \label{ch11:eq:AR}
\begin{aligned}
 & \bmy(t + L) + P_{L-1}\bmy(t + L -1) + \cdots  + P_0\bmy(t)  =\\
 & Q_{L-1}\bmu(t + L-1) + \cdots  + Q_0\bmu(t) + 
 \bmv(t).
\end{aligned}
\end{equation}
Here $L$ is a positive integer, called the {\em order} of the system. The input $\bmu(t)$ and output $\bmy(t)$ are assumed to take their values in $\mathbb{R}^m$ and $\mathbb{R}^p$, respectively. The term $\bmv(t)$ represents unknown noise. The parameters of the model are real $p \times p$ matrices $P_0, P_1, \ldots, P_{L - 1}$ and $p \times m$ matrices $Q_0, Q_1, \ldots ,Q_{L - 1}$.  Using the 
shift operator $(\sigma \bmf)(t) = \bmf(t +1)$, \eqref{ch11:eq:AR} can be written as 
\begin{equation} \label{ch11:eq:AR short}
P(\sigma) \bmy = Q(\sigma) \bmu + \bmv,
\end{equation}
where $P(\xi)$ and $Q(\xi)$ are the real $p \times p$ and $p \times m$ polynomial matrices defined by
\begin{equation} \label{ch11:eq:polmats}
\begin{aligned}
P(\xi) &= I \xi^L + P_{L - 1} \xi^{L -1} + \cdots +P_1 \xi + P_0,  \\
Q(\xi) & = Q_{L -1} \xi^{L -1} + \cdots + Q_1 \xi + Q_0. 
\end{aligned}	
\end{equation}
Since the leading coefficient matrix of $P(\xi)$ is the $p \times p$ identity matrix, $P$ is invertible as a rational matrix and $P^{-1}(\xi) Q(\xi)$ is strictly proper. Thus, indeed, \eqref{ch11:eq:AR short} represents a causal input-output system with control input $\bmu$, noise input $\bmv$ and output $\bmy$. A feedback controller for the input-output system \eqref{ch11:eq:AR short} with 
$P(\xi)$ and $Q(\xi)$ of the form \eqref{ch11:eq:polmats} will be taken to be of the form
\begin{equation} \label{ch11:eq:FBC}
G(\sigma) \bmu = F(\sigma) \bmy,
\end{equation}
with 
\[
\begin{aligned}
G(\xi) &= I \xi^L + G_{L - 1} \xi^{L -1} + \cdots +G_1 \xi + G_0, \\
F(\xi) & =  F_{L -1} \xi^{L -1} + \cdots + F_1 \xi + F_0.
\end{aligned}	
\]
The leading coefficient matrix of $G(\xi)$ is assumed to be the $m \times m$ identity matrix and $G_i \in \mathbb{R}^{m \times m}$, 
$F_i \in \mathbb{R}^{m \times p}$ for $i = 0,1, \ldots, L-1$. The closed loop system obtained by interconnecting a system of the form \eqref{ch11:eq:AR short} and the controller is represented by
\begin{equation} \label{ch11:eq:closedloop}
\bbm     G(\sigma) &  -F(\sigma)  \\
             -Q(\sigma) &   P(\sigma)  \ebm   \bbm \bmu  \\ \bmy \ebm = \bbm 0 \\ I_p  \ebm \bmv.
\end{equation}
Note that the leading coefficient matrix is the $q \times q$ identity matrix. We call the controller \eqref{ch11:eq:FBC} a stabilizing controller for \eqref{ch11:eq:AR short} if the corresponding autonomous system 
\begin{equation} \label{auto}
\bbm     G(\sigma) &  -F(\sigma)  \\
             -Q(\sigma) &   P(\sigma)  \ebm   \bbm \bmu  \\ \bmy \ebm = 0
\end{equation}
is stable, in the sense that all solutions  $\bmu$ and $\bmy$ of \eqref{auto} tend to zero as time tends to infinity.
The problem that we consider is to find a feedback controller of the form \eqref{ch11:eq:FBC} that stabilizes the unknown true system 
\begin{equation} \label{unknown true system}
P_s(\sigma) \bmy = Q_s(\sigma) \bmu + \bmv.
\end{equation}
For this, we assume that the order $L$ is known, and that only data obtained from the true system can be used. These data are the input-output data given by $u(0),u(1), \ldots, u(T)$, $y(0), y(1), \ldots, y(T)$ on the interval $[0,T]$ with $T \geq L$.  These are samples of $\bmu$ and $\bmy$ satisfying the system equation \eqref{unknown true system} for some noise signal $\bmv$. 
The noise $\bmv$ is unknown, but its samples are assumed to satisfy an assumption analogously to Assumption \ref{ch2:assumption on noise samples}:
\begin{assumption}[Assumption on the noise] \label{A}
The noise samples $v(0), v(1), \ldots ,v(T - L)$, collected in the real $p \times (T - L +1)$ matrix 
\[
V : = \bbm v(0) & v(1) & \cdots & v(T -L) \ebm
\]
satisfy the quadratic matrix inequality
\begin{equation} \label{ch11:eq:noiseQMI}
\bbm I \\ V^\top \ebm^\top \Pi \bbm I \\ V^\top \ebm \geq 0,
\end{equation}
where $\Pi \in \S{p + T -L +1}$ is a known partitioned matrix
\[
\Pi = \bbm \Pi_{11} & \Pi_{12} \\ \Pi_{21}   & \Pi_{22} \ebm,
\]
with $\Pi_{11} \in \S{p}$, $\Pi_{12} \in \mathbb{R}^{p \times (T - L +1)}$, $\Pi_{21}  = \Pi_{12}^\top$ and $\Pi_{22} \in \S{T - L + 1}$. We assume that $\Pi \in \bpi_{p,T - L +1}$ By Proposition \ref{ch0:t:Z-r nonempty} the set of matrices $V$ that satisfy \eqref{ch11:eq:noiseQMI} is nonempty. 
\end{assumption}

As noted before in this paper, in general the given data $u(0),u(1), \ldots, u(T)$, $y(0), y(1), \ldots, y(T)$ do not determine the true system uniquely. In fact, the data determine a whole set of systems that are consistent with the data. As a consequence, finding a stabilizing controller for the true system based only on the data requires finding a controller that stabilizes all systems that are consistent with the data. If, for given data, such controller exists, then we call the input-output data informative for stabilization. This will now be made precise. In order to do this, first the set of all systems that are consistent with the data will be specified.

After denoting $q : = p + m$, $R(\xi) = \bbm  -Q(\xi)  & P(\xi) \ebm$ and $\bmz = \col(\bmu,\bmy)$,
\eqref{ch11:eq:AR short} can be rewritten as 
\begin{equation} \label{ch11:eq:ARkernel}
R(\sigma) \bmz = \bmv.
\end{equation}
Collect the (unknown) coefficient matrices of the polynomial matrix $R(\xi)$ in the $p \times qL $ matrix
\begin{equation} \label{ch11:eq:CoeffR}
R:= \bbm -Q_0 & P_0 & -Q_1 & P_1 & \cdots & -Q_{L - 1} & P_{L - 1} \ebm 
\end{equation}
Note that, with a slight abuse of notation, we denote both the polynomial matrix and its coefficient matrix by $R$.
We call \eqref{ch11:eq:CoeffR} the coefficient matrix of the system \eqref{ch11:eq:ARkernel}. Arrange the data $u(0),u(1), \ldots, u(T), y(0), y(1), \ldots, y(T)$ into the vectors 
\[
z(t) = \bbm u(t) \\ y(t) \ebm,~~ (t = 0,1, \ldots ,T )
\]
and define the associated depth $L + 1$ Hankel matrix by 
\begin{equation} \label{ch11:eq:originalH}
H(z) := \bbm z(0) & z(1) & \cdots & z(T - L) \\
	         z(1) & z(2) & \cdots & z(T - L +1) \\
                     \vdots  & \vdots &    &  \vdots  \\
             z(L-1)  & z(L) &    \cdots  & z(T-1) \\
                     y(L)  &   y(L + 1) & \cdots & y(T)
                     \ebm.
\end{equation}
Furthermore, partition 
\begin{equation} \label{ch11:eq:partH}
H(z)= \bbm H_1(z) \\ H_2(z) \ebm,
\end{equation}
where $H_1(z)$ contains the first $qL$ rows and $H_2(z)$ the last $p$ rows. It is then easily verified that any input-output system \eqref{ch11:eq:ARkernel} for which the coefficient matrix $R$ defined in \eqref{ch11:eq:CoeffR} satisfies 
\begin{equation} \label{ch11:eq:iocomp}
\bbm R & I    \ebm \bbm H_1(z) \\ H_2(z) \ebm = V
\end{equation}
for some $V \in \calZ_{T - L +1}(\Pi)$,
could have generated the given input-output data. In other words, $z(0),z(1),$ $\ldots, z(T)$ are also samples on the interval $[0,T]$ of a $\bmz$ that satisfies 
$
R(\sigma) \bmz = \bmv
$
for some $\bmv$ satisfying Assumption \ref{A}.
Therefore, $R$ satisfies \eqref{ch11:eq:iocomp} for some $V \in \calZ_{T - L +1}(\Pi)$ if and only if the AR system with coefficient matrix $R$ is consistent with the data. Recall that, in particular, the true system is consistent with the data.
Now define
\begin{equation} \label{ch11:eq:defN}
N := \bbm I  & H_2(z) \\
                 0  &  H_1(z) 
                      \ebm \Pi
        \bbm I  & H_2(z) \\
                 0  &  H_1(z) \\
                         \ebm^\top.
\end{equation}
Then by combining \eqref{ch11:eq:noiseQMI} and \eqref{ch11:eq:iocomp} we see that the system with coefficient matrix $R$ is consistent with the data if and only if $R^\top$ satisfies the QMI
\begin{equation} \label{ch11:eq:QMIcomp}
\bbm I \\ R^\top  \ebm^\top \! \! N \bbm I \\ R^\top  \ebm \geq 0.
\end{equation}
Thus we have succeeded in finding an explicit expression for the set of systems that are consistent with the data. Indeed, this set is equal to 
\[
\Sigma_{\cal D} = \{ R \in \mathbb{R}^{p \times qL} \mid R^\top \in \calZ_{qL}(N)\}.
\]
Since the true system is consistent with the data, this set is nonempty.

Our aim is to find a single controller of the form \eqref{ch11:eq:FBC} that stabilizes all input-output systems \eqref{ch11:eq:ARkernel} that are consistent with the data, so all systems in $\Sigma_{\calD}$. In order to investigate the existence of such controller, we will now first study stability of autonomous systems in AR form. 

Given a nonsingular $p \times p$ polynomial matrix $P(\xi)$, the corresponding autonomous AR system $P(\sigma) \bmy = 0$ is called stable if $\bmy(t) \to 0$ as $t \to \infty$ for all solutions $\bmy : \mathbb{Z}_+ \rightarrow \mathbb{R}^p$. This space of all solutions on $\mathbb{Z}_+$ is called the behavior of the system and is denoted by ${\calB}(P)$.
Stability of autonomous AR systems can be characterized in terms of quadratic difference forms on behaviors. For details on QDFs, see: `Quadratic Difference Forms`. In a continuous-time context, the connection between stability and QDFs was studied in \cite{Willems1998}, while the discrete-time version was considered in \cite{Kojima2005}. The following proposition holds:

\begin{sidebar}{Quadratic difference forms}
\label{sidebar:QDFs}

\setcounter{sequation}{0}
\renewcommand{\thesequation}{S\arabic{sequation}}
\setcounter{stable}{0}
\renewcommand{\thestable}{S\arabic{stable}}
\setcounter{sfigure}{0}
\renewcommand{\thesfigure}{S\arabic{sfigure}}

\sdbarinitial{A} crucial instrument in studying stability of systems is the notion of Lyapunov function. Studying stability of autonomous systems in AR form requires the notion of Lyapunov functions given by {\em quadratic difference forms} (QDFs). In this sidebar we review the basic material on QDFs and establish some useful preliminary results. For more details, we refer to \cite{Willems1998, Willems2002, Kojima2005, Kojima2006}. 

Let  $N$ and $q$ be positive integers and for $i,j = 0,1,\ldots ,N$ let $\Phi_{i,j} \in \mathbb{R}^{q \times q}$ be such that $\Phi_{i,i} \in \S{q}$ and $\Phi_{i,j} = \Phi_{j,i}^\top$ for all $i \neq j$. Arrange these matrices into the partitioned matrix $\Phi \in \S{(N + 1)q}$ given by
\begin{equation*}
\Phi := 
\begin{bmatrix}
\Phi_{0,0} & \Phi_{0,1} & \cdots & \Phi_{0,N} \\

\Phi_{1,0} &\Phi_{1,1}& \cdots &\Phi_{1,N} \\
\vdots & \vdots & \ddots & \vdots \\
\Phi_{N,0} & \Phi_{N,1}& \cdots &\Phi_{N,N}
\end{bmatrix}.
\end{equation*}
Then the quadratic difference form associated with $\Phi$ is the operator $Q_{\Phi}$ that maps $\mathbb{R}^q$-valued functions $\bmz$ on $\mathbb{Z}_+$ to  $\mathbb{R}$-valued functions $Q_\Phi(\bmz)$ on $\mathbb{Z}_+$ defined by
\begin{sequation}  \label{ch0:eq:QDF}
Q_{\Phi}(\bmz)(t) := \sum_{k,\ell=0}^{N} \bmz(t + k)^\top \Phi_{k,\ell} ~\bmz(t + \ell).
\end{sequation}
In terms of the matrix $\Phi$ this can be written as 
\[
Q_{\Phi}(\bmz)(t) = \bbm \bmz(t ) \\ \bmz(t + 1) \\ \vdots \\  \bmz(t+N) \ebm^\top \!\!\Phi   \bbm \bmz(t ) \\ \bmz(t + 1) \\ \vdots \\  \bmz(t+N) \ebm.
\]
Thus, vector valued functions are mapped to quadratic expressions in terms of these function and their time shifts up to a certain dergree.

Obviously, some of the matrices $\Phi_{i,j}$, or even an entire block row or column of $\Phi$ could be zero. We define the {\em degree} of the QDF \eqref{ch0:eq:QDF} as the smallest integer $d$ such that $\Phi_{ij} = 0$ for all $i > d$ or $j >d$. This degree is denoted by $\deg(Q_{\Phi})$. The matrix $\Phi$ is called a coefficient matrix of the QDF. Note that a given QDF does not determine the coefficient matrix uniquely. However, if the degree of the QDF is $d$, it allows a coefficient matrix $\Phi \in \S{(d + 1)q}$.

The QDF $Q_\Phi$ is called nonnegative if $Q_{\Phi}(\bmz) \geq 0$ for all $\bmz: \mathbb{Z}_+ \rightarrow \mathbb{R}^q$. We denote this as $Q_{\Phi} \geq 0$. Clearly, this holds if and only if $\Phi \geq 0$. The QDF is called positive if it is nonnegative and, in addition, $Q_{\Phi}(\bmz)= 0$ if and only if $\bmz = 0$.  This is denoted as $Q_{\Phi} > 0$. Likewise we define nonpositivity and negativity.

For a given QDF $Q_\Phi$, its {\em rate of change} along a given $\bmz: \mathbb{Z}_+ \rightarrow \mathbb{R}^q$ is given by
$Q_\Phi(\bmz)(t+1) - Q_\Phi(\bmz)(t)$. It turns out that the rate of change defines a QDF itself. Indeed, by defining the matrix 
$\nabla \Phi \in \S{(N + 2) q}$ by
\begin{sequation} \label{ch0:eq:nabla phi}
\nabla \Phi := \bbm 0_{q \times q} & 0 \\
                               0  &  \Phi \ebm -
                     \bbm \Phi & 0 \\
                               0  &  0_{q \times q} \ebm,                              
\end{sequation}
it is easily verified that 
\[
Q_{\nabla \Phi}(\bmz)(t) = Q_\Phi(\bmz)(t+1) - Q_\Phi(\bmz)(t)
\]
for all $\bmz: \mathbb{Z}_+ \rightarrow \mathbb{R}^q$ and $t \in \mathbb{Z}_+$.

Quadratic difference forms are particularly relevant in combination with behaviors defined by AR systems. Let $R(\xi)$ be a real $p \times q$ polynomial matrix and consider the AR system represented by
$R(\sigma) \bmz = 0$.
Let ${\cal B}(R)$ be the behavior of this system, i.e., the space of all solutions as given by
\[
{\cal B}(R): = \{ \bmz : \mathbb{Z}_+ \rightarrow \mathbb{R}^q \mid  R(\sigma) \bmz = 0 \}.
\]
The QDF $Q_\Phi$ is called nonnegative on ${\cal B}(R)$ if $Q_{\Phi}(\bmz) \geq 0$ for all $\bmz \in {\cal B}(R)$.  It is called positive on ${\cal B}(R)$ if, in addition, $Q_{\Phi}(\bmz)= 0$ if and only if $\bmz = 0$. We denote this as $Q_\Phi \geq 0$ on ${\cal B}(R)$ and $Q_\Phi > 0$ on ${\cal B}(R)$, respectively. Likewise we define nonpositivity and negativity on ${\cal B}(R)$.

\end{sidebar}

\begin{proposition}[A QDF as Lyapunov function \cite{vanWaarde2022-QDF}] \label{ch0:prop:ARstability}
Let $P(\xi)$ be a nonsingular polynomial matrix. The corresponding autonomous system $P(\sigma) \bmy = 0$ is stable if and only if there exists a QDF $Q_\Psi$ such that $Q_\Psi \geq 0$ on ${\calB}(P)$ and $Q_{\nabla \Psi} < 0$ on ${\calB}(P)$.
\end{proposition}
For obvious reasons, we refer to $Q_\Psi$  as a \emph{Lyapunov function}. In principle, the above theorem does not specify the degree of $Q_{\Psi}$. However, it turns out that if $P(\xi)$ is of the form 
\begin{equation} \label{ch11:eq:Pxi}
P(\xi) = I \xi^L + P_{L -1} \xi^{L - 1} + \ldots P_1 \xi + P_0
\end{equation}
(with leading coefficient matrix the identity matrix) and the corresponding autonomous system $P(\sigma) \bmy = 0$ of order $L$ is stable, there exists a Lyapunov function of degree at most $L- 1$. Indeed, we have
\begin{lemma}[A degree bound on the Lyapunov QDF \cite{vanWaarde2022-QDF}] \label{ch11:lem:Lyapunovdegree}
Let $P(\xi)$ be a polynomial matrix of the form \eqref{ch11:eq:Pxi}. The corresponding autonomous system $P(\sigma) \bmy = 0$ order $L$ is stable if and only if there exists a  QDF $Q_\Psi(\bmy)$ of degree at most $L - 1$ such that $Q_\Psi \geq 0$ and $Q_{\nabla \Psi} < 0$ on ${\calB}(P)$.
\end{lemma}
The fact that the degree of the QDF defining the Lyapunov function can be bounded from above by the order of the system is crucial for enabling us to express stability of the system $P(\sigma)\bmy = 0$ in terms of a quadratic matrix inequality. This QMI involves a symmetric matrix $\Psi$ of dimensions $pL \times pL$ leading to a Lyapunov function $Q_\Psi$, and the matrix $P = \bbm P_0 & P_1 & \cdots & P_{L - 1} \ebm$. Again, for ease of notation we denote both the polynomial matrix and its coefficient matrix by $P$. Then we have:
\begin{theorem}[A QMI condition for stability \cite{vanWaarde2022-QDF}] \label{ch11:th:LyapunovQMI}
Let $P(\xi) = I \xi^L + P_{L - 1} \xi^{L -1} + \ldots + P_1 \xi + P_0$ and let $P(\sigma)\bmy = 0$ be the corresponding autonomous system. This system is stable if and only if there exists $\Psi \in \S{pL}$ such that $\Psi \geq 0$ and
	\begin{equation}
	\label{ch11:eq:LyapunovQMI}
	\begin{bmatrix}
	I \\ -P
	\end{bmatrix}^\top \left( \begin{bmatrix}
	0_p& 0 \\ 0 & \Psi
	\end{bmatrix} - \begin{bmatrix}
	\Psi & 0 \\ 0 & 0_p
	\end{bmatrix} \right) \begin{bmatrix}
	I  \\ -P
	\end{bmatrix} < 0.
	\end{equation}
Any such $\Psi$ defines a Lyapunov function $Q_{\Psi}$.
\end{theorem}

Next, we turn to the data-driven stabilization problem. For a given controller of the form \eqref{ch11:eq:FBC}, denote
\[
 C(\xi) : =  \bbm G(\xi)  &   -F(\xi)  \ebm, 
\]
and recall that $\bmz = \col(\bmu, \bmy)$.
Then \eqref{ch11:eq:closedloop} can equivalently be written as  
\begin{equation} \label{ch11:eq:RoverC}
\bbm C(\sigma) \\
          R(\sigma)   \ebm \bmz = \bbm 0 \\ I_p  \ebm \bmv.
 \end{equation}
Collect the coefficient matrices of $F(\xi)$ and $G(\xi)$ in the matrix $C$ defined by
\begin{equation} \label{ch11:eq:controller coeff}
C := \bbm G_0 & -F_0  & G_1  & -F_1 & \cdots & G_{L-1}  & -F_{L -1}\ebm
\end{equation}
and recall the definition \eqref{ch11:eq:CoeffR} of the coefficient matrix $R$ associated likewise with $R(\xi)$. Note that the leading coefficient matrix of the polynomial matrix
$\bbm C(\xi)^\top & R(\xi)^\top \ebm^\top$ is the $q \times q$ identity matrix. Furthermore, its coefficient matrix is
$
\bbm C^\top & R^\top \ebm^\top.
$
Recall that the controller \eqref{ch11:eq:FBC} is a stabilizing controller for the input-ouput system \eqref{ch11:eq:ARkernel}
if and only if the autonomous system \eqref{auto}
is stable. 
As an immediate consequence of Theorem \ref{ch11:th:LyapunovQMI} we then have
\begin{lemma}[A QMI condition for stabilization \cite{vanWaarde2022-QDF}] \label{ch11:lem:LyapunovQMIcontrolled}
The controller $C(\sigma) \bmz = 0$ is a stabilizing controller for the system \eqref{ch11:eq:ARkernel} if and only if there exists $\Psi \in \S{qL}$ such that $\Psi \geq 0$ and
	\begin{equation}
	\label{ch11:eq:LyapunovQMIcontrolled}
	\begin{bmatrix}
	I_{qL} \\ -C \\ -R
	\end{bmatrix}^\top \left( \begin{bmatrix}
	0_q& 0 \\ 0 & \Psi
	\end{bmatrix} - \begin{bmatrix}
	\Psi & 0 \\ 0 & 0_q
	\end{bmatrix} \right) \begin{bmatrix}
	I_{qL}  \\ -C  \\ -R
	\end{bmatrix} < 0.
	\end{equation}
	Moreover, if $\Psi \geq 0$ satisfies \eqref{ch11:eq:LyapunovQMIcontrolled}, then $\Psi > 0$.
\end{lemma}

Now recall that our aim is to find a single stabilizing controller for all systems in  $\Sigma_{\cal D} = \{ R \in \mathbb{R}^{p \times qL} \mid R^\top \in \calZ_{qL}(N)\}$, i.e. for all systems whose coefficient matrix satisfies the QMI \eqref{ch11:eq:QMIcomp}.
This leads to the following definition of {\em informativity for quadratic stabilization}.
\begin{definition}[Informativity for quadratic stabilization] \label{ch11:def:info auto AR control}
The input-output data $u(0),  \ldots, u(T), y(0),\ldots, y(T)$ are called {\em informative for quadratic stabilization} if there exist $C \in \mathbb{R}^{m \times qL}$ and $\Psi \in \S{qL}$ with $\Psi \geq 0$ such that the QMI \eqref{ch11:eq:LyapunovQMIcontrolled} holds for all $R$ that satisfy the QMI \eqref{ch11:eq:QMIcomp}, with $N$  defined by \eqref{ch11:eq:defN}.
\end{definition}

Informativity for quadratic stabilization thus means that there exist a controller $C(\sigma) \bmz = 0$ (equivalently, $G(\sigma) \bmu = F(\sigma) \bmy$) and a matrix $\Psi \in \S{qL}$ such that the QDF $Q_\Psi$  is a common Lyapunov function for all closed loop systems obtained by interconnecting the controller with an arbitrary system that is consistent with the data.

We now aim at finding necessary and sufficient conditions on the given data to be informative for quadratic stabilization. 
Define the $q(L -1) \times qL$ matrix $J$ by
\begin{equation} \label{ch11:eq:Zcontrolled}
 J:= \bbm 0_{q(L -1) \times q} & I_{q(L -1)} \ebm.
 \end{equation}
It can be proven that  \eqref{ch11:eq:LyapunovQMIcontrolled} holds if and only if $\Psi > 0$ and 
\begin{equation} \label{ch11:eq:bigQMIcontrolled}
\bbm I_{qL}  \\ R^\top \bbm 0 & 0  & -I_p \ebm \ebm^\top 
M
 \bbm I_{qL}  \\ R^\top \bbm 0 &0  &  -I_p  \ebm \ebm > 0.
 \end{equation}
where the $2qL \times 2qL$ matrix $M$ is defined by
\begin{equation} \label{ch11:eq:Mmatrixcontrolled}
M : = \bbm \Psi^{-1} - \bbm J  \\ -C \\  0 \ebm \Psi^{-1} \bbm J \\ -C   \\  0 \ebm^\top &  \bbm J   \\ -C \\ 0 \ebm \!\Psi^{-1} \\
      \!\!    \Psi^{-1} \!\!\bbm J   \\ -C \\0 \ebm^\top         & -\Psi^{-1} \ebm.
\end{equation}
This means that informativity for quadratic stabilization is equivalent to the existence of an $m \times qL$ matrix $C$ and a matrix $\Psi \in \S{qL}$, $\Psi > 0$ such that the QMI  \eqref{ch11:eq:bigQMIcontrolled} holds for all matrices $R$ that satisfy the QMI \eqref{ch11:eq:QMIcomp}. The matrix $C$ is then the coefficient matrix of a suitable controller.
In terms of solutions sets of QMIs this can be restated as
\[
R^\top \in {\calZ}_{qL}(N) ~\Longrightarrow ~ R^\top \bbm 0 & 0 & -I_p \ebm \in {\calZ}^+_{qL}(M),
\]
or equivalently,
\begin{equation} \label{ch1:eq:notyetcontrolled}
 {\calZ}_{qL}(N) \bbm 0   & 0   & -I_p \ebm \subseteq {\calZ}^+_{qL}(M).
\end{equation}
In order to be able to apply the strict matrix S-lemma in Theorem \ref{t:strictS-lemmaN22}, we will express the set on the left in \eqref{ch1:eq:notyetcontrolled} as the solution set of a QMI. Define the $2qL \times 2qL$ matrix $\bar{N}$ by 
\begin{equation} \label{ch11:eq:Nbarcontrolled}
\bar{N} : = \bbm \bbm 0 &0  &-I_p \ebm  & 0 \\
                           0        &    I_{qL}    \ebm^\top
                           N
                   \bbm \bbm 0 &  0&  -I_p  \ebm  & 0 \\
                           0      &      I_{qL}    \ebm.    
\end{equation}
Then indeed the following can be proven:
\begin{lemma}[An instrumental lemma \cite{vanWaarde2022-QDF}]\label{ch11:lem:projectedQMIcontrolled}
Assume that the Hankel matrix $H_1(z)$  has full row rank. Then ${\calZ}_{qL}(N) \bbm 0   & 0   &-I_p  \ebm = {\calZ}_{qL}(\bar{N})$.
\end{lemma}

From the above we see that, under the assumption that $H_1(z)$ has full row rank, informativity for quadratic stabilization requires the existence of $C$ and $\Psi > 0$ such that the inclusion $ {\calZ}_{qL}(\bar{N})  \subseteq {\calZ}^+_{qL}(M)$. holds. This inclusion is dealt with in Theorem \ref{t:strictS-lemmaN22}. 
\begin{lemma}[A condition for informativity \cite{vanWaarde2022-QDF}]  \label{ch11:lem:stillnonlinear controlled}
Let $\Psi > 0$, $C \in \mathbb{R}^{m \times qL}$ and let $M$ be given by \eqref{ch11:eq:Mmatrixcontrolled}. Assume that $H_1(z)$ has full row rank. Then $ {\calZ}_{qL}(\bar{N})  \subseteq {\calZ}^+_{qL}(M)$ if and only if there exists a scalar $\alpha \geq 0$ such that 
\begin{equation} \label{ch11:eq:stillnonlinear controlled} 
M - \alpha \bar{N} > 0.
\end{equation}
\end{lemma}
Note that the unknowns $C$ and $\Psi$ appear in the matrix $M$ in a nonlinear way, and even in the form of an inverse. However, by putting $\Phi : = \Psi^{-1}$ we can get rid of the inverse, and rewrite the condition $M - \alpha \bar{N} > 0$ as
\begin{equation}  \label{ch0:th:NS}
\bbm ~~~\Phi - \bbm J  \\ -C \\  0 \ebm \Phi \bbm J \\ -C   \\  0 \ebm^\top &  \bbm J   \\ -C \\ 0 \ebm \!\Phi \\
      \!\!    \Phi  \!\!\bbm J   \\ -C \\0 \ebm^\top         & -\Phi \ebm   - \alpha \bar{N} > 0.
\end{equation}
Thus, informativity for quadratic stabilization holds if and only if there exists $\Phi > 0$, a matrix $C$ and a scalar $\alpha \geq 0$ such that \eqref{ch0:th:NS} holds. Note that $\alpha$ must be positive due to the negative definite lower right block in \eqref{ch0:th:NS}. By scaling $\Phi$ we can therefore take $\alpha = 1$. Finally, by introducing the new variable $D := - C \Phi$ and taking a suitable Schur complement, \eqref{ch0:th:NS} can be reformulated as the following LMI in the unknowns $\Phi$ and $D$:
\begin{equation} \label{ch11:eq:mainLMI}
\bbm
\Phi  &   \bbm   J \Phi \\ D \\ 0 \ebm    &  \bbm J \Phi \\ D \\ 0 \ebm   \\
   \bbm   J \Phi \\ D \\ 0 \ebm^\top &    -\Phi     &   0     \vspace{1mm}               \\  
 \bbm   J \Phi \\ D\\ 0 \ebm^\top    &    0       &   \Phi 
\ebm   
  - \bbm \bar{N}    &     0  \\
            0                                   &    0_{qL}   \ebm >0.
\end{equation}    
This then immediately leads to the following characterization of informativity for quadratic stabilization and a method to compute a suitable feedback controller together with a common Lyapunov function.               
\begin{theorem}[An LMI condition for informativity \cite{vanWaarde2022-QDF}] \label{ch11:th:mainNS} 
Suppose that the data $u(0), u(1), \ldots, u(T)$, $y(0), y(1), \ldots, y(T)$ are collected from system \eqref{unknown true system} with noise as in Assumption~\ref{A}. In addition, assume that $H_1(z)$ has full row rank. Let the matrix $\bar{N}$ be given by \eqref{ch11:eq:Nbarcontrolled}, with $N$ defined by \eqref{ch11:eq:defN}.
Then the input-output data are informative for quadratic stabilization if and only if there exist  matrices $D \in \mathbb{R}^{m \times qL }$ and $\Phi \in \S{qL}$ such that $\Phi > 0$ and the LMI \eqref{ch11:eq:mainLMI} holds.

In that case, the feedback controller with coefficient matrix $C: = -D \Phi^{-1}$  stabilizes  all systems that are consistent with the input-output data. Moreover, the QDF $Q_{\Psi}$ with $\Psi : = \Phi^{-1}$ is a common Lyapunov function for all closed loop systems. 
\end{theorem}

Thus, in order to compute a controller that stabilizes all systems consistent with the data and which gives a common Lyapunov function, first compute the matrix $\bar{N}$ using the Hankel matrix associated with the data. Next, check feasibility of the LMI \eqref{ch11:eq:mainLMI} and, if it is feasible, compute solutions  $D$ and $\Phi$. 
An AR representation of the controller with coefficient matrix $C = -D \Phi^{-1}$ is then obtained as follows: partition
\[
C := \bbm    G_0 & -F_0  & G_1 & -F_1 &  \cdots &  G_{L-1} & -F_{L -1} \ebm,
\]
with $F_i \in \mathbb{R}^{m \times p}$ and $G_i \in \mathbb{R}^{m \times m}$. Next define $F(\xi):= F_{L -1} \xi^{L - 1} + \cdots + F_0$ and $G(\xi) := I \xi^L + G_{L-1}\xi^{L -1} + \cdots + G_0$. The corresponding controller is then given in AR representation by $G(\sigma) \bmu = F(\sigma) \bmy$.

\subsection{Simulation example }

In this example, we consider a model of a magnetic suspension system, where an electromagnet is used to levitate a magnetic mass. We assume that we can measure the vertical position of the mass and control the current of the electromagnet with the aim of stabilizing the mass at a pre-determined position. Of course, in the context of this paper, we will develop such a controller on the basis of collected measurements. 

For a detailed derivation of the model, see \cite[Example 1.18]{Khalil2015}. Following \cite[Example 12.8]{Khalil2015}, we let $x_1$ denote the vertical position and $x_2$ the vertical velocity of the ball. Moreover, $x_3$ denotes the current and $\hat{u}$ the voltage of the circuit. The model is then given by:
\begin{align*}
	\dot{x}_1 & = x_2, \\ 
	\dot{x}_2 & = g- \frac{k}{m}x_2 - \frac{L_0a x_3^2}{2m(a+x_1)^2}, \\ 
	\dot{x}_3 &= \frac{1}{L(x_1)} \left( - Rx_3 +\frac{L_0ax_2x_3}{(a+x_1)^2}+\hat{u}\right),
\end{align*}
where $L(x_1) = L_1+\frac{L_0a}{a+x_1}$. Defining the function $f$ accordingly, we write this system as $\dot{x}= f(x,\hat{u})$. 

As noted, we are interested in stabilizing the ball at $x_1 =r>0$, on the basis of measurements of $x_1$. In order to apply the results of this paper, we will first linearize the model around the corresponding equilibrium point. After this, we will discretize and rewrite it to an AR model of the form considered in this paper. 

For the physical quantities we will use the following values: 
\begin{align*}
	m	& = 0.1 \hspace{1em}\textrm{kg}, &	L_0 & = 0.01 \hspace{1em}\textrm{H}, \\  
	k	& = 0.001 \hspace{1em}\textrm{N}/\textrm{m}/\textrm{sec}, &	L_1 & = 0.02 \hspace{1em}\textrm{H}, \\
	g	& = 9.81 \hspace{1em} \textrm{m}/\textrm{sec}^2, &	R 	& = 1 \hspace{1em}\Omega, \\
	a 	& = 0.05 \hspace{1em}\textrm{m}, &	r 	& = 0.05 \hspace{1em}\textrm{m}.
\end{align*}

First, we solve $f(x,\hat{u})=0$ with $ x_1 =r$ in order to obtain the equilibrium point of interest of the system. This yields the solution
\[ x^0 = \begin{pmatrix} r \\ 0 \\  \sqrt{\frac{2gm(a+r)^2}{L_0a}}\end{pmatrix}, \quad\hat{u}^0 = R\sqrt{\frac{2gm(a+r)^2}{L_0a}}. \] 
We can shift the equilibrium point to the origin by defining $\bar{x}= x-x^0$, and $u=\hat{u}-\hat{u}^0$, obtaining in the new variables:
\[ \dot{\bar{x}} = f(\bar{x}+x_0,u+\hat{u}_0) \] 
Linearizing this around the origin yields $\dot{\bar{x}} = A\bar{x}+ Bu$, where
\[ A := \begin{bmatrix} 
	0 & 1 & 0 \\ 
	\frac{2g}{a+r} & -\frac{k}{m} & -\frac{L_0ax_3^0}{m(a+r)^2} \\
	0 & \frac{L_0ax_3^0}{L(r)(a+r)^2} & -\frac{R}{L(r)}
	\end{bmatrix}, 
\quad B := \begin{bmatrix} 0 \\ 0 \\ \frac{1}{L(r)}\end{bmatrix}. \]
Since we want to control the system on the basis of measurements of $\bar{x}_1$, we add an output $y=C\bar{x}$, where $C := \begin{bmatrix} 1 & 0&0\end{bmatrix}$. Now, we can discretize this with step size $\delta>0$ and obtain 
\[ \bar{x}(t+1) = (I+\delta A)\bar{x}(t) + \delta B u(t),\quad y(t)=C\bar{x}(t).\] 

In order to obtain an AR model, note that for each $s\geq 1$ we have that
\[y(t+s) = C(I+\delta A)^s \bar{x}(t) + \delta\sum_{i=0}^{s-1} C(I+\delta A)^i Bu(t+s-i-1).\]
The characteristic polynomial of $I+\delta A$ is denoted
\[ \chi (\lambda) = \lambda^3 + P_2\lambda^2 + P_1\lambda + P_0.\]
The Cayley-Hamilton theorem states that $\chi(I+\delta A) =0$. Therefore, we obtain that:
\begin{subequations}\label{eq:armamodel}
\begin{align}
	&y(t+3)+P_2y(t+2) +P_1 y(t+1) +P_0y(t) \\
	&= Q_2 u(t+2) + Q_1 u(t+1) + Q_0 u(t),
\end{align}
\end{subequations}
where the matrices $Q_0,Q_1,$ and $Q_2$ are given by:
\begin{align*}
	Q_0 & =	\delta C(I+\delta A)^2B + \delta P_2 C(I+\delta A) B +\delta P_1 CB, \\
	Q_1 & = \delta C(I+\delta A)B + \delta P_2 CB,	\\		
	Q_2 & = \delta  CB.
\end{align*}
This brings the model into the form considered in this paper. 

In this simulation example, we will perform measurements on the (discretized) nonlinear system. We will treat this nonlinear system as an AR model of the form \eqref{unknown true system} where the additive noise term $\bmv(t)$ captures the nonlinearities. Using the methods of this paper, we will find a stabilizing controller for all such systems consistent with the measurements and a noise model of the form $VV^\top \leq \epsilon$. 

We obtain measurements of the system close to the equilibrium point. To be precise we take $\delta=0.005$, $T=28$, and generate random inputs from the interval $10^{-5}[-1,1]$. These are applied to the nonlinear system with given initial conditions. The measurements resulting from this can be seen in \eqref{measurementsnonl}. For these measurements, we observe that $VV^\top\leq 10^{-17}$.

We will use Theorem~\ref{ch11:th:mainNS} to show that these measurements are informative for quadratic stabilization. For this, we first form the matrices $H'_1$, $H'_2$ and $\bar{N}$. It is straightforward to see that $H'_1$ has full row rank. We now use Yalmip with Mosek as a solver in order to find matrices $D \in \mathbb{R}^{1 \times 6 }$, and $\Phi \in \S{6}$, such that $\Phi > 0$ and the LMI \eqref{ch11:eq:mainLMI} holds. Indeed, such matrices exist, and therefore the data are informative for quadratic stabilization. We can find a stabilizing controller by taking $C= -D\Phi^{-1}$, which results in
\[C= \begin{bmatrix} -0.90 &35886.96&-0.94&-88735.92&-0.73&53966.58\end{bmatrix}\]
That is, a controller of the form: 
\begin{equation}\label{eq:controller}\small\begin{array}{lllllll}
	&\bmu(t+3) + 0.73\bmu(t+2) + 0.94\bmu(t+1)+0.90 \bmu(t) \\
	&= -53966.58 \bmy(t+2) +88735.92\bmy(t+1) +35886.96\bmy(t).\normalsize
\end{array}\end{equation}

\begin{figure}[t]
	\input{plots.tex}
	\caption{The results of interconnecting the controller \eqref{eq:controller} with the linearized system, two other systems consistent with the data, and the original nonlinear system.}\label{fig:plots}
\end{figure}
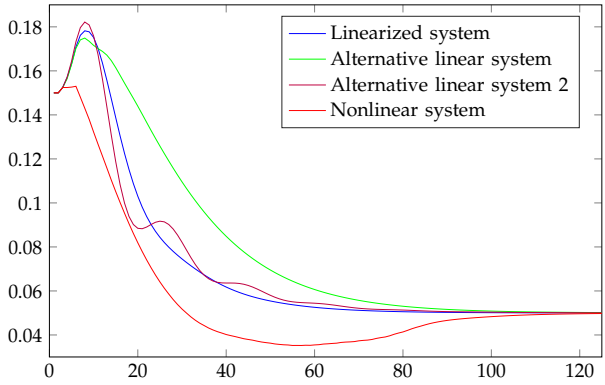
By definition, this means that the controller $C$ stabilizes all linear systems of the form \eqref{eq:armamodel} that are compatible with the measurements. A few trajectories of compatible systems interconnected with the controller are shown in Figure~\ref{fig:plots}. More specifically, the linearization derived earlier is consistent with the measurements, and is therefore stabilized by the found controller. As a last remark, we can interconnect the controller with the discrete-time nonlinear plant. Given that the controller stabilizes the linearization, it locally stabilizes the nonlinear system. This is also illustrated in Figure~\ref{fig:plots}.

\begin{figure*}
\hrulefill
\begin{equation}\label{measurementsnonl}
	\footnotesize\begin{array}{ll}		 
		\begin{bmatrix} y(0)\hphantom{11}\!\!\!\! &\cdots &\!\! y(15) \end{bmatrix}&\!=10^{-6}\left\lbrack\begin{array}{ccccccccccccccccc}0.100&0.100&0.100&0.101&0.101&0.103&0.105&0.107&0.109&0.112&0.115&0.118&0.122&0.126&0.130&0.135\end{array}\right\rbrack  \\
		\begin{bmatrix} y(16)\hphantom{0}\!\!\!\! &\cdots&\!\!  y(31) \end{bmatrix} 
		&\!=10^{-6} \left\lbrack\begin{array}{cccccccccccccccc}	0.141&0.48&0.156&0.164&0.174&0.184&0.195&0.207&0.220&0.233&0.247&0.262&0.277&0.294&0.313\end{array}\right\rbrack \\
		\begin{bmatrix} u(0)\hphantom{10}\!\!\!\! &\cdots&\!\!  u(14)\end{bmatrix} &\!=10^{-5}\left\lbrack\begin{array}{cccccccccccccccc}-0.75&\hphantom{-}0.61&\hphantom{-}0.71&-0.06&0.13&0.92&0.66&0.30&0.26&\hphantom{-}0.95&-0.86&-0.74&-0.88&-0.45&-0.11 \end{array}\right\rbrack \\		
		\begin{bmatrix} u(15)\hphantom{0}\!\!\!\! &\cdots&\!\! u(29)   \end{bmatrix} &\!=10^{-5}\left\lbrack\begin{array}{cccccccccccccccc}-0.59&-0.09&-0.40&\hphantom{-}0.85&0.95&0.57&0.49&0.67&0.51&-0.84&\hphantom{-}0.92&-0.57&-0.46&-0.31&\hphantom{-}0.86\end{array}\right\rbrack
	\end{array} \normalsize
\end{equation}
\hrulefill
\end{figure*}

\section{Conclusions and discussion} 
In this paper, we have given an introduction to the informativity approach to data-driven control. By using a combination of a new viewpoint, classical methods, and novel technical results, we have illustrated the framework by providing a number of solutions to problems with different model classes of linear systems, different types of measurements, and various control objectives. There remain, however, certain limitations to the results. While a number of these limitations yield interesting directions for future research, some of them are inherent to the approach.

First of all, we have been interested in providing necessary and sufficient conditions for informativity. Of course, such conditions are in a certain sense the gold standard, as they precisely characterize the information contained in the data. On the other hand, it might well be the case that the data contain more information than required. As such, a potentially interesting variant of these problems is to provide condition which are easier to check, but only sufficient. In many cases, such an approach might computationally outperform the methods of this paper. In a similar vein, a number of heuristic methods can drastically outperform the design methods of this paper, albeit without strong theoretical guarantees. In particular, in the case of very small noise samples, the set of consistent systems may be small. In this case, an intuitive method of performing data-driven control is a certainty-equivalent approach: Find any compatible system and solve the control problem for that system. 

The informativity approach has a number of moving parts: The control objective, model class, and noise model. In this paper, we have mainly varied our choice of control objective or analysis problem. In particular, we have focused entirely on model classes consisting of different flavors of linear systems. An extension towards nonlinear systems could improve the applicability of the results. For well-behaved nonlinear systems, we can draw certain conclusions on the basis of the behavior of its linearization. However, a more natural approach would be to investigate informativity problems for certain classes of nonlinear systems directly. Of course, when changing the model class one needs to balance the benefits of more general model classes and the tractability of the resulting robust control problems. Some classes of systems have shown a favorable trade-off in this regard, such as bilinear systems \cite{Bisoffi2020,Markovsky2022}, polynomial systems \cite{Guo2020,Guo2022}, rational systems \cite{Strasser2021} and systems with quadratic or sector bounded nonlinearities \cite{Luppi2022,vanWaarde2022}. As was shown in the aforementioned works, a thorough understanding of the linear case often remains invaluable for the proposal of nonlinear extensions. 

Another problem of interest is changing the way the noise acts on our system and measurements. Measurement noise, that is, noise which acts only on the measurements but not on the system, can be modeled in a similar manner as in this paper. However, an open problem is to provide conditions for data informativity in this setting, because the structure of the set of consistent systems appears to be more complicated than the ones studied here. So far, we are only aware of sufficient conditions for quadratic stabilization with measurement noise \cite[Sec. VA]{DePersis2020} that rely on somewhat conservative bounds. 

The noise models considered in this paper can be applied to treat different scenarios such as energy bounds and sample covariance bounds on the noise. An advantage of these noise models is that the resulting informativity conditions take the form of LMIs with a complexity that is independent of the number of measurements. Clearly, such limited computational complexity is desirable in any control problem. However, the assumption that the noise signal can be described by the solution set of a QMI also comes with certain limitations. 

In particular, with the noise models described in this paper, it is not possible to treat the situation of sample bounds without conservatism. More generally, combining different sets of measurements is a nontrivial problem in this setting, given that the intersection of solution sets of QMIs can generally not be described as the solution set of a single QMI. Without either developing tools that can deal with such intersections or alternative noise models, two important problems are difficult to tackle. First of all, the question of incremental informativity: Does adding more measurements lead to more informative data? Moreover, this lack of scalability inhibits the development of online or adaptive methods as compared to the offline methods of this paper. This motivates the development for new technical results for sample-bounded noise. In \cite[Sec. VII]{vanWaarde2022b} a simple sufficient LMI condition was proposed for quadratic stabilization in the presence of sample-bounded noise, which was further studied in \cite{Bisoffi2021e}. Although this approach appears to be less conservative than describing sample-bounded noise by a QMI, it is not well-understood from a theoretical perspective.  	

One of the strengths of methods based on the fundamental lemma \cite{Willems2005} (see also the Sidebar ``Willems' fundamental lemma") is the following: For controllable linear systems, we can \emph{guarantee} that the input-output data have favorable rank properties by injecting inputs that are persistently exciting. The fundamental lemma is thus an \emph{experiment design} result, that provides a guide for choosing the inputs of the experiment in order to generate informative data (for system identification). It can be shown that the persistency of excitation condition can be replaced by an \emph{online} design of the inputs \cite{vanWaarde2021}, which uses less data samples. An important topic for future work, however, is to develop experiment design methods corresponding to the various problems studied in this paper, especially those for noisy data. Although this is a largely unexplored area of research, we believe that the conditions provided in this paper will form the basis for an experiment design theory. Indeed, to be able to guarantee that the data are informative requires a thorough understanding of informativity in the first place.

\bibliographystyle{IEEEtran}
\bibliography{references}

\newpage

\end{document}

%% file: ka-newcommands.tex
\DeclareMathOperator{\im}{im}
\DeclareMathOperator{\col}{col}

\DeclareMathOperator{\rank}{rank}
\DeclareMathOperator{\trace}{tr}

\let\leq\leqslant
\let\geq\geqslant

\newcommand{\calB}{\ensuremath{\mathcal{B}}}

\newcommand{\calD}{\ensuremath{\mathcal{D}}}

\newcommand{\calH}{\ensuremath{\mathcal{H}}}

\newcommand{\calK}{\ensuremath{\mathcal{K}}}
\newcommand{\calL}{\ensuremath{\mathcal{L}}}
\newcommand{\calM}{\ensuremath{\mathcal{M}}}

\newcommand{\calO}{\ensuremath{\mathcal{O}}}
\newcommand{\calP}{\ensuremath{\mathcal{P}}}

\newcommand{\calS}{\ensuremath{\mathcal{S}}}

\newcommand{\calZ}{\ensuremath{\mathcal{Z}}}

\newcommand{\hatF}{\ensuremath{\hat{F}}}
\newcommand{\hatG}{\ensuremath{\hat{G}}}
\newcommand{\hatH}{\ensuremath{\hat{H}}}

\newcommand{\bmat}{\begin{matrix}}
\newcommand{\emat}{\end{matrix}}
\newcommand{\bbm}{\begin{bmatrix}}
\newcommand{\ebm}{\end{bmatrix}}
\newcommand{\bbma}{\begin{bmatrix*}[r]}
\newcommand{\ebma}{\end{bmatrix*}}
\newcommand{\bpm}{\begin{pmatrix}}
\newcommand{\epm}{\end{pmatrix}}
\newcommand{\bvm}{\begin{vmatrix}}
\newcommand{\evm}{\end{vmatrix}}
\newcommand{\bse}{\begin{subequations}}
\newcommand{\ese}{\end{subequations}}
\newcommand{\beq}{\begin{equation}}
\newcommand{\eeq}{\end{equation}}
\newcommand{\ben}{\renewcommand{\labelenumi}{\arabic{enumi}.}
\renewcommand{\theenumi}{\arabic{enumi}}\begin{enumerate}}

\newcommand{\een}{\end{enumerate}}

\newcommand{\beni}{\renewcommand{\labelenumi}{\roman{enumi}.}
\renewcommand{\theenumi}{\roman{enumi}}\begin{enumerate}}

\newcommand{\eeni}{\end{enumerate}}

\newcommand{\bena}{\renewcommand{\labelenumi}{\alph{enumi}.}
\renewcommand{\theenumi}{\alph{enumi}}\begin{enumerate}}

\newcommand{\eena}{\end{enumerate}}

\newcommand{\bit}{\begin{itemize}}
\newcommand{\eit}{\end{itemize}}
\newcommand{\bthe}{\begin{theorem}}
\newcommand{\ethe}{\end{theorem}}
\newcommand{\blem}{\begin{lemma}}
\newcommand{\elem}{\end{lemma}}
\newcommand{\bprop}{\begin{proposition}}
\newcommand{\eprop}{\end{proposition}}
\newcommand{\bex}{\begin{example}}
\newcommand{\eex}{\end{example}}
\newcommand{\bas}{\begin{assumption}}
\newcommand{\eas}{\end{assumption}}
\newcommand{\bre}{\begin{remark}}
\newcommand{\ere}{\end{remark}}
\newcommand{\bcor}{\begin{corollary}}
\newcommand{\ecor}{\end{corollary}}
\newcommand{\bdfn}{\begin{definition}}
\newcommand{\edfn}{\end{definition}}
\newcommand{\bcon}{\begin{conjecture}}
\newcommand{\econ}{\end{conjecture}}

\newcommand{\half}{\ensuremath{\frac{1}{2}}}
\newcommand{\inv}{\ensuremath{^{-1}}}
\newcommand{\pset}[1]{\ensuremath{\{#1\}}}

\newcommand{\set}[2]{\ensuremath{\{#1\mid #2\}}}

\newcommand{\abs}[1]{\ensuremath{| #1 |}}
\newcommand{\norm}[1]{\ensuremath{\| #1 \|}}

\newcommand{\R}{\ensuremath{\mathbb R}}
\newcommand{\C}{\ensuremath{\mathbb C}}

\newcommand{\Z}{\ensuremath{\mathbb Z}}



%% file: plots.tex
\centering
\tikzset{every picture/.style={scale=0.82}}
\subfigure{
	\label{plot:position}\begin{tikzpicture}
		\begin{axis}[
			width=4.3in,
			height=3.566in/1.3,
			at={(0.758in,0.481in)},
			scale only axis,
			xmin=0,
			xmax=125,
			ymin=0.03,
			ymax=0.19,
			yticklabel style={/pgf/number format/fixed,/pgf/number format/precision=5},
			scaled y ticks=false,
			axis background/.style={fill=white},
			legend style={at={(axis cs:120,0.185)}, legend cell align=left, align=left, draw=white!15!black}
			]
			
			\addplot[color=blue] table[row sep=crcr] {%
1	0.150000000000000 \\
2	0.150000000000000 \\
3	0.152437200000000 \\
4	0.156922138140000 \\
5	0.163152803605893 \\
6	0.170898008587214 \\
7	0.176174022995042 \\
8	0.178196352465183 \\
9	0.177749996078542 \\
10	0.175188275838864 \\
11	0.170314418835600 \\
12	0.163559172886239 \\
13	0.155721946780477 \\
14	0.147360136666056 \\
15	0.138797597189297 \\
16	0.130388213437168 \\
17	0.122471654324185 \\
18	0.115237950989931 \\
19	0.108752908759253 \\
20	0.103040689587828 \\
21	0.0980915073334777\\
22	0.0938444967042829\\
23	0.0902068067011053\\
24	0.0870825856886305\\
25	0.0843806056596108\\
26	0.0820132290065666\\
27	0.0799019730164358\\
28	0.0779837114458279\\
29	0.0762105763290801\\
30	0.0745471837928820\\
31	0.0729692474958623\\
32	0.0714624035251134\\
33	0.0700197895697154\\
34	0.0686393924503186\\
35	0.0673221906803033\\
36	0.0660707955347108\\
37	0.0648881998564629\\
38	0.0637768188754898\\
39	0.0627379938792419\\
40	0.0617718082043912\\
41	0.0608770509784998\\
42	0.0600513135852438\\
43	0.0592912117459680\\
44	0.0585926613573662\\
45	0.0579511457469385\\
46	0.0573619588120677\\
47	0.0568204186301974\\
48	0.0563220367951826\\
49	0.0558626354363524\\
50	0.0554384175795795\\
51	0.0550459997381670\\
52	0.0546824130785910\\
53	0.0543450800989288\\
54	0.0540317755438186\\
55	0.0537405793197386\\
56	0.0534698268518257\\
57	0.0532180609275852\\
58	0.0529839882026203\\
59	0.0527664422918491\\
60	0.0525643540846839\\
61	0.0523767291330186\\
62	0.0522026315255802\\
63	0.0520411733153294\\
64	0.0518915083451155\\
65	0.0517528293015059\\
66	0.0516243669420740\\
67	0.0515053905910004\\
68	0.0513952091670755\\
69	0.0512931722004299\\
70	0.0511986704820318\\
71	0.0511111361466875\\
72	0.0510300421137083\\
73	0.0509549009044118\\
74	0.0508852629213496\\
75	0.0508207143107466\\
76	0.0507608745431291\\
77	0.0507053938445482\\
78	0.0506539505969385\\
79	0.0506062488048827\\
80	0.0505620157016150\\
81	0.0505209995429261\\
82	0.0504829676157406\\
83	0.0504477044696394\\
84	0.0504150103651880\\
85	0.0503846999227857\\
86	0.0503566009495087\\
87	0.0503305534184940\\
88	0.0503064085751033\\
89	0.0502840281457439\\
90	0.0502632836281052\\
91	0.0502440556451025\\
92	0.0502262333485325\\
93	0.0502097138619922\\
94	0.0501944017557473\\
95	0.0501802085488372\\
96	0.0501670522357096\\
97	0.0501548568361136\\
98	0.0501435519679033\\
99	0.0501330724428954\\
100	0.0501233578860900\\
101	0.0501143523784880\\
102	0.0501060041235176\\
103	0.0500982651367862\\
104	0.0500910909585631\\
105	0.0500844403881059\\
106	0.0500782752387113\\
107	0.0500725601121945\\
108	0.0500672621914006\\
109	0.0500623510493091\\
110	0.0500577984733084\\
111	0.0500535783032778\\
112	0.0500496662822010\\
113	0.0500460399181485\\
114	0.0500426783565812\\
115	0.0500395622620481\\
116	0.0500366737084631\\
117	0.0500339960772486\\
118	0.0500315139627228\\
119	0.0500292130841847\\
120	0.0500270802042145\\
121	0.0500251030527562\\
122	0.0500232702565924\\
123	0.0500215712738532\\
124	0.0500199963332259\\
125	0.0500185363775538\\};'
			\addlegendentry{Linearized system}
			\addplot[color=green] table[row sep=crcr] {%
1	0.150000000000000 \\
2	0.150000000000000 \\
3	0.152437200000000 \\
4	0.156397066680849 \\
5	0.162830098621997 \\
6	0.170226323074362 \\
7	0.174094396343809 \\
8	0.174901546045952 \\
9	0.173447266890721 \\
10	0.171585585494951 \\
11	0.170133425254568 \\
12	0.168843406289300 \\
13	0.167154232151877 \\
14	0.164667188215226 \\
15	0.161488232363287 \\
16	0.157970955476278 \\
17	0.154427004512901 \\
18	0.150959850274299 \\
19	0.147501883251999 \\
20	0.143958387085635 \\
21	0.140305632652485 \\
22	0.136600563595762 \\
23	0.132925919670788 \\
24	0.129336265031060 \\
25	0.125842490082376 \\
26	0.122430703071438 \\
27	0.119089486428568 \\
28	0.115823181990590 \\
29	0.112647648779623 \\
30	0.109578421120182 \\
31	0.106622679799877 \\
32	0.103779327314834 \\
33	0.101044078494062 \\
34	0.0984141667110067\\
35	0.0958896124604241\\
36	0.0934715885217248\\
37	0.0911602304558894\\
38	0.0889537068911720\\
39	0.0868487684440069\\
40	0.0848418655929304\\
41	0.0829298748222162\\
42	0.0811101144596046\\
43	0.0793799479885364\\
44	0.0777364412990985\\
45	0.0761763171991114\\
46	0.0746961441954605\\
47	0.0732925538688772\\
48	0.0719623384057845\\
49	0.0707024191303017\\
50	0.0695097687073794\\
51	0.0683813681638432\\
52	0.0673142223342550\\
53	0.0663054058431021\\
54	0.0653520999789570\\
55	0.0644516015612140\\
56	0.0636013109587373\\
57	0.0627987174568597\\
58	0.0620413942802140\\
59	0.0613270034571992\\
60	0.0606533032598832\\
61	0.0600181515642980\\
62	0.0594195035777203\\
63	0.0588554066439023\\
64	0.0583239955956351\\
65	0.0578234902020356\\
66	0.0573521940383003\\
67	0.0569084932537183\\
68	0.0564908542898905\\
69	0.0560978206424721\\
70	0.0557280093524976\\
71	0.0553801078136083\\
72	0.0550528710256608\\
73	0.0547451190685926\\
74	0.0544557345221207\\
75	0.0541836597289376\\
76	0.0539278939836748\\
77	0.0536874907914300\\
78	0.0534615552816197\\
79	0.0532492417713792\\
80	0.0530497514239014\\
81	0.0528623299577610\\
82	0.0526862654021247\\
83	0.0525208859210170\\
84	0.0523655577313008\\
85	0.0522196831226740\\
86	0.0520826985718094\\
87	0.0519540729377259\\
88	0.0518333057305140\\
89	0.0517199254529990\\
90	0.0516134880184693\\
91	0.0515135752461272\\
92	0.0514197934323153\\
93	0.0513317719931012\\
94	0.0512491621737055\\
95	0.0511716358216348\\
96	0.0510988842216807\\
97	0.0510306169912745\\
98	0.0509665610342025\\
99	0.0509064595500558\\
100	0.0508500710965382\\
101	0.0507971687019371\\
102	0.0507475390254261\\
103	0.0507009815631189\\
104	0.0506573078978539\\
105	0.0506163409906299\\
106	0.0505779145115790\\
107	0.0505418722084105\\
108	0.0505080673103834\\
109	0.0504763619659958\\
110	0.0504466267126742\\
111	0.0504187399768101\\
112	0.0503925876025324\\
113	0.0503680624076662\\
114	0.0503450637653946\\
115	0.0503234972102243\\
116	0.0503032740669268\\
117	0.0502843111011979\\
118	0.0502665301908297\\
119	0.0502498580162505\\
120	0.0502342257693400\\
121	0.0502195688794862\\
122	0.0502058267559054\\
123	0.0501929425452979\\
124	0.0501808629039632\\
125	0.0501695377835429\\};'
\addlegendentry{Alternative linear system}
			\addplot[color=purple] table[row sep=crcr] {%
1	0.150000000000000 \\
2	0.150000000000000 \\
3	0.152437200000000 \\
4	0.157104639078601 \\
5	0.164287451114761 \\
6	0.173248650601364 \\
7	0.179656327252243 \\
8	0.182192504970362 \\
9	0.180835572246098 \\
10	0.175871626614687 \\
11	0.167405744121693 \\
12	0.156262499139148 \\
13	0.143670799665630 \\
14	0.130847695800404 \\
15	0.118775667434878 \\
16	0.108264480904641 \\
17	0.0998811801790493\\
18	0.0938732713277902\\
19	0.0901692580846508\\
20	0.0884588646871357\\
21	0.0882695769920548\\
22	0.0890387929753967\\
23	0.0901906394982656\\
24	0.0912106659573027\\
25	0.0916990303124562\\
26	0.0913989545259820\\
27	0.0902042564444340\\
28	0.0881490608388554\\
29	0.0853817783372526\\
30	0.0821286370038855\\
31	0.0786534378143440\\
32	0.0752191985886657\\
33	0.0720557248232643\\
34	0.0693360494938311\\
35	0.0671633717262656\\
36	0.0655685356255345\\
37	0.0645166995985075\\
38	0.0639209688705458\\
39	0.0636603017803173\\
40	0.0635988813016974\\
41	0.0636043705309832\\
42	0.0635630012764080\\
43	0.0633901487423850\\
44	0.0630357959883619\\
45	0.0624849977412111\\
46	0.0617540430208929\\
47	0.0608834331803881\\
48	0.0599290087818000\\
49	0.0589525802449730\\
50	0.0580132719392186\\
51	0.0571605198765937\\
52	0.0564293202605400\\
53	0.0558379626509627\\
54	0.0553881449096504\\
55	0.0550670940758921\\
56	0.0548511312154502\\
57	0.0547100284763160\\
58	0.0546115094150069\\
59	0.0545253247370696\\
60	0.0544264732069463\\
61	0.0542973063714900\\
62	0.0541284307039698\\
63	0.0539184796828745\\
64	0.0536729543394809\\
65	0.0534024136501564\\
66	0.0531203323207442\\
67	0.0528409357035817\\
68	0.0525772774166538\\
69	0.0523397555671875\\
70	0.0521351806300039\\
71	0.0519664240367958\\
72	0.0518326017509723\\
73	0.0517296892018789\\
74	0.0516514273912586\\
75	0.0515903659600098\\
76	0.0515388957976008\\
77	0.0514901474095844\\
78	0.0514386663125614\\
79	0.0513808172198020\\
80	0.0513149089811047\\
81	0.0512410673188489\\
82	0.0511609088830896\\
83	0.0510770861256274\\
84	0.0509927776222491\\
85	0.0509111937768740\\
86	0.0508351553800865\\
87	0.0507667849555619\\
88	0.0507073311072180\\
89	0.0506571268988994\\
90	0.0506156669003276\\
91	0.0505817754837356\\
92	0.0505538320562768\\
93	0.0505300172194260\\
94	0.0505085467742168\\
95	0.0504878669846813\\
96	0.0504667932434430\\
97	0.0504445838553055\\
98	0.0504209497791741\\
99	0.0503960088075311\\
100	0.0503701981152641\\
101	0.0503441620464719\\
102	0.0503186324544215\\
103	0.0502943171895811\\
104	0.0502718089806859\\
105	0.0502515226273156\\
106	0.0502336637919954\\
107	0.0502182283527801\\
108	0.0502050277321676\\
109	0.0501937331661878\\
110	0.0501839306501192\\
111	0.0501751782575983\\
112	0.0501670585026357\\
113	0.0501592201278130\\
114	0.0501514058366105\\
115	0.0501434647208287\\
116	0.0501353501798262\\
117	0.0501271057672146\\
118	0.0501188424957514\\
119	0.0501107116321731\\
120	0.0501028769499639\\
121	0.0500954898729270\\
122	0.0500886700722560\\
123	0.0500824930308874\\
124	0.0500769850154963\\
125	0.0500721249324230\\};'
\addlegendentry{Alternative linear system 2}			
			\addplot[color=red] table[row sep=crcr] {%
1	0.150000000000000 \\
2	0.150000000000000 \\
3	0.152437200000000 \\
4	0.152437200000000 \\
5	0.152622604930932 \\
6	0.152993405522550 \\
7	0.148055508596784 \\
8	0.143036330517787 \\
9	0.137930568782675 \\
10	0.132196582040095 \\
11	0.126684520502193 \\
12	0.121297842042165 \\
13	0.115869622516075 \\
14	0.110403100654353 \\
15	0.105169891469695 \\
16	0.100163609303440 \\
17	0.0953037287682851\\
18	0.0905903973477975\\
19	0.0861225511445921\\
20	0.0818977578067498\\
21	0.0778995105003008\\
22	0.0741123863493187\\
23	0.0705630025501969\\
24	0.0672178348888205\\
25	0.0641179155087522\\
26	0.0612456126245331\\
27	0.0586039014928192\\
28	0.0560981446765926\\
29	0.0538326172042824\\
30	0.0517886158354041\\
31	0.0499844159598622\\
32	0.0482358614940641\\
33	0.0467153845465923\\
34	0.0453829586209279\\
35	0.0442894604069313\\
36	0.0431827237627374\\
37	0.0422886113456243\\
38	0.0415119178383263\\
39	0.0408742803082074\\
40	0.0402002442711850\\
41	0.0397171882841609\\
42	0.0392739250162930\\
43	0.0387834860815341\\
44	0.0382892249093285\\
45	0.0379524482969524\\
46	0.0376284436153376\\
47	0.0371855660921496\\
48	0.0367995638558818\\
49	0.0365199672845077\\
50	0.0362811148998878\\
51	0.0359272816958627\\
52	0.0357147799715228\\
53	0.0355396941696385\\
54	0.0354481508717610\\
55	0.0352314397127628\\
56	0.0352445675908893\\
57	0.0352250220802923\\
58	0.0353317352912974\\
59	0.0352935084523967\\
60	0.0354724222179886\\
61	0.0355876653532488\\
62	0.0358685135446094\\
63	0.0359850762579762\\
64	0.0361444282231884\\
65	0.0362767467430077\\
66	0.0366029669536525\\
67	0.0367617535595984\\
68	0.0368499857392931\\
69	0.0369892868373839\\
70	0.0373348564047146\\
71	0.0375366402987575\\
72	0.0376705915184115\\
73	0.0379550849347087\\
74	0.0384387253506689\\
75	0.0387996869599990\\
76	0.0390921974173296\\
77	0.0396224494669650\\
78	0.0403093954326743\\
79	0.0408780663472082\\
80	0.0413695951903007\\
81	0.0420614304746792\\
82	0.0428102552834023\\
83	0.0434446584164671\\
84	0.0439990633001417\\
85	0.0445853965789272\\
86	0.0451268713286700\\
87	0.0455795993172644\\
88	0.0459689688393570\\
89	0.0463156930470467\\
90	0.0466133908899579\\
91	0.0468673128543515\\
92	0.0470900030032973\\
93	0.0472892802387404\\
94	0.0474691127789141\\
95	0.0476338457700265\\
96	0.0477875306684476\\
97	0.0479325050403540\\
98	0.0480698914646275\\
99	0.0482004926864680\\
100	0.0483248292420085\\
101	0.0484429936754278\\
102	0.0485548538370173\\
103	0.0486602971668193\\
104	0.0487592669920500\\
105	0.0488517410232050\\
106	0.0489377706671058\\
107	0.0490175158802167\\
108	0.0490912279831590\\
109	0.0491592192034654\\
110	0.0492218472365728\\
111	0.0492795016243772\\
112	0.0493325829453570\\
113	0.0493814840146823\\
114	0.0494265789316139\\
115	0.0494682162914222\\
116	0.0495067138391913\\
117	0.0495423558627593\\
118	0.0495753938416198\\
119	0.0496060488272484\\
120	0.0496345144111921\\
121	0.0496609601908587\\
122	0.0496855355800816\\
123	0.0497083734547094\\
124	0.0497295933095655\\
125	0.0497493039170189\\};
			\addlegendentry{Nonlinear system}
		\end{axis}
	\end{tikzpicture}
}